\documentclass{amsart}
\usepackage[hmargin=2.1cm,vmargin=2cm]{geometry}
\usepackage{adjustbox}
\usepackage{graphicx}
\usepackage[utf8]{inputenc}
\usepackage{tikz}
\usepackage{pgfplots}
\usetikzlibrary{backgrounds,calc,fit,intersections}
\usepgfplotslibrary{fillbetween}
\pgfplotsset{compat=1.18}
\usepackage{amsthm}
\usepackage{amsmath}

\usepackage{xcolor}

\newcommand{\compactfiginput}[1]{\adjustbox{max width=0.7\linewidth,max height=5.4cm}{%
    \begingroup\setlength{\textwidth}{\linewidth}\input{#1}\endgroup}}
\newcommand{\compactsubfiginput}[1]{\adjustbox{max width=0.85\linewidth,max height=4.4cm}{%
    \begingroup\setlength{\textwidth}{\linewidth}\input{#1}\endgroup}}
\newcommand{\typepanel}[1]{\adjustbox{scale=0.66}{%
    \begingroup\setlength{\textwidth}{\linewidth}\input{#1}\endgroup}}

\tikzset{
    xingHalo/.style={},  
    Spath/.style   ={draw=black,          line width=0.9pt, line cap=round, line join=round, xingHalo},
    Tpath/.style   ={draw=green!55!black, line width=0.9pt, line cap=round, line join=round, xingHalo},
    Hpath/.style   ={draw=red!70!black,   line width=0.9pt, line cap=round, line join=round, xingHalo},
    Upath/.style   ={draw=blue!60!black,  line width=0.9pt, line cap=round, line join=round, xingHalo},
    rootedge/.style={draw=orange!90!black, line width=0.7pt, line cap=round, xingHalo},
}

\definecolor{bsOne}{RGB}{31,119,180}     
\definecolor{bsTwo}{RGB}{44,160,44}      
\definecolor{bsThree}{RGB}{214,39,40}    
\definecolor{bsFour}{RGB}{255,127,14}    
\definecolor{bsFive}{RGB}{148,103,189}   
\definecolor{bsSix}{RGB}{140,86,75}      

\tikzset{
    bsHalo/.style={line width=7pt, opacity=0.32,
                   line cap=round, line join=round},
    setR1/.style={bsHalo, draw=bsThree},
    setR2/.style={bsHalo, draw=bsFour},
    setR3/.style={bsHalo, draw=bsFive},
    setR4/.style={bsHalo, draw=bsSix},
}

\newtheorem{theorem}{Theorem}
\newtheorem{lemma}[theorem]{Lemma}

\newtheorem{remark}[theorem]{Remark}

\theoremstyle{definition}
\newtheorem{definition}[theorem]{Definition}
\newtheorem{conjecture}[theorem]{Conjecture}
\newtheorem*{theorem*}{Theorem}
\newenvironment{restated}[1]{\par\medskip\noindent\textbf{Theorem~#1.}\itshape\ }{\par\medskip}

\title{A Rooted Structure Theorem for $5$-Connected Graphs and $K_6^-$-Minors}

\author{Ken-ichi Kawarabayashi}
\thanks{K.~Kawarabayashi: National Institute of Informatics \& The University of Tokyo, Tokyo, Japan, \texttt{k\_keniti@nii.ac.jp}. Supported by the Alexander von Humboldt Foundation (Humboldt Research Award), JSPS Kakenhi 26K21777 and JP25K24465, and JST ASPIRE JPMJAP2302.}

\author{Alexandre Simon}
\thanks{A.~Simon: National Institute of Informatics, Tokyo, Japan, and The Graduate University for Advanced Studies (SOKENDAI), Hayama, Japan, \texttt{alexandresimon@nii.ac.jp}. This work was supported by JST, the establishment of university fellowships toward the creation of science technology innovation, Grant Number JPMJFS2136.}


\begin{document}

\begin{abstract}
J{\o}rgensen conjectured that every $6$-connected graph with no $K_6$-minor is apex. For a $6$-connected graph containing a $K_4$ on vertices $a_1,\dots,a_4$, being apex is equivalent to $G-a_i$ having a planar embedding with the other three roots on a common face for some $i$, so that for such graphs the conjecture becomes a rooted planarity statement. Robertson, Seymour, and Thomas showed that the failure of the corresponding three-rooted planarity is witnessed by a legless tripod~\cite{RST1993}. We introduce a four-root structure, called a \emph{frame}, that extends this certificate, and prove a structural theorem: every $5$-connected graph containing a frame rooted at four prescribed vertices contains a subgraph of one of five explicit types, denoted $A$--$E$. This theorem does not depend on the edges among the roots. As an application, we show that when the four roots span a $K_4^-$, each of the five configurations yields a $K_6^-$-minor. Together with a planar argument for the case without a frame and an induction based on a theorem of Kawarabayashi~\cite{Kawarabayashi2001}, this gives a proof that every non-planar $5$-connected graph contains a $K_6^-$-minor. This statement also follows from a recent theorem of Lo~\cite{Lo2026}; our proof is independent and is driven by the rooted structure theorem, whose method we expect to apply to J{\o}rgensen's conjecture with other target configurations.
\end{abstract}

\maketitle

\section{Introduction}

Hadwiger's conjecture predicts that every graph with no $K_t$-minor is $(t-1)$-colorable~\cite{Hadwiger1943}. Robertson, Seymour, and Thomas proved the case $t=6$~\cite{RST1993} by showing that a minimal counterexample is apex, that is, becomes planar after deleting one vertex. This motivated J{\o}rgensen~\cite{Jorgensen1994} to conjecture that apex graphs are the only obstruction to $K_6$-minors in highly connected graphs.

\begin{conjecture}[J{\o}rgensen]\label{conj:jorgensen}
Every $6$-connected graph with no $K_6$-minor is apex.
\end{conjecture}

Kawarabayashi, Norine, Thomas, and Wollan~\cite{KNTW2018} proved Conjecture~\ref{conj:jorgensen} for all sufficiently large graphs; the general case is open.

The starting point of our approach is that, in the presence of a $K_4$, the conjecture can be stated in terms of four prescribed roots. Let $G$ be $6$-connected and suppose that four vertices $a_1,a_2,a_3,a_4$ span a $K_4$. If $G-v$ is planar, then $G-v$ is $5$-connected, so every triangle of $G-v$ bounds a face; a plane graph containing a $K_4$ whose four triangles all bound faces is $K_4$ itself, so $v=a_i$ for some $i$, and the other three roots, which span a triangle, lie on a common face of $G-a_i$. Hence, for such $G$, Conjecture~\ref{conj:jorgensen} is equivalent to the following rooted statement: \emph{either $G-a_i$ has a planar embedding with the other three roots on the outer face for some $i$, or $G$ has a $K_6$-minor.} This form of planarity does not depend on the edges among the roots, and its failure has a structural certificate: by a theorem of Robertson, Seymour, and Thomas (Lemma~\ref{Seymour}), if $G-a_4$ has no planar embedding with $a_1,a_2,a_3$ on the outer face, then $G-a_4$ contains a legless tripod on $a_1,a_2,a_3$.

Our first step is to extend this certificate from three roots to four. A \emph{frame} rooted at $\{a_1,a_2,a_3,a_4\}$ consists of two trees, disjoint except at the roots, one reaching all four roots and the other reaching three of them (Definition~\ref{def:frame}); it exists as soon as rooted planarity fails at $a_3$ or at $a_4$ (Lemma~\ref{lem:legless-tripod-to-frame}). Our main result describes what a frame forces in a $5$-connected graph.

\begin{restated}{\ref{structure} (Frame structure theorem)}
If a $5$-connected graph contains a frame rooted at $\{a_1,a_2,a_3,a_4\}$, then it contains a subgraph of one of five explicit types $A$--$E$ rooted at these vertices.
\end{restated}

The theorem is purely structural: it uses neither the edges among the roots nor any assumption on the minors contained in the graph, and it shows that one of five explicit configurations always occurs. In the setting of Conjecture~\ref{conj:jorgensen}, type $A$ is already decisive. A configuration of type $A$ consists of two trees, each reaching all four roots and disjoint except at the roots; since $\{a_1,a_2,a_3,a_4\}$ is not a separator, they are joined by a path avoiding the roots, and together with a $K_4$ on the roots this yields a $K_6$-minor. The other types are tailored to the $K_4^-$ case of Theorem~\ref{thm:k6minus}; for $K_6$, the same frame analysis can be directed toward a different list of terminal configurations (see the concluding remarks).

To test the method, we carry it out one step below $K_6$: we replace the $K_4$ by a $K_4^-$, $6$-connectivity by $5$-connectivity, and the $K_6$-minor by a $K_6^-$-minor, where $K_6^-$ is $K_6$ minus an edge.

\begin{restated}{\ref{thm:k6minus}}
Every non-planar $5$-connected graph contains a $K_6^-$-minor.
\end{restated}

The proof follows the rooted reformulation above. Suppose first that $G$ contains a $K_4^-$, on $a_1,a_2,a_3,a_4$ with $a_3a_4$ missing. If rooted planarity fails at $a_3$ or at $a_4$, then $G$ contains a frame, and each of the five configurations of Theorem~\ref{structure} yields a $K_6^-$-minor (Lemma~\ref{K4minusgivesK6minus}). Otherwise both $G-a_3$ and $G-a_4$ are planar with the other three roots on the outer face, and we exclude this case by an apex-pair argument and a planar counting argument; when the four vertices span a $K_4$, we first choose the pair playing the role of the missing edge (Lemma~\ref{lem:k4-pair}). Finally, if $G$ contains no $K_4^-$, a theorem of Kawarabayashi~\cite{Kawarabayashi2001} gives an edge whose contraction preserves $5$-connectivity, and we show that such a contraction cannot make $G$ planar unless $G$ already has a $K_6^-$-minor; this closes an induction.

Theorem~\ref{thm:k6minus} is related to the Kelmans--Seymour theorem that every $5$-connected non-planar graph contains a subdivision of $K_5$~\cite{HeWangYu2020}; here we work with minors instead of subdivisions. Lo~\cite{Lo2026} recently obtained Theorem~\ref{thm:k6minus} in a stronger form, for $4$-connected non-planar graphs of minimum degree at least five and order at least seven (a $6$-vertex graph of minimum degree five is $K_6$), as a consequence of a characterization of graphs with no $K_{3,4}$-minor. Our proof is independent and gives rooted information that an excluded-minor characterization does not: it produces, through the frame structure theorem, explicit configurations located at four prescribed roots, which is the tool needed for Conjecture~\ref{conj:jorgensen}. Theorem~\ref{thm:k6minus} is the frame approach carried out one step below $K_6$.

Theorem~\ref{structure} is close in spirit to a result of Hayashi, Kawarabayashi, and Yoo~\cite{HKY2025}, who also start from a tripod on prescribed terminals and use it to construct rooted configurations. Under their connectivity and separation hypotheses, if $Z$ is a stable set of four vertices and $G$ cannot be drawn in a disk with $Z$ on the boundary, then, for some ordering $Z=\{v_1,v_2,v_3,v_4\}$, $G$ contains a nondegenerate diamond on $(v_1,v_2;v_3,v_4)$, a $K_{2,3}^+$-subdivision with feet $v_1,v_2,v_3$ avoiding $v_4$, or a $K_5^-$-subdivision with feet $v_1,v_2$~\cite[Theorems~1.2 and~1.4]{HKY2025}. Their motivation is Haj\'os' conjecture for $K_5$, and their alternatives are tailored to that problem. Our setting differs in three respects: we assume $5$-connectivity, impose no condition on the edges among the roots, and tailor our five configurations to minors, so that the edges among the roots turn them into $K_6^-$-minors, or into $K_6$-minors in the case of type~$A$ when the roots span a $K_4$.

The paper is organized as follows. Section~\ref{sec:prelim} fixes the terminology and the forms of Menger's theorem and Perfect's theorem that we use. Section~\ref{sec:structure} introduces frames and proves the frame structure theorem. Section~\ref{sec:k6minus} derives Theorem~\ref{thm:k6minus}: the framed $K_4^-$ case, the case without a frame, and the $K_4^-$-free case. We conclude with the open step toward Conjecture~\ref{conj:jorgensen}.

\section{Preliminaries}\label{sec:prelim}

We use standard graph-theoretic terminology; see Diestel~\cite{Diestel}. All graphs are finite, simple, and undirected. For a graph $G$, we write $V(G)$ and $E(G)$ for its vertex and edge sets. Since all graphs are simple, we write $uv$ for the edge joining two adjacent vertices $u$ and $v$.

If $v \in V(G)$, a vertex $u$ is a \emph{neighbor} of $v$ if $uv \in E(G)$. The set of all neighbors of $v$ in $G$ is denoted by $N_G(v)$, or simply $N(v)$ when the graph is clear from context.

A graph is \emph{planar} if it can be drawn in the plane so that edges meet only at common endpoints. A \emph{planar embedding} is such a drawing with a fixed cyclic order of the edges around each vertex.

A graph is \emph{$k$-connected} if it has more than $k$ vertices and remains connected after deleting any set of fewer than $k$ vertices. A set $X \subseteq V(G)$ is a \emph{$k$-separator} if $|X| = k$ and $G - X$ is disconnected. Thus a $k$-connected graph contains no $\ell$-separator with $\ell < k$; we use this fact throughout in the form ``$X$ is not a separator'' (or ``$X$ is not an $\ell$-separator'') for a set $X$ with $|X| = \ell < k$.

A graph $F$ is a \emph{minor} of $G$ if it can be obtained from $G$ by deleting vertices, deleting edges, and contracting edges. Equivalently, $F$ can be obtained from a subgraph of $G$ by contracting edges. All contractions are followed by deleting loops and parallel edges. If $e = uv$ is an edge of $G$, then $G/e$, also written $G/uv$, denotes the graph obtained by contracting $e$. We write $K_r^-$ for the graph obtained from $K_r$ by deleting one edge. A graph is \emph{$K_4^-$-free} if it has no subgraph isomorphic to $K_4^-$, induced or not; in particular, a $K_4^-$-free graph contains no $K_4$. Accordingly, a $K_4^-$ of $G$ is a subgraph of $G$ isomorphic to $K_4^-$.

If $H$ is a subgraph of $G$, an \emph{$H$-path} is a path in $G$ whose ends lie in $V(H)$ and whose internal vertices lie in $V(G)\setminus V(H)$. Two paths are \emph{internally disjoint} if they have no common internal vertex. When several $H$-paths are considered together, they are always internally disjoint from $H$; any additional disjointness between the paths themselves is stated explicitly.

For a path $Q$ and vertices $x,y\in V(Q)$, we write $Q[x,y]$ for the $(x,y)$-subpath of $Q$. The notations $Q(x,y)$, $Q[x,y)$, and $Q(x,y]$ have their usual meanings, obtained by deleting the indicated endpoints.

For a cycle $C$ and distinct vertices $x,y\in V(C)$, whenever one of the two $(x,y)$-arcs of $C$ is specified by context, we denote this arc by $C[x,y]$. The notations $C(x,y)$, $C[x,y)$, and $C(x,y]$ are obtained by deleting the indicated endpoints.

For $x \in V(G)$ and $B \subseteq V(G)\setminus\{x\}$, an \emph{$(x,B)$-fan of size $k$} is a set of $k$ paths from $x$ to $B$ that pairwise intersect only in $x$, whose ends in $B$ are pairwise distinct, and whose internal vertices avoid $B$.

We use the following fan form of Menger's theorem.

\begin{lemma}[Menger's theorem, fan form]\label{lem:menger-fan}
    Let $G$ be a $k$-connected graph, let $x \in V(G)$, and let $B \subseteq V(G)\setminus\{x\}$ with $|B|\ge k$.
    Then $G$ contains an $(x,B)$-fan of size $k$.
\end{lemma}

We also use the following augmentation form of Perfect's theorem; we include a short proof for completeness.

\begin{lemma}[Perfect's theorem~\cite{Perfect1968}]\label{perfect}
    Let $G$ be a graph, let $x \in V(G)$, and let $B \subseteq V(G)\setminus\{x\}$.
    Let $1 \le j \le k$. If $G$ contains an $(x,B)$-fan of size $k$, and an
    $(x,B)$-fan of size $j$ ending at distinct vertices $b_1,\dots,b_j \in B$,
    then $G$ contains an $(x,B)$-fan of size $k$ in which $j$ of the paths end
    at $b_1,\dots,b_j$.
\end{lemma}

\begin{proof}
    If $j=k$, the given fan of size $j$ is the required fan. Otherwise, it suffices to treat the case $k=j+1$; indeed, a subset of an $(x,B)$-fan is again an $(x,B)$-fan, so the case $k=j+1$ applied successively to fans of sizes $j,j+1,\dots,k-1$ gives the lemma. Let $\mathcal{P}$ be the given fan of size $j$, let $\mathcal{Q}$ be an $(x,B)$-fan of size $j+1$, and put $B'=\{b_1,\dots,b_j\}$. The $j+1$ paths of $\mathcal{Q}$ end at distinct vertices of $B$, so $B\setminus B'\neq\emptyset$. Let $G'$ be the graph obtained from $G$ by identifying all vertices of $B\setminus B'$ into a single new vertex $z$, deleting loops and parallel edges; thus $z$ is adjacent to every vertex of $V(G)\setminus(B\setminus B')$ having a neighbor in $B\setminus B'$. Put $T=B'\cup\{z\}$.

    We claim that $G'$ contains an $(x,T)$-fan of size $j+1$. Otherwise, by the fan form of Menger's theorem~\cite[Section~3.3]{Diestel}, some set $S\subseteq V(G')\setminus\{x\}$ with $|S|\le j$ meets every path from $x$ to $T$. The paths of $\mathcal{P}$ avoid $B\setminus B'$, so they are paths from $x$ to $T$ in $G'$ that pairwise intersect only in $x$; hence $S$ consists of exactly one vertex of each of them, and $z\notin S$. In particular $S$ misses some path $Q\in\mathcal{Q}$. As the internal vertices of $Q$ avoid $B$, the image of $Q$ in $G'$ is a path from $x$ to $b_i$ or to $z$, which avoids $S$, a contradiction.

    Let $\mathcal{R}$ be an $(x,T)$-fan of size $j+1$ in $G'$. Its paths end at the $j+1$ vertices of $T$, and their internal vertices avoid $T$, so $z$ occurs only as the last vertex of one path $R$. Hence the paths ending in $B'$ are paths of $G$ whose internal vertices avoid $B$, and if $yz$ is the last edge of $R$, replacing it by an edge $yb$ of $G$ with $b\in B\setminus B'$ gives a path of $G$ ending at $b$ whose internal vertices avoid $B$. This is the required $(x,B)$-fan of size $j+1$.
\end{proof}

We use $5$-connectivity mostly in the following form.

\begin{lemma}\label{lem:hpath-across}
    Let $H$ be a subgraph of a graph $G$, let $Z \subseteq V(G)$ be such that $G - Z$ is connected, and let $X \subseteq V(H) \setminus Z$ be nonempty, with $V(H) \setminus (X \cup Z) \neq \emptyset$ and every neighbor in $H$ of a vertex of $X$ lying in $X \cup Z$. Then there is an $H$-path from $X$ to $V(H) \setminus (X \cup Z)$.
\end{lemma}

\begin{proof}
    Take a path of $G - Z$ from $X$ to $V(H) \setminus (X \cup Z)$, let $y$ be its first vertex in $V(H) \setminus (X \cup Z)$, and let $x$ be its last vertex in $X$ before $y$. The subpath from $x$ to $y$ has no internal vertex in $V(H)$, as $Z$ is deleted, and it is not an edge of $H$, as $x$ has no neighbor in $H$ outside $X \cup Z$. So it is an $H$-path.
\end{proof}

In particular, if $G$ is $5$-connected and $|Z| \le 4$, then $G - Z$ is connected; this is what we mean when we say that a set $Z$ \emph{is not a separator} to obtain an $H$-path. Typically, $X$ is a union of interiors of legs of $H$ attached to the rest of $H$ only through $Z$, and the cited lemmas then restrict the ends of the $H$-path obtained; when $X$ also contains internal vertices of $H$-paths, as in (C8) of Lemma~\ref{crossed} and in Step~4 there, the argument is given explicitly.

\section{The Frame Structure Theorem}\label{sec:structure}

We first introduce frames, the main structural object used throughout the paper.

\begin{definition}[Frame]\label{def:frame}
    Let $a_1,a_2,a_3,a_4$ be distinct vertices of $G$. A subgraph $H = S \cup T$ of $G$ is a \emph{frame} rooted at $\{a_1,a_2,a_3,a_4\}$ if there exist vertices $s_1,s_2,t$ with $t\notin\{s_1,s_2\}$ and $s_1,s_2,t\notin\{a_1,a_2,a_3,a_4\}$ and indices $i_1,i_2,i_3,i_4\in\{1,2\}$, with exactly two of the $i_r$ equal to $1$, such that
    \begin{itemize}
        \item $S=\bigcup_{r=0}^4 S_r$, where $S_0$ is a possibly trivial $(s_1,s_2)$-path and $S_r$ is a nonempty $(s_{i_r},a_r)$-path for each $r\in\{1,2,3,4\}$;
        \item $T=\bigcup_{r=1}^3 T_r$, where $T_r$ is a nonempty $(t,a_r)$-path for each $r\in\{1,2,3\}$;
        \item the paths $S_0,S_1,S_2,S_3,S_4,T_1,T_2,T_3$ only meet at their explicitly prescribed common endpoints.
    \end{itemize}

    A \emph{pointed frame} is a triple $(H,a_3,a_4)$.
    Here $H$ is rooted at $\{a_1,a_2,a_3,a_4\}$; the root $a_4$ is the
    \emph{off-tripod root}, namely the unique root not joined to the apex $t$;
    and $a_3$ is a chosen \emph{distinguished tripod root}.
    The two remaining roots $a_1,a_2$ form the \emph{symmetric pair}.

    The frame is \emph{reduced} if $s_1=s_2$.
    If $s_1 \neq s_2$, the frame is \emph{classical} when $i_3=i_4$, and \emph{crossed} otherwise.

    When a frame is fixed in the proof below, we use the following convention.
    We relabel $s_1$ and $s_2$, if necessary, so that $S_4$ starts at $s_2$
    (equivalently, after this relabelling, $i_4=2$). We then keep the notation
    $s_1,s_2$ throughout, without introducing an additional shorthand.
    In the reduced case this simply means $s_1=s_2$.

    A frame is \emph{$S_0$-minimal} if it minimizes $|E(S_0)|$ over all frames rooted at $\{a_1,a_2,a_3,a_4\}$. See Figure~\ref{fig:frame}.

    \begin{figure}[htbp]
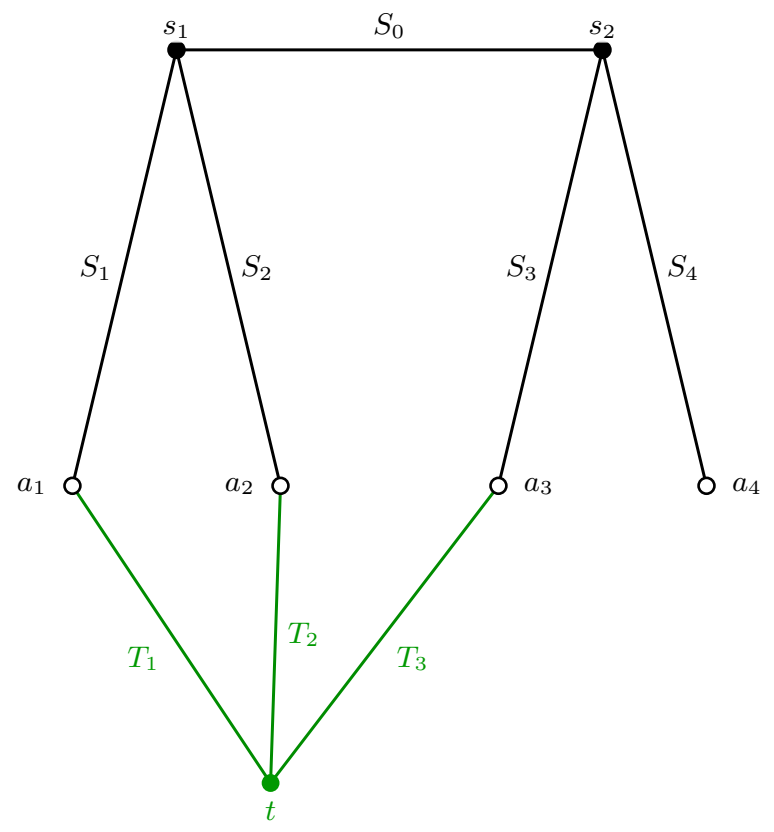

        \centering
        \compactfiginput{figures/definition_graph.tex}
        \caption{A pointed frame and its naming conventions. Throughout, the $S$-side is drawn in black and the $T$-side in green; hubs are solid black, the apex $t$ solid green, and the roots $a_1,\dots,a_4$ hollow. Auxiliary paths are red; in types $D$ and $E$, the paths of $U$ are blue.}
        \label{fig:frame}
    \end{figure}

\end{definition}

\begin{definition}[Legless tripod]
    Let $a_1,a_2,a_3$ be distinct vertices of a graph $G$.
    A subgraph $L$ of $G$ is a \emph{legless tripod} on $a_1,a_2,a_3$ with \emph{branch vertices} $s,t$ if $s,t,a_1,a_2,a_3$ are distinct and
    \[
        L = S \cup T,
        \qquad
        S = \bigcup_{i=1}^3 S_i,
        \qquad
        T = \bigcup_{i=1}^3 T_i,
    \]
    where
    \begin{itemize}
        \item $S_i$ is an $(s,a_i)$-path and $T_i$ is a $(t,a_i)$-path for each $i \in \{1,2,3\}$;
        \item the paths $S_1,S_2,S_3$ meet pairwise only in $s$, and the paths $T_1,T_2,T_3$ meet pairwise only in $t$;
        \item $S \cap T = \{a_1,a_2,a_3\}$.
    \end{itemize}
\end{definition}

\begin{lemma}[Robertson--Seymour--Thomas~\cite{RST1993}]\label{Seymour}
    Let $G$ be a $4$-connected graph, and let $a_1,a_2,a_3$ be distinct vertices of $G$.
    Then $G$ contains a legless tripod on $a_1,a_2,a_3$ or $G$ has a planar embedding with $a_1,a_2,a_3$ on the outer face.
\end{lemma}
This follows from (3.5) and (3.4) of~\cite{RST1993}: since $G$ is $4$-connected, it has no $(\le 2)$-separation $(A,B)$ with $a_1,a_2,a_3 \in A$ and $|B \setminus A| \ge 2$, so by (3.5) either $G$ can be drawn in a disk with $a_1,a_2,a_3$ on the boundary, or $G$ contains a tripod on $a_1,a_2,a_3$; and since $G$ has no $3$-separation, (3.4) turns such a tripod into a legless one. A legless tripod in the sense of~\cite{RST1993} is exactly a legless tripod as defined above, with $s,t$ the two vertices joined by the three paths $S_i \cup T_i$.

\begin{lemma}\label{lem:legless-tripod-to-frame}
    Let $G$ be a $5$-connected graph, and let $a_1,a_2,a_3,a_4$ be distinct vertices of $G$. Suppose that $G - \{a_3\}$ has no planar embedding with $a_1,a_2,a_4$ on the outer face, or that $G - \{a_4\}$ has no planar embedding with $a_1,a_2,a_3$ on the outer face.
    Then, after possibly exchanging $a_3$ and $a_4$, $G$ contains a pointed frame $(H,a_3,a_4)$ whose underlying frame $H$ is rooted at $\{a_1,a_2,a_3,a_4\}$.
\end{lemma}

\begin{proof}
    Up to relabelling, we may assume that $G - \{a_4\}$ has no planar embedding with $a_1,a_2,a_3$ on the outer face.
    Since $G-\{a_4\}$ is $4$-connected, Lemma~\ref{Seymour} and the hypothesis imply that $G-\{a_4\}$ contains a legless tripod $L=S\cup T$ on $a_1,a_2,a_3$.
    Since $L\subseteq G-\{a_4\}$, we have $a_4\notin V(L)$, and $|V(L)|\ge5$.
    As $G$ is $5$-connected, Lemma~\ref{lem:menger-fan} yields an
    $(a_4,V(L))$-fan of size four. At most three of its paths end in
    $\{a_1,a_2,a_3\}$, so at least one, call it $P$, ends at a vertex
    $x\in V(L)\setminus\{a_1,a_2,a_3\}$ and is internally disjoint from $L$.
    By symmetry between the two sides of the tripod, we may assume that $x\in V(S)\setminus\{a_1,a_2,a_3\}$.
    Therefore, $L\cup P$ forms a pointed frame $(H,a_3,a_4)$.
\end{proof}

The following definition collects the five configuration types produced by the proof. Type~$A$ replaces the tripod by a second tree reaching all four roots. In types~$B$ and~$C$, the tree $S$ keeps its four legs, while one leg of $T$ ends on an interior vertex of a leg of $S$ instead of at the corresponding root, and further paths join $S$ and $T$. In types~$D$ and~$E$, a third tree $U$, whose legs end on legs of $S$ and $T$, is added, and two or three paths link the three trees. Each type is designed so that, together with a $K_4^-$ on the roots, it yields a $K_6^-$-minor (Lemma~\ref{K4minusgivesK6minus}).

\begin{definition}[Types relative to $H$]
    Let $(H,a_3,a_4)$ be a pointed frame rooted at $\{a_1,a_2,a_3,a_4\}$. A subgraph $K$ of $G$ is of \emph{type $X$}, where $X \in \{A,B,C,D,E\}$, if it satisfies the conditions listed under Type $X$ below; see Figure~\ref{fig:types}. In each type, the named vertices are pairwise distinct, except that $s_1 = s_2$ if $S_0$ is trivial, $t_1 = t_2$ if $T_0$ is trivial, and two named vertices may coincide when a membership condition below allows it (for instance $a'_3 = a_3$ or $u_1 = a_4$ in type~$B$, $a'_2 = a_2$, $u = a_3$ or $v = a_4$ in type~$C$, $u_r = a_r$ in types~$D$ and~$E$, and the common ends of the paths $V_i$ described in types~$D$ and~$E$). In particular, no hub ($s_1,s_2,s,t,t_1,t_2$, and $u$ in types~$D$ and~$E$) is a root. A nonempty path has at least one edge, and ``disjoint except at the specified endpoints'' means that two of the listed paths meet only in vertices prescribed by the definition, namely common ends of both and ends of one of them that are required to lie on the other.
    \begin{description}
        \item[Type $A$] A subgraph $K$ is of \emph{type $A$} if
        \[
            K = S \cup T,
        \]
        and there exist vertices $s_1,s_2,t_1,t_2$ and indices $i_1,\dots,i_4,j_1,\dots,j_4 \in \{1,2\}$, with exactly two of the $i_r$ and exactly two of the $j_r$ equal to $1$, for which
        \begin{itemize}
            \item $S = \bigcup_{r=0}^4 S_r$, where $S_0$ is a possibly trivial $(s_1,s_2)$-path and $S_r$ is a nonempty $(s_{i_r},a_r)$-path for each $r \in \{1,2,3,4\}$;
            \item $T = \bigcup_{r=0}^4 T_r$, where $T_0$ is a possibly trivial $(t_1,t_2)$-path and $T_r$ is a nonempty $(t_{j_r},a_r)$-path for each $r \in \{1,2,3,4\}$.
        \end{itemize}
        Furthermore, the paths $S_0,S_1,S_2,S_3,S_4,T_0,T_1,T_2,T_3,T_4$ are pairwise disjoint except at the specified endpoints.

        \item[Type $B$] A subgraph $K$ is of \emph{type $B$} if
        \[
            K = S \cup T \cup U,
        \]
        and there exist vertices $s_1,s_2,t,a'_3,u_1,u_2,v_1,v_2$ and indices $i_1,\dots,i_4 \in \{1,2\}$, with exactly two of the indices equal to $1$, for which, after possibly interchanging the labels $a_1$ and $a_2$, the following holds:
        \begin{itemize}
            \item $S = \bigcup_{r=0}^4 S_r$, where $S_0$ is a possibly trivial $(s_1,s_2)$-path and $S_r$ is a nonempty $(s_{i_r},a_r)$-path for each $r \in \{1,2,3,4\}$;
            \item $T = T_2 \cup T_3 \cup T_4$, where $T_r$ is a nonempty $(t,a_r)$-path for $r \in \{2,4\}$, and $T_3$ is a nonempty $(t,a'_3)$-path with $a'_3 \in V(S_3(s_{i_3},a_3])$;
            \item $U = U_1 \cup U_2$, where $U_1$ and $U_2$ are nonempty $(u_1,v_1)$- and $(u_2,v_2)$-paths such that $u_1 \in V(S_4(s_{i_4},a_4])$, $v_1 \in V(T_3(t,a'_3] \cup S_3(a'_3,a_3])$, $u_2 \in V(S_0 \cup S_1 \cup S_2 \cup S_3[s_{i_3},a'_3) \cup S_4[s_{i_4},u_1)) \setminus \{a_1,a_2\}$ and $v_2 \in V(T_2 \cup T_3[t,w) \cup T_4) \setminus \{a_2,a_4\}$, where $w = v_1$ if $v_1 \in V(T_3)$ and $w = a'_3$ otherwise.
        \end{itemize}
        Furthermore, the paths $S_0$, $S_1$, $S_2$, $S_3$, $S_4$, $T_2$, $T_3$, $T_4$, $U_1$, $U_2$ are pairwise disjoint except at the specified endpoints.

        \item[Type $C$] A subgraph $K$ is of \emph{type $C$} if
        \[
            K = S \cup T \cup U,
        \]
        and there exist vertices $s_1,s_2,t,u,v,a'_2$ and indices $i_1,\dots,i_4 \in \{1,2\}$, with exactly two of the $i_r$ equal to $1$, for which, after possibly interchanging the labels $a_1$ and $a_2$, the following holds:
        \begin{itemize}
            \item $S = \bigcup_{r=0}^4 S_r$, where $S_0$ is a possibly trivial $(s_1,s_2)$-path and $S_r$ is a nonempty $(s_{i_r},a_r)$-path for each $r \in \{1,2,3,4\}$;
            \item $T = T_1 \cup T_2 \cup T_3$, where $T_1$ and $T_3$ are nonempty $(t,a_1)$- and $(t,a_3)$-paths, and $T_2$ is a nonempty $(t,a'_2)$-path with $a'_2 \in V(S_2(s_{i_2},a_2])$;
            \item $U$ is a nonempty $(u,v)$-path, where $u \in V(S_3(s_{i_3},a_3])$ and $v \in V(S_4(s_{i_4},a_4])$.
        \end{itemize}
        Furthermore, the paths $S_0$, $S_1$, $S_2$, $S_3$, $S_4$, $T_1$, $T_2$, $T_3$, $U$ are pairwise disjoint except at the specified endpoints.

        \item[Type $D$] A subgraph $K$ is of \emph{type $D$} if
        \[
            K = S \cup T \cup U \cup V,
        \]
        and there exist vertices $s,t,u,u_2,u_3$ and $v_1,v'_1,v_2,v'_2,v_3,v'_3$ such that, after possibly interchanging the labels $a_1$ and $a_2$, the following holds:
        \begin{itemize}
            \item $S = S_1 \cup S_2 \cup S_4$, where $S_r$ is a nonempty $(s,a_r)$-path for $r = 1,2,4$;
            \item $T = T_2 \cup T_3 \cup T_4$, where $T_r$ is a nonempty $(t,a_r)$-path for $r = 2,3,4$;
            \item $U = U_2 \cup U_3$, where $U_2$ and $U_3$ are nonempty $(u,u_2)$- and $(u,u_3)$-paths, with $u_2 \in V(S_2(s,a_2])$ and $u_3 \in V(T_3(t,a_3])$;
            \item $V = V_1 \cup V_2 \cup V_3$, where $V_1$, $V_2$, and $V_3$ are nonempty $(v_1,v'_1)$-, $(v_2,v'_2)$-, and $(v_3,v'_3)$-paths with $v_1, v_2 \in V(S_1[s,a_1) \cup S_2[s,u_2) \cup S_4[s,a_4))$, $v'_1, v_3 \in V(T_2[t,a_2) \cup T_3[t,u_3) \cup T_4[t,a_4))$, and $v'_2, v'_3 \in V(U_2[u,u_2) \cup U_3[u,u_3))$.
        \end{itemize}
        Furthermore, the paths $S_1$, $S_2$, $S_4$, $T_2$, $T_3$, $T_4$, $U_2$, $U_3$ are pairwise disjoint except at the specified endpoints, and each of $V_1$, $V_2$, $V_3$ is internally disjoint from $S \cup T \cup U$.
        The paths $V_1$, $V_2$, $V_3$ may intersect only as follows: $V_1 \cap V_2$ is empty or a common subpath of $V_1$ and $V_2$ starting at $v_1 = v_2$; $V_1 \cap V_3$ is empty or a common subpath of $V_1$ and $V_3$ starting at $v'_1 = v_3$; $V_2 \cap V_3$ is empty or a common subpath of $V_2$ and $V_3$ starting at $v'_2 = v'_3$; and no vertex lies on all three.

        \item[Type $E$] A subgraph $K$ is of \emph{type $E$} if
        \[
            K = S \cup T \cup U \cup V,
        \]
        and there exist vertices $s,t,u,u_2,u_3,u_4,v_1,v'_1,v_2,v'_2$ such that, after possibly interchanging the labels $a_1$ and $a_2$, the following holds:
        \begin{itemize}
            \item $S = S_1 \cup S_2 \cup S_4$, where $S_r$ is a nonempty $(s,a_r)$-path for $r = 1,2,4$;
            \item $T = T_2 \cup T_3 \cup T_4$, where $T_r$ is a nonempty $(t,a_r)$-path for $r = 2,3,4$;
            \item $U = U_2 \cup U_3 \cup U_4$, where $U_r$ is a nonempty $(u,u_r)$-path for $r = 2,3,4$, with $u_r \in V(S_r(s,a_r] \cup T_r(t,a_r])$ for $r = 2,4$ and $u_3 \in V(T_3(t,a_3])$;
            \item $V = V_1 \cup V_2$, where $V_1$ and $V_2$ are nonempty $(v_1,v'_1)$- and $(v_2,v'_2)$-paths joining two distinct pairs among the sets
            \[
                \begin{gathered}
                    V(S_1[s,a_1) \cup S_2[s,u_2) \cup S_4[s,u_4)), \qquad V(T_2[t,u_2) \cup T_3[t,u_3) \cup T_4[t,u_4)),\\
                    V(U) \setminus \{u_2,u_3,u_4\},
                \end{gathered}
            \]
            where $S_r[s,u_r)$ stands for $S_r[s,a_r)$ if $u_r \notin V(S_r)$, and $T_r[t,u_r)$ for $T_r[t,a_r)$ if $u_r \notin V(T_r)$.
        \end{itemize}
        Furthermore, the paths $S_1$, $S_2$, $S_4$, $T_2$, $T_3$, $T_4$, $U_2$, $U_3$, and $U_4$ are pairwise disjoint except at the specified endpoints, each of $V_1$, $V_2$ is internally disjoint from $S \cup T \cup U$, and $V_1 \cap V_2$ is empty or a common end.
    \end{description}

    In types $D$ and $E$, we may also specify the connectors by \emph{routes}: for each $i$, a route $R_i$ is a path from the first of the two prescribed sets of $V_i$ to the second that meets $S \cup T \cup U$ only in these two sets, and the routes are required to satisfy only the pairwise intersection conditions imposed above on the $V_i$, not internal disjointness from $S \cup T \cup U$. Given such routes, we let $V_i$ be the subpath of $R_i$ from its last vertex in the first set to the next vertex in the second. The paths $V_i$ so obtained satisfy all the conditions of the definition: they are internally disjoint from $S \cup T \cup U$ by construction, and they still satisfy the intersection conditions, since a common subpath of two routes meets $S \cup T \cup U$ only in the prescribed set common to both.

    \begin{figure}[tbp]
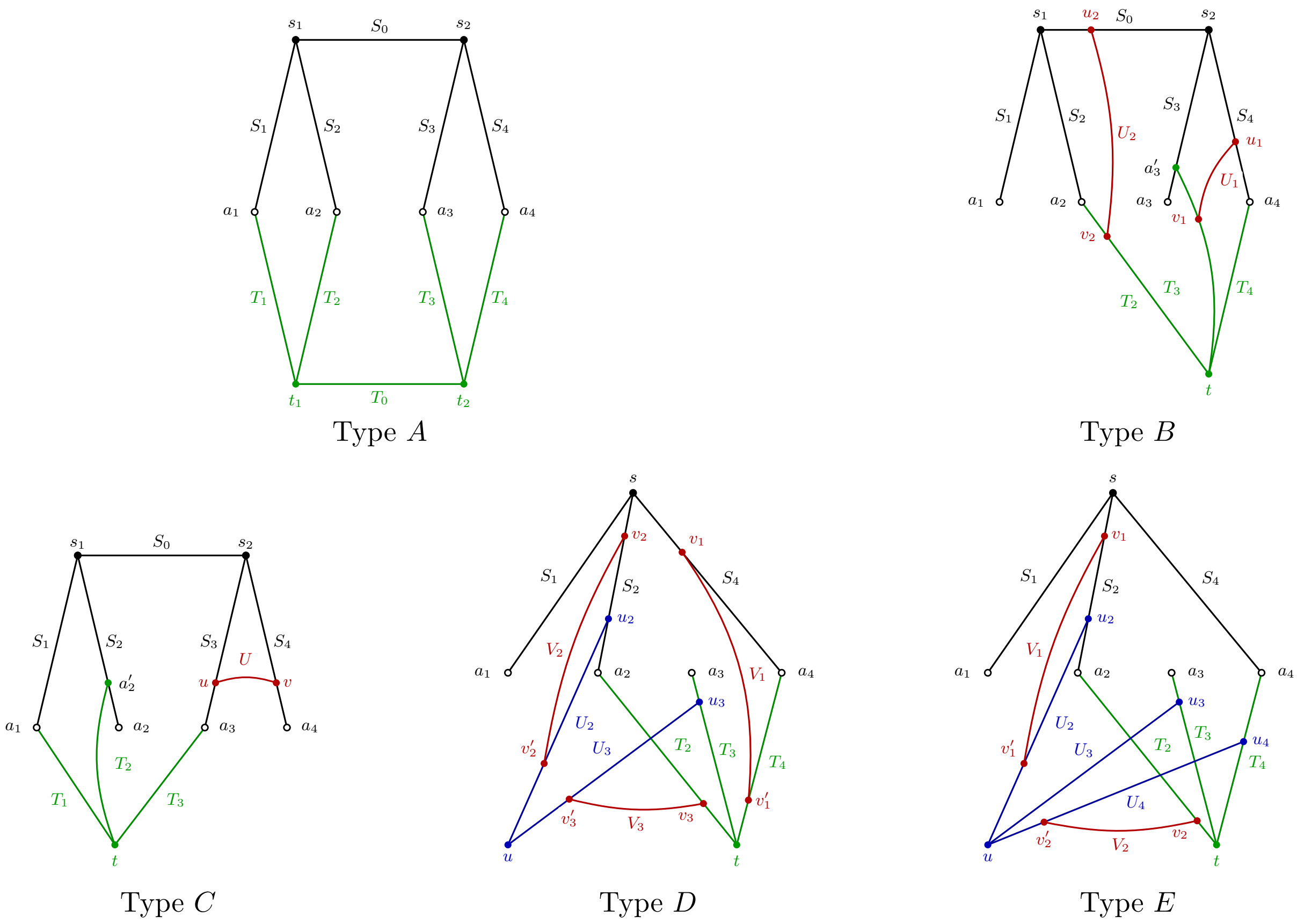

        \centering
        \begin{minipage}[b]{0.62\linewidth}\centering\typepanel{figures/type_a_graph.tex}\par\vspace{3pt}Type $A$\end{minipage}\hfill\begin{minipage}[b]{0.32\linewidth}\centering\typepanel{figures/type_b_graph.tex}\par\vspace{3pt}Type $B$\end{minipage}

\vspace{0.8em}

\begin{minipage}[b]{0.32\linewidth}\centering\typepanel{figures/type_c_graph.tex}\par\vspace{3pt}Type $C$\end{minipage}\hfill\begin{minipage}[b]{0.32\linewidth}\centering\typepanel{figures/type_d_graph.tex}\par\vspace{3pt}Type $D$\end{minipage}\hfill\begin{minipage}[b]{0.32\linewidth}\centering\typepanel{figures/type_e_graph.tex}\par\vspace{3pt}Type $E$\end{minipage}

        \caption{The five unavoidable configurations of types $A$--$E$.}
        \label{fig:types}
        \vspace{0.5em}
    \end{figure}

\end{definition}

In the remainder of this section, the goal is to prove the following result.

\begin{theorem}[Frame structure theorem]\label{structure}
    Let $G$ be a $5$-connected graph, and let $a_1,a_2,a_3,a_4$ be distinct vertices of $G$.
    If $G$ admits a pointed frame $(H,a_3,a_4)$ whose underlying frame $H$ is rooted at $\{a_1,a_2,a_3,a_4\}$, then $G$ contains a subgraph of type $X$ for some $X \in \{A,B,C,D,E\}$.
\end{theorem}

Fix a pointed frame $(H,a_3,a_4)$ whose underlying frame $H$ is rooted at $\{a_1,a_2,a_3,a_4\}$. The roots are ordered and $a_3$ is distinguished, while $a_1$ and $a_2$ remain symmetric; this is why, in types $B$--$E$, we allow ourselves to interchange their labels. No edge among the four roots is assumed in Theorem~\ref{structure}; edge information enters only in Section~\ref{sec:k6minus}, where the theorem is applied to a $K_4^-$ whose missing edge is $a_3a_4$.

We argue by contradiction: until the end of this section, assume that $G$ contains no subgraph of types $A$--$E$. We first use $5$-connectivity to attach two $H$-paths to the frame, one leaving $S_0$ and one leaving $a_4$, and call the result an augmented frame; we choose it so that the underlying frame is $S_0$-minimal (Lemma~\ref{lem:minimal-augmented-frame}). Lemma~\ref{classical} shows that a minimal augmented frame is crossed or reduced. Lemmas~\ref{S0min}--\ref{S3toT} restrict the $H$-paths leaving the legs of $S$, and Lemma~\ref{crossed} uses these restrictions to constrain the crossed case. Lemmas~\ref{S1pairs}--\ref{S1toT} treat the reduced case, together with the crossed case in which the second added path ends on $S_1$. Each lemma either excludes a configuration of $H$-paths or produces a subgraph of one of types $A$--$E$, and the proof of Theorem~\ref{structure} combines them.

\begin{definition}[Augmented frame]
    Let $(H,a_3,a_4)$ be a pointed frame with $H = S \cup T$, and with $s_1,s_2$
    fixed by the convention in the definition of a frame. An \emph{augmented
    frame} is a tuple
    \[
        H^* = \bigl(H, P(s_0,v_0), P(a_4,b_4)\bigr),
    \]
    where the two paths are disjoint $H$-paths, and
    \begin{itemize}
        \item $s_0 \in V(S_0)$, with $s_0 \neq s_1$ whenever $s_1 \neq s_2$
              (equivalently, $s_0 \in V(S_0(s_1,s_2])$ in the non-reduced case,
              and $s_0 = s_1 = s_2$ in the reduced case);
        \item $v_0 \in V(T_1\cup T_2\cup T_3)\setminus\{a_1,a_2,a_3\}$;
        \item $b_4 \in V(S_1 \cup S_2)\setminus\{a_1,a_2,s_1,s_2\}$.
    \end{itemize}
    The augmented frame is classical, crossed, or reduced according to the underlying frame $H$.
    The augmented frame is \emph{minimal} if its underlying frame is $S_0$-minimal.
\end{definition}

The first step produces the frame on which the whole argument runs.

\begin{lemma}[Minimal augmented frame]\label{lem:minimal-augmented-frame}
    The graph $G$ contains a minimal augmented frame.
\end{lemma}

\begin{proof}
    \emph{Step 1. There is an $S_0$-minimal frame $H$ rooted at $\{a_1,a_2,a_3,a_4\}$ admitting an $H$-path $P(s_0,v_0)$ with $s_0 \in V(S_0)$ and $v_0 \in V(T)\setminus\{a_1,a_2,a_3\}$. Moreover, if $s_1 \neq s_2$, then $s_0 \in V(S_0(s_1,s_2])$, i.e.\ $s_0 \neq s_1$.}

    Let $H$ be an $S_0$-minimal frame. We distinguish three cases according to $|E(S_0)|$.

    \emph{Case 1: $|E(S_0)| \ge 2$.}
    Then the set $V(S_0(s_1,s_2))$ is nonempty. Furthermore, as $\{s_1,s_2\}$ is not a separator, there is a path in $G-\{s_1,s_2\}$ from $V(S_0(s_1,s_2))$ to $V(H)\setminus V(S_0)$.
    Choosing such a path of minimum length gives an $H$-path $P(s_0,y)$ with $s_0 \in V(S_0(s_1,s_2))$ and $y \in V(H)\setminus V(S_0)$.

    If $y \in V(T) \setminus \{a_1,a_2,a_3\}$, we are done. Otherwise $y \in V(S_r(s_{i_r},a_r])$ for some $r \in \{1,2,3,4\}$.
    Suppose first that $s_{i_r} = s_1$, and let $S_{r'}$ be the other leg of
    $H$ at $s_1$. Then the subgraph with hubs $s_0$ and $s_2$, connecting path
    $S_0[s_0,s_2]$, legs
    \[
        P(s_0,y) \cup S_r[y,a_r]
        \qquad\text{and}\qquad
        S_0[s_0,s_1] \cup S_{r'}
    \]
    at $s_0$, and the two legs of $H$ at $s_2$, is a frame rooted at
    $\{a_1,a_2,a_3,a_4\}$ with $|E(S_0[s_0,s_2])| < |E(S_0)|$, since
    $s_0 \neq s_1$. This contradicts the minimality of $|E(S_0)|$.
    The case $s_{i_r} = s_2$ is symmetric, with hubs $s_1$ and $s_0$ and
    connecting path $S_0[s_1,s_0]$. Therefore $y \in V(T) \setminus \{a_1,a_2,a_3\}$, and
    $P(s_0,y)$ is the desired path.

    \emph{Case 2: $|E(S_0)| = 1$.}
    Let $j \in \{1,2,3\}$ be the index for which $S_j$ starts at $s_2$, and
    let $\ell,m$ be the two indices in $\{1,2,3\}\setminus\{j\}$. Set
    \[
        B = \{s_1,a_1,a_2,a_3,a_4\} \cup \bigl(V(T) \setminus \{a_1,a_2,a_3\}\bigr)
            \cup V(S_\ell(s_1,a_\ell)) \cup V(S_m(s_1,a_m)).
    \]
    The paths $S_j$, $S_4$, and the edge $S_0$ are three internally disjoint
    paths from $s_2$ to $a_j$, $a_4$, $s_1$ whose internal vertices avoid $B$.
    Since $G$ is $5$-connected and $|B| \ge 5$, Lemma~\ref{lem:menger-fan}
    together with Lemma~\ref{perfect} yields an $(s_2,B)$-fan of size
    four, three of whose paths end at $a_j$, $a_4$ and $s_1$; call the first
    two $S'_j$ and $S'_4$ (the third is not used, the edge $S_0$ playing its role),
    while the fourth, $R$, ends at some
    \[
        x \in \{a_\ell,a_m\} \cup V(S_\ell(s_1,a_\ell)) \cup V(S_m(s_1,a_m))
            \cup \bigl(V(T) \setminus \{a_1,a_2,a_3\}\bigr).
    \]
    All four paths have their internal vertices outside $B$; in particular
    they avoid $V(T)$, the vertex $s_1$, and $V(S_\ell \cup S_m)\setminus\{s_1\}$.

    If $x \in V(S_\ell(s_1,a_\ell])$ (the case $x \in V(S_m(s_1,a_m])$ being
    symmetric), then the four paths
    \[
        S'_j, \qquad S'_4, \qquad R \cup S_\ell[x,a_\ell],
        \qquad S_0 \cup S_m
    \]
    pairwise meet only in $s_2$ and meet $T$ only in roots, so together with
    $T$ they form a reduced frame rooted at $\{a_1,a_2,a_3,a_4\}$,
    contradicting the minimality of $|E(S_0)| = 1$.
    Hence, $x \in V(T) \setminus \{a_1,a_2,a_3\}$. Then
    $H' := S_0 \cup S'_j \cup S'_4 \cup S_\ell \cup S_m \cup T$ is a frame
    with the same connecting path $S_0$, hence $S_0$-minimal, and $R$ is an
    $H'$-path from $s_2 \in V(S_0)$ to $V(T) \setminus \{a_1,a_2,a_3\}$.

    \emph{Case 3: $|E(S_0)| = 0$.}
    Set $B = \{a_1,a_2,a_3,a_4\} \cup (V(T) \setminus \{a_1,a_2,a_3\})$. The legs $S_1,\dots,S_4$ are
    four internally disjoint paths from $s_1 = s_2$ to $B$ whose internal
    vertices avoid $B$. As above, Lemma~\ref{lem:menger-fan} and
    Lemma~\ref{perfect} yield an $(s_1,B)$-fan of size five, four of whose
    paths, $S'_1,\dots,S'_4$, end at $a_1,\dots,a_4$, while the
    fifth, $R$, ends at some $v_0 \in V(T) \setminus \{a_1,a_2,a_3\}$. All five paths have their
    internal vertices outside $B$, hence outside $V(T)$, so
    $H' := S'_1 \cup S'_2 \cup S'_3 \cup S'_4 \cup T$ is a reduced frame,
    and $R$ is an $H'$-path from $s_1 \in V(S_0)$ to $V(T) \setminus \{a_1,a_2,a_3\}$.

    \emph{Step 2.} By Step~1, there is an $S_0$-minimal frame $H$ and an $H$-path
    $P(s_0,v_0)$. Here $s_0 \in V(S_0)$, $s_0 \neq s_1$ when $s_1 \neq s_2$, and
    $v_0 \in V(T)\setminus\{a_1,a_2,a_3\}$.

    Write $s := s_{i_4}$, and, if $s_1 \neq s_2$, let $s'$ be the other hub. Set
    \[
        B := \bigl(V(H)\setminus V(S_4(s,a_4])\bigr) \cup V\bigl(P(s_0,v_0)\bigr).
    \]
    Thus, $B$ contains $s$, the paths $S_0,S_1,S_2,S_3$, all of $T$, and every
    vertex of $P(s_0,v_0)$, but neither $a_4$ nor any interior vertex of the
    leg $S_4$. In particular $|B| \ge 5$.

    The leg $S_4$, traversed from $a_4$, is an $(a_4,B)$-fan of size one ending
    at $s$: its interior is $V(S_4(s,a_4))$, which is disjoint from $V(H)\setminus
    V(S_4(s,a_4])$ by definition of $B$, and disjoint from $V(P(s_0,v_0))$ since
    the latter meets $V(H)$ only in $s_0$ and $v_0$. Since $G$ is $5$-connected
    and $|B|\ge5$, Lemma~\ref{lem:menger-fan} gives an $(a_4,B)$-fan of size five;
    combined with the fan of size one via Lemma~\ref{perfect} (case $j=1$, $k=5$),
    we obtain an $(a_4,B)$-fan $R_1,\dots,R_5$ of size five in which $R_1$ ends
    at $s$. These paths meet pairwise only in $a_4$, end at five distinct vertices
    of $B$, and have all their interior vertices outside $B$.

    Put
    \[
        H' := S_0 \cup S_1 \cup S_2 \cup S_3 \cup T \cup R_1.
    \]
    Then $H'$ is a frame rooted at $\{a_1,a_2,a_3,a_4\}$: the path $R_1$ joins
    $a_4$ to $s$ and its interior avoids $B$, hence avoids $S_0,S_1,S_2,S_3$ and
    $T$, so $R_1$ may serve as the leg to $a_4$ at $s$. The connecting path $S_0$
    is unchanged, so $H'$ is again $S_0$-minimal, and the indices $i_1,\dots,i_4$
    are unchanged, so $H'$ is reduced, classical, or crossed exactly as $H$ is.
    Moreover, the interior of $R_1$ avoids $V(P(s_0,v_0))$, so $P(s_0,v_0)$ is
    still an $H'$-path. We rename $H'$ as $H$ and its leg $R_1$ as $S_4$.

    The paths $R_2,\dots,R_5$ end at four distinct vertices of $B$, none of them
    equal to $s$, at most two of which lie in $\{a_1,a_2\}$; so at least one of
    them, say $R$, ends at a vertex $b_4 \in B\setminus\{s,a_1,a_2\}$.
    The interior of $R$ avoids $B$ and $R_1$, hence avoids $H$ and $P(s_0,v_0)$; thus $R$ is an $H$-path from $a_4$ to $b_4$, disjoint from $P(s_0,v_0)$ except possibly at $b_4$.

    Suppose that $H$ is not reduced and that $b_4 \in V(S_0[s',s))$. Let $S_r$ be the leg of $H$ at $s$ other than $S_4$.
    If $b_4 = s'$, then the two legs of $H$ at $s'$ and the legs $R$ and $S_0 \cup S_r$, all at $s'$, together with $T$, form a reduced frame rooted at $\{a_1,a_2,a_3,a_4\}$.
    Otherwise, the hubs $s'$ and $b_4$, the connecting path $S_0[s',b_4]$, the two legs of $H$ at $s'$, and the legs $R$ and $S_0[b_4,s] \cup S_r$ at $b_4$, together with $T$, form a frame rooted at $\{a_1,a_2,a_3,a_4\}$.
    In both cases the connecting path is shorter than $S_0$, contrary to the $S_0$-minimality of $H$.

    If $b_4 \in V(T) \setminus \{a_1,a_2,a_3\}$, then $S$ and $T \cup R$ form a subgraph of type $A$.
    If $b_4$ is an interior vertex of $P(s_0,v_0)$, then $S$ and $T \cup P[v_0,b_4] \cup R$ form a subgraph of type $A$.
    If $b_4 \in V(S_3(s_{i_3},a_3])$, then $S \cup T \cup R$ is of type $C$, with $U = R$, $u = b_4$, $v = a_4$ and $a'_2 = a_2$.

    Since $V(P(s_0,v_0)) \cap V(H) = \{s_0,v_0\}$, with $s_0 \in V(S_0)$ and $v_0 \in V(T)$, we conclude that
    \[
        b_4 \in V(S_1\cup S_2)
              \setminus\{a_1,a_2,s_1,s_2\};
    \]
    in particular, $b_4 \notin V(P(s_0,v_0))$, so $R$ and $P(s_0,v_0)$ are disjoint. Writing $P(a_4,b_4) = R$,
    \[
        H^* = \bigl(H, P(s_0,v_0), P(a_4,b_4)\bigr)
    \]
    is an augmented frame whose underlying frame is $S_0$-minimal, thus a minimal augmented frame.
\end{proof}

From now on, $H^* = \bigl(H, P(s_0,v_0), P(a_4,b_4)\bigr)$ denotes an arbitrary minimal augmented frame; one exists by Lemma~\ref{lem:minimal-augmented-frame}. Each lemma below is proved for every minimal augmented frame satisfying its hypotheses, so in a proof we may apply the previous lemmas to the minimal augmented frames obtained from $H^*$ by re-routing some of its paths.

For a vertex $v$ of $T = T_1 \cup T_2 \cup T_3$ and a root $a_i$ with $i \in \{1,2,3\}$, we write $T[v,a_i]$ for the $(v,a_i)$-path in the tree $T$.

In the figures of this section, a constructed subgraph $K$ is shown by shading: $S_K$ in grey, $T_K$ in green, $U_K$ in red for types $B$ and $C$ and in blue for types $D$ and $E$, and $V_K$ in red. The hubs of $S_K$, $T_K$ and $U_K$ are their branch vertices; in general they differ from the hubs of $H$, and when we check that $K$ is of a given type, the legs $S_r$, $T_r$ and the intervals such as $S_r(s_{i_r},a_r]$ in the definition refer to the legs of $K$, not to those of $H$.

If two paths $P$ and $Q$ meet and $x$ is an end of $Q$, we say that we \emph{cut $Q$ at $P$ from $x$} when we replace $Q$ by its subpath from $x$ to its first vertex on $P$; the new path meets $P$ only in its last vertex.

\begin{lemma} \label{classical}
    $H^*$ is not classical.
\end{lemma}

\begin{proof}
    Suppose, for a contradiction, that $H$ is classical, so that $s_1 \neq s_2$ and $i_1 = i_2 = 1$.
    Let
    \[
        \begin{aligned}
            S_K ={}& S_0 \cup S_1 \cup S_2 \cup S_3 \cup P[b_4,a_4], \\
            T_K ={}& T_1 \cup T_2 \cup T_3, \\
            U_K ={}& S_4.
        \end{aligned}
    \]
    Then $K = S_K \cup T_K \cup U_K$ is of type $C$ as in Figure~\ref{fig:classical-type-c}(a). Indeed, the hubs of $S_K$ are $s_1$ and $b_4$, and the leg of $S_K$ to $a_3$ is $S_0 \cup S_3$, which contains $s_2$ as an internal vertex; so $U_K = S_4$ joins $u = s_2$ on this leg to $v = a_4$, and $a'_2 = a_2$.

    \begin{figure}[htbp]
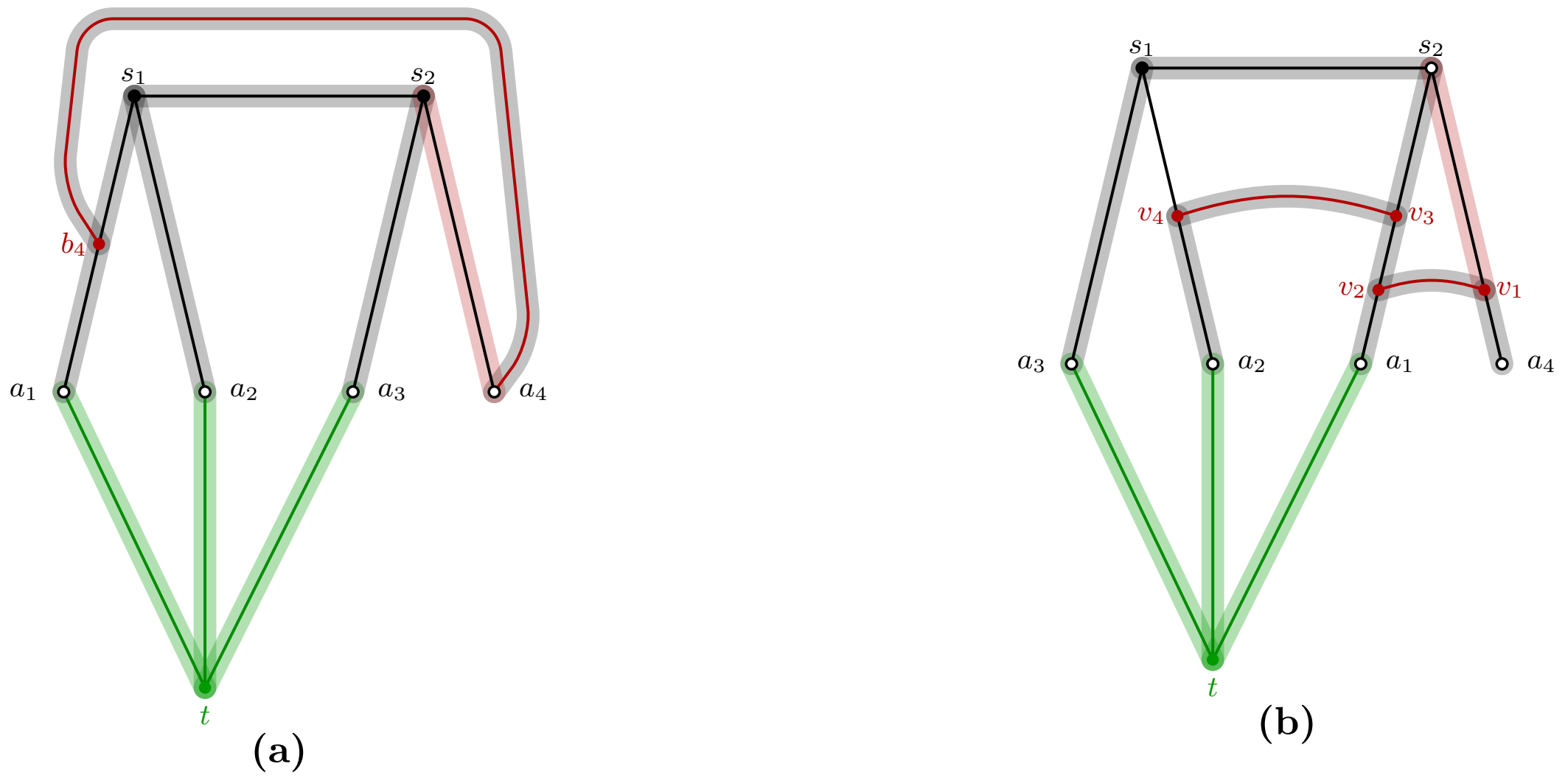

        \centering
        \begin{minipage}{0.48\textwidth}
            \centering
            \compactsubfiginput{figures/sequalssp.tex}\par
            \textbf{(a)}
        \end{minipage}\hfill
        \begin{minipage}{0.48\textwidth}
            \centering
            \compactsubfiginput{figures/S1toS2_graph.tex}\par
            \textbf{(b)}
        \end{minipage}
        \caption{(a) The type-$C$ configuration in Lemma~\ref{classical}. (b) The type-$C$ configuration in Lemma~\ref{S1toS2}.}
        \label{fig:classical-type-c}
    \end{figure}

\end{proof}

\begin{lemma}\label{S0min}
    Suppose that $H^*$ is crossed. Let $q \in \{1,2,3,4\}$, let $s \in \{s_1,s_2\}$ be the hub at which $S_q$ starts, and let $s'$ be the other hub. There is no $H$-path $P(v_1,v_2)$ with $v_1 \in V(S_q(s,a_q])$ and $v_2 \in V(S_0) \setminus \{s\}$.
\end{lemma}

\begin{proof}
    Let $S_r$ be the other leg of $H$ at $s$.
    If $v_2 = s'$, then the legs $P[s',v_1] \cup S_q[v_1,a_q]$ and $S_0 \cup S_r$ at $s'$, together with the two legs of $H$ at $s'$ and $T$, form a reduced frame rooted at $\{a_1,a_2,a_3,a_4\}$, contrary to the $S_0$-minimality of $H$, since $S_0$ is not trivial.
    Otherwise $v_2 \in V(S_0(s_1,s_2))$, and the hubs $v_2$ and $s'$, the connecting path $S_0[v_2,s']$, the legs $P[v_2,v_1] \cup S_q[v_1,a_q]$ and $S_0[v_2,s] \cup S_r$ at $v_2$, and the two legs of $H$ at $s'$, together with $T$, form a frame rooted at $\{a_1,a_2,a_3,a_4\}$ whose connecting path is shorter than $S_0$, contrary to the $S_0$-minimality of $H$.
\end{proof}

By Lemma~\ref{classical}, after interchanging the labels of $S_1$ and $S_2$ if necessary, $H^*$ is crossed or reduced, i.e., $i_1 = i_4 = 2$ and $i_2 = i_3 = 1$.

\begin{lemma}\label{S4toS3T1T2T3}
    There is no $H$-path $P(v_1,v_2)$ with $v_1 \in V(S_4(s_2,a_4])$ and $v_2 \in V(S_3 \cup T_1 \cup T_2 \cup T_3) \setminus \{s_1,a_1,a_2\}$.
\end{lemma}

\begin{proof}
    If $P(v_1,v_2)$ meets $P(a_4,b_4)$, let $y$ be its first vertex on $P(a_4,b_4)$ from $v_2$, and let $Q = P[v_2,y] \cup P[y,a_4]$, an $H$-path from $v_2$ to $a_4$. If $v_2 \in V(T_1 \cup T_2 \cup T_3) \setminus \{a_1,a_2,a_3\}$, then $S$ and $T \cup Q$ form a subgraph of type $A$; if $v_2 \in V(S_3(s_1,a_3])$ then $S \cup T \cup Q$ is of type $C$ with $U = Q$. 
    Hence, $P(v_1,v_2)$ and $P(a_4,b_4)$ are disjoint. 

    First assume that $v_2 \in V(T_1 \cup T_2 \cup T_3) \setminus \{a_1,a_2,a_3\}$.
    Let
    \[
        \begin{aligned}
            S_K ={}& S_0 \cup S_1 \cup S_2 \cup S_3 \cup P[b_4,a_4], \\
            T_K ={}& T_1 \cup T_2 \cup T_3 \cup P[v_1,v_2] \cup S_4[v_1,a_4].
        \end{aligned}
    \]
    Then $K = S_K \cup T_K$ is of type $A$, as in Figure~\ref{fig:s4-to-s3-t}(a).

    Now suppose that $v_2 \in V(S_3(s_1,a_3])$.
    Let
    \[
        \begin{aligned}
            S_K ={}& S_0 \cup S_1 \cup S_2 \cup S_3 \cup P[b_4,a_4], \\
            T_K ={}& T_1 \cup T_2 \cup T_3, \\
            U_K ={}& P[v_2,v_1] \cup S_4[v_1,a_4].
        \end{aligned}
    \]
    Then $K = S_K \cup T_K \cup U_K$ is of type $C$, as in Figure~\ref{fig:s4-to-s3-t}(b).

    \begin{figure}[htbp]
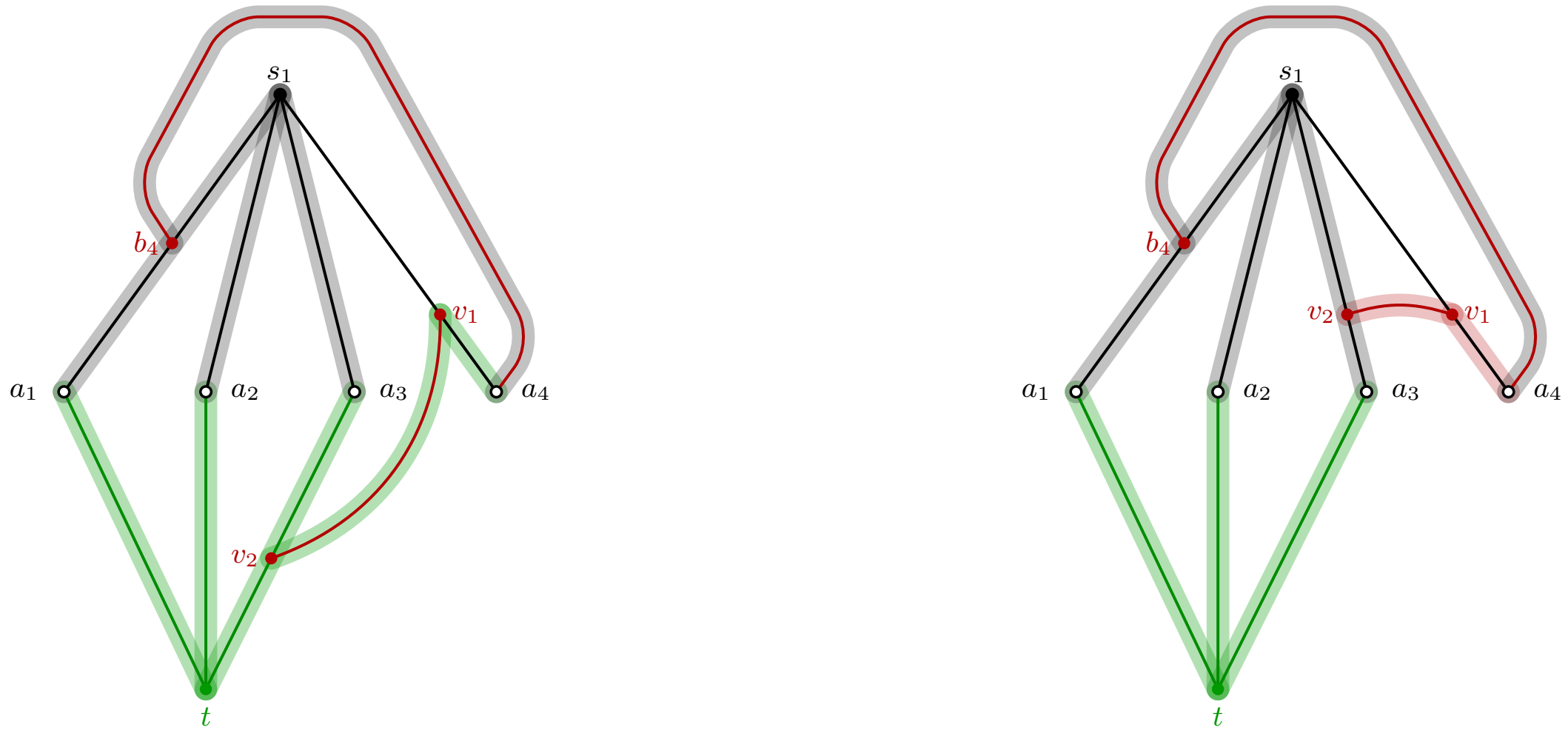

        \centering
        \begin{minipage}{0.48\textwidth}
            \centering
            \compactsubfiginput{figures/S4toS3T1T2T3_graph.tex}\par
            \textbf{(a)} type $A$
        \end{minipage}\hfill
        \begin{minipage}{0.48\textwidth}
            \centering
            \compactsubfiginput{figures/S4toS3T1T2T3_graphn2.tex}\par
            \textbf{(b)} type $C$
        \end{minipage}
        \caption{The configurations in Lemma~\ref{S4toS3T1T2T3}.}
        \label{fig:s4-to-s3-t}
        \vspace{0.5em}
    \end{figure}
\end{proof}

\begin{lemma} \label{S1toS2}
    Let $m, n \in \{1,2\}$, $m \neq n$. The following two $H$-paths cannot both exist:
    \begin{enumerate}
        \item $P(v_1,v_2)$ with $v_1 \in V(S_4(s_2,a_4])$ and $v_2 \in V(S_m(s_{i_m},a_m))$,
        \item $P(v_3, v_4)$ with $v_3 \in V(S_m(s_{i_m},a_m))$ and $v_4 \in V(S_n(s_{i_n},a_n] \cup T_n(t,a_n])$.
    \end{enumerate}
    Note that $P(v_1,v_2)$ and $P(v_3,v_4)$ might intersect.
\end{lemma}

\begin{proof}
    If $P(v_3,v_4)$ meets $P(v_1,v_2)$, cut it at $P(v_1,v_2)$ from $v_4$; thus $P(v_3,v_4)$ meets $S_m \cup P(v_1,v_2)$ only in $v_3$, and $v_3 \in V(S_m(s_{i_m},a_m) \cup P(v_1,v_2))$.
    We assume that $v_4 \in V(S_n(s_{i_n},a_n])$; if $v_4 \in V(T_n(t,a_n])$, the constructions below work after replacing $S_n[v_4,a_n]$ by $T_n[v_4,a_n]$ in $S_K$ and $T_n$ by $T_n[t,v_4]$ in $T_K$.
    Let

    \[
        S' =
        \begin{cases}
            S_0 & \text{if } m = 1,\\
            \emptyset & \text{if } m = 2,
        \end{cases}
        \qquad
        U' =
        \begin{cases}
            \emptyset & \text{if } m = 1,\\
            S_0 & \text{if } m = 2,
        \end{cases}
    \]
    and
    \[
        \begin{aligned}
            S_K ={}& S_m \cup P[v_3,v_4] \cup P[v_1,v_2] \cup S_4[v_1,a_4] \cup S_n[v_4,a_n] \cup S_3 \cup S', \\
            T_K ={}& T_1 \cup T_2 \cup T_3, \\
            U_K ={}& U' \cup S_4[s_2,v_1].
        \end{aligned}
    \]
    Then, $K = S_K \cup T_K \cup U_K$ is of type $C$, as in Figure~\ref{fig:classical-type-c}(b), with $a'_n = a_n$. Indeed, in $S_K$ the hub of $H$ at which $S_m$ starts has degree $2$, so the leg of $S_K$ to $a_3$ is the path from $a_3$ along $S_3$, through $S'$ if $m=1$, and along $S_m$ to its first branch vertex; it contains $s_2$ as an internal vertex if $m = 1$ and $s_1$ if $m = 2$. Hence $U_K$ joins this internal vertex to $v_1$, which lies on the leg of $S_K$ to $a_4$ and is not its first vertex (possibly $v_1 = a_4$).
    
\end{proof}

\begin{lemma} \label{crossedS1toS2S3}
    Suppose that $H^*$ is crossed.
    There is no $H$-path $P(v_1,v_2)$ with $v_1 \in V(S_1(s_2,a_1))$ and $v_2 \in V(S_2 \cup S_3) \setminus \{s_1,a_2,a_3\}$.
    Furthermore, if some $H$-path $P(a,b)$ goes from $S_4(s_2,a_4]$ to $S_2(s_1,a_2)$, then no $H$-path goes from $S_3(s_1,a_3)$ to $a_1$.
\end{lemma}

\begin{proof}
    Suppose, for a contradiction, that the lemma fails. Then either an $H$-path $P(v_1,v_2)$ as in the first assertion exists, or the first assertion holds, $P(a,b)$ exists, and some $H$-path $Q = P(c,a_1)$ has $c \in V(S_3(s_1,a_3))$.

    First suppose that $P(v_1,v_2)$ exists. If $v_2 \in V(S_2)$, then $P(a_4,b_4)$ and $P(v_1,v_2)$ contradict Lemma~\ref{S1toS2}, with $m$ such that $b_4 \in V(S_m)$; hence $v_2 \in V(S_3(s_1,a_3))$. Moreover, $P(v_1,v_2)$ is disjoint from $P(s_0,v_0)$: otherwise, with $x$ the vertex of $P(v_1,v_2)$ on $P(s_0,v_0)$ closest to $v_2$, $P[v_2,x] \cup P[x,s_0]$ contradicts Lemma~\ref{S0min}. Let $c = v_2$, $Q = P[v_2,v_1] \cup S_1[v_1,a_1]$, let $z$ be the first vertex of $P(a_4,b_4)$ from $a_4$ on $Q \cup S_1 \cup S_2$ (it exists since $b_4 \in V(S_1 \cup S_2)$), and let
    \[
        T'_4 =
        \begin{cases}
            P[a_4,z] \cup S_1[z,v_1] & \text{if } z \in V(S_1[s_2,v_1)),\\
            P[a_4,z] & \text{otherwise},
        \end{cases}
        \qquad S'_4 = S_4.
    \]

    Now suppose that the first assertion holds, and that $P(a,b)$ and $Q = P(c,a_1)$ exist. If $b_4 \in V(S_2)$, we may take $P(a,b) = P(a_4,b_4)$, so we may assume that $a = a_4$ or $b_4 \in V(S_1)$. Let
    \[
        T'_4 = P[b,a] \cup S_4[a,a_4], \qquad
        S'_4 =
        \begin{cases}
            S_4 & \text{if } a = a_4,\\
            S_1[s_2,b_4] \cup P[b_4,a_4] & \text{otherwise}.
        \end{cases}
    \]
    By Lemma~\ref{S4toS3T1T2T3}, $P(a,b)$ is disjoint from $P(s_0,v_0)$; by the first assertion, $P(a,b)$ and $P(a_4,b_4)$ are internally disjoint if $a \neq a_4$; and by Lemmas~\ref{S4toS3T1T2T3} and~\ref{S0min}, $Q$ is internally disjoint from $P(a,b) \cup P(a_4,b_4) \cup P(s_0,v_0)$.

    In both cases, let
    \[
        S_K = S_2 \cup S_3 \cup Q \cup T'_4, \qquad
        T_K = T \cup P[s_0,v_0] \cup S_0[s_0,s_2] \cup S'_4.
    \]
    Then $K = S_K \cup T_K$ is of type $A$, as in Figure~\ref{fig:crossed-s1-to-s2-s3}(a), a contradiction; note that $S_K$ avoids $S_0[s_0,s_2]$ since $s_0 \neq s_1$, and that, by the above, $P(s_0,v_0)$ is disjoint from the $H$-paths in $S_K$.

    \begin{figure}[htbp]
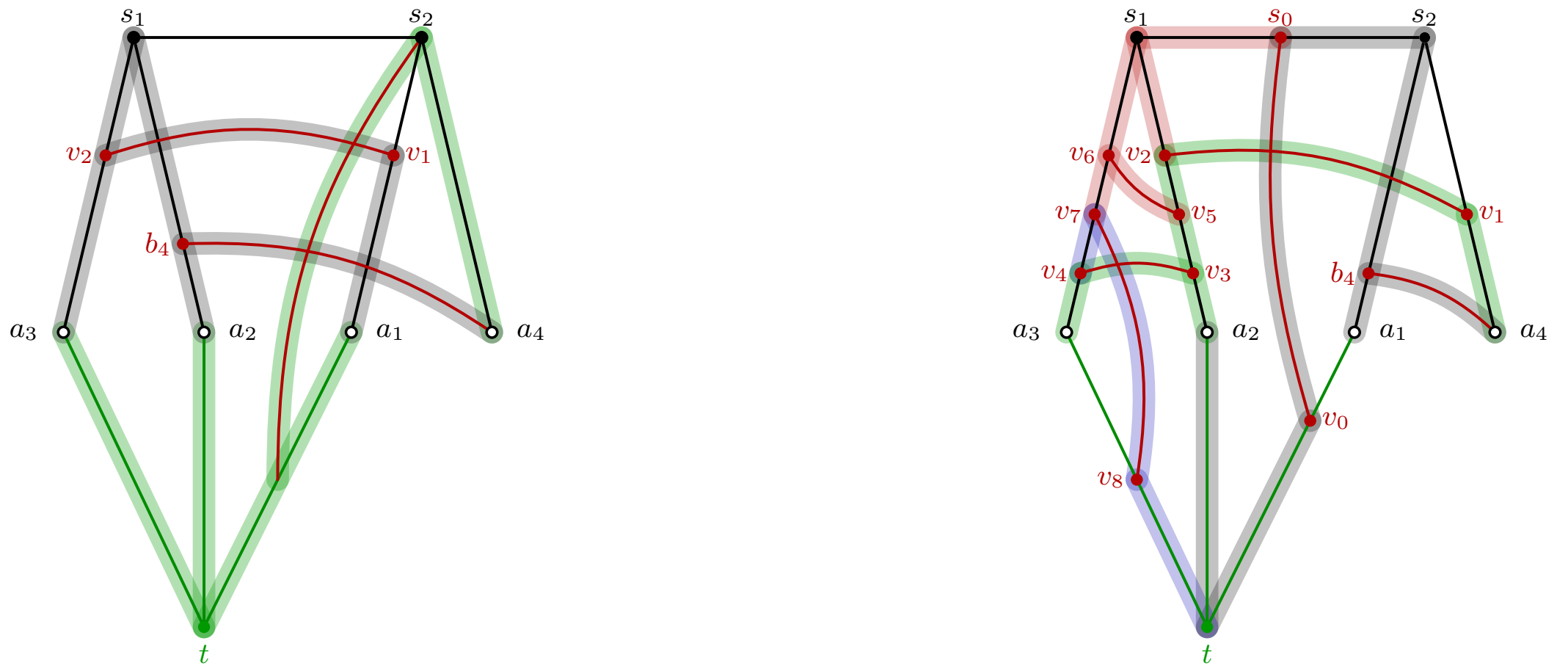

        \centering
        \begin{minipage}{0.48\textwidth}
            \centering
            \compactsubfiginput{figures/crossedS1toS2S3.tex}\par
            \textbf{(a)}
        \end{minipage}\hfill
        \begin{minipage}{0.48\textwidth}
            \centering
            \compactsubfiginput{figures/crossS3.tex}\par
            \textbf{(b)}
        \end{minipage}
        \caption{(a) The type-$A$ configuration in Lemma~\ref{crossedS1toS2S3}. (b) The type-$D$ configuration in Lemma~\ref{S3toT}.}
        \label{fig:crossed-s1-to-s2-s3}
    \end{figure}
\end{proof}

\begin{lemma} \label{S3toT}
    The following four $H$-paths cannot all exist:
    \begin{enumerate}
        \item $P(v_1,v_2)$ with $v_1 \in V(S_4(s_2,a_4])$ and $v_2 \in V(S_2(s_1,a_2))$;
        \item $P(v_3,v_4)$ with $v_3 \in V(S_2(s_1,a_2))$ and $v_4 \in V(S_3(s_1,a_3])$;
        \item $P(v_5,v_6)$ with $v_5 \in V(S_2(s_1,a_2))$ and $v_6 \in V(S_3(s_1,v_4))$;
        \item $P(v_7,v_8)$ with $v_7 \in V(S_3(s_1,v_4))$ and
              $v_8 \in V(T_1 \cup T_2 \cup T_3) \setminus \{a_1,a_2,a_3\}$.
    \end{enumerate}
    The four paths, $P(a_4,b_4)$ and $P(s_0,v_0)$ are pairwise internally disjoint, except that
    $P(v_5,v_6)$ may meet $P(v_3,v_4)$ and $P(v_7,v_8)$, $P(v_7,v_8)$ may
    meet $P(s_0,v_0)$, and $P(v_1,v_2)$ may meet $P(a_4,b_4)$.
\end{lemma}

\begin{proof}
    Suppose, for a contradiction, that the four paths exist.
    If $P(v_1,v_2)$ meets $P(a_4,b_4)$, let $y$ be the first vertex of $P(v_1,v_2)$ met when traversing $P(a_4,b_4)$ from $a_4$, and replace $P(v_1,v_2)$ by $P[a_4,y] \cup P[y,v_2]$.
    Otherwise, if $b_4 \in V(S_2(s_1,a_2))$, replace $P(v_1,v_2)$ by $P(a_4,b_4)$.
    In both cases the new path is contained in the union of the old one and $P(a_4,b_4)$, so the hypotheses still hold, and now either $v_1 = a_4$, or $b_4 \in V(S_1)$ and $P(v_1,v_2)$ is disjoint from $P(a_4,b_4)$.

    Let $z$ and $w$ be the vertices of $\{v_2,v_3\}$ closer to and farther from $s_1$ on $S_2$ (so $z = w$ if $v_2 = v_3$), and let
    \[
        S'_4 =
        \begin{cases}
            S_4 & \text{if } v_1 = a_4,\\
            S_1[s_2,b_4] \cup P[b_4,a_4] & \text{otherwise},
        \end{cases}
    \]
    \[
        T'_4 = P[v_2,v_1] \cup S_4[v_1,a_4], \qquad
        P'_0 = P[s_0,v_0] \cup T[v_0,a_2].
    \]
    Let $U_2$ be the subpath of $P(v_7,v_8)$ from $v_7$ to its first vertex on $P'_0$, or, if $P(v_7,v_8)$ does not meet $P'_0$, the path $P[v_7,v_8]$ followed by the path of $T$ from $v_8$ to its first vertex on $T[v_0,a_2]$; let $u_2$ be its last vertex. Let
    \[
        \begin{aligned}
            S_K ={}& S_1 \cup S'_4 \cup S_0[s_2,s_0] \cup P'_0, \\
            T_K ={}& S_2[z,a_2] \cup T'_4 \cup P[v_3,v_4] \cup S_3[v_4,a_3], \\
            U_K ={}& U_2 \cup S_3[v_7,v_4], \\
            V_1 ={}& S_0[s_0,s_1] \cup S_2[s_1,z], \\
            V_2 ={}& S_0[s_0,s_1] \cup S_3[s_1,v_7],
        \end{aligned}
    \]
    and let $V_3 = X_3 \cup P[v_5,v_6] \cup Y_3$, where
    \[
        X_3 =
        \begin{cases}
            S_2[z,v_5] & \text{if } v_5 \in V(S_2(s_1,z)),\\
            \emptyset & \text{otherwise},
        \end{cases}
    \]
    \[
        Y_3 =
        \begin{cases}
            S_3[v_6,v_7] & \text{if } v_6 \in V(S_3(s_1,v_7)),\\
            \emptyset & \text{otherwise}.
        \end{cases}
    \]
    Then $K = S_K \cup T_K \cup U_K \cup V_1 \cup V_2 \cup V_3$ is of type $D$, as in Figure~\ref{fig:crossed-s1-to-s2-s3}(b), a contradiction; note that $s_0$ precedes $u_2$ on $P'_0$ by the choice of $U_2$. Here $V_3$ is given in the sense of the convention following the definition of the types: by the hypotheses on the four paths, $P(v_5,v_6)$ meets $S_K \cup T_K \cup U_K$ only in its ends, in interior vertices of $P(v_3,v_4)$, which lie on the leg of $T_K$ to $a_3$ before $u_3 = v_4$, and in vertices of $U_2 - u_2$; moreover, $V_1 \cap V_3 \subseteq X_3$ and $V_2 \cap V_3 \subseteq Y_3$ are common subpaths ending at $z$ and at $v_7$, respectively.

\end{proof}

In a reduced frame, the labels $a_1$ and $a_2$ may be interchanged, so the lemma also holds with $S_1$ and $T_1$ in place of $S_2$ and $T_2$.

\begin{lemma} \label{crossed}
    Suppose that $H^*$ is crossed.
    There is no $H$-path $P(a,b)$ with $a \in V(S_4(s_2,a_4])$ and $b \in V(S_2(s_1,a_2))$.
\end{lemma}

\begin{proof}
    Suppose, for a contradiction, that the lemma fails. Let $\mathcal{H}$ be the set of augmented frames obtained from $H^*$ by replacing $S_3$ and $S_4$ so that $P(s_0,v_0)$ and $P(a_4,b_4)$ remain $H$-paths; they have the same $S_0$ as $H^*$, hence are minimal. Among the members of $\mathcal{H}$ for which the lemma fails, choose one, still denoted $H^*$, and an $H$-path $P(a,b)$ with $a \in V(S_4(s_2,a_4])$ and $b \in V(S_2(s_1,a_2))$, such that:
    \begin{itemize}
        \item[(M1)] $|E(S_4[a,a_4])|$ is minimal;
        \item[(M2)] subject to (M1), if some $H$-path goes from $S_2(s_1,a_2)$ to $S_3(s_1,a_3]$, then such a path $P(v_1,v_2)$ is chosen with $|E(S_3[v_2,a_3])|$ minimal;
        \item[(M3)] subject to (M1) and (M2), if some $H$-path goes from $S_2(s_1,a_2)$ to $S_4(s_2,a)$, then such a path $P(v_3,v_4)$ is chosen with $|E(S_4[s_2,v_4])|$ minimal.
    \end{itemize}

    Let $T'_4 = P[b,a] \cup S_4[a,a_4]$. If $a = a_4$, let $S'_4 = S_4$; otherwise $b_4 \in V(S_1)$ by (M1), and we let $S'_4 = S_1[s_2,b_4] \cup P[b_4,a_4]$. By Lemmas~\ref{S4toS3T1T2T3} and~\ref{crossedS1toS2S3}, $P(a,b)$ is disjoint from $P(s_0,v_0)$ and, if $a \neq a_4$, from $P(a_4,b_4)$. Hence $T'_4$ and $S'_4$ meet only in $a_4$ and are disjoint from $P(s_0,v_0)$. We also write $P'_0 = P[s_0,v_0] \cup T[v_0,a_2]$.

    \emph{(C1) Suppose that $s_0 \neq s_2$, or that some $H$-path $P(v_2,v_3)$ goes from $S_4(s_2,a_4]$ to $S_1(s_2,a_1)$. Then there is no $H$-path $P(v_1,a_3)$ with $v_1 \in V(S_1(s_2,a_1))$.}
    Suppose otherwise. By Lemmas~\ref{S4toS3T1T2T3} and~\ref{S0min}, $P(v_1,a_3)$ is internally disjoint from $P(a,b) \cup P(a_4,b_4) \cup P(s_0,v_0)$; by Lemmas~\ref{S4toS3T1T2T3} and~\ref{crossedS1toS2S3}, $P(v_2,v_3)$ is internally disjoint from $P(v_1,a_3)$, $P(s_0,v_0)$ and $T'_4$.
    Let $U_1$ be the subpath of $Q \cup P[v_1,a_3]$ from its last vertex on $S'_4$ to $a_3$, where
    \[
        Q =
        \begin{cases}
            P[v_2,v_3] \cup S_1[v_3,v_1] & \text{if } s_0 = s_2 \text{ and } a = a_4,\\
            S_1[s_2,v_1] & \text{otherwise}.
        \end{cases}
    \]
    Then
    \[
        S_K = T \cup P[s_0,v_0] \cup S_0[s_0,s_2] \cup S'_4, \qquad T_K = S_2 \cup S_3 \cup T'_4, \qquad U_K = U_1 \cup S_0[s_0,s_1]
    \]
    form a subgraph of type $B$ (Figure~\ref{fig:crossed-configurations}(a)), with $U_1$ starting after $s_0$ on the leg of $S_K$ to $a_4$, and $U_2 = S_0[s_0,s_1]$. Indeed, this leg is $P[v_0,s_0] \cup S_0[s_0,s_2] \cup S'_4$, and $Q$ meets it in $S'_4$, which starts at $s_2$: in the first case of $Q$, $Q$ starts at $v_2 \in V(S_4(s_2,a_4])$ and $S'_4 = S_4$; in the second, $Q$ starts at $s_2$, and if $s_0 = s_2$ then $a \neq a_4$, $S'_4 = S_1[s_2,b_4] \cup P[b_4,a_4]$, and the last vertex of $Q$ on $S'_4$ is $v_1$ or $b_4$, both different from $s_2$.

    \smallskip
    \emph{Step 1.} \emph{(C2) Every $H$-path from $S_2(s_1,a_2)$ to $V(H) \setminus V(S_2)$ ends in $V(S_3(s_1,a_3] \cup S_4(s_2,a]) \cup (V(T) \setminus \{a_1,a_2\})$,} by Lemma~\ref{S0min}, Lemma~\ref{S1toS2} (with $m = 2$, applied to $P(a,b)$), and (M1).

    \emph{Step 2.} \emph{(C3) There is no $H$-path $P(v_1,v_2)$ with $v_1 \in V(S_2(s_1,a_2))$ and $v_2 \in V(T) \setminus \{a_1,a_2,a_3\}$.}
    By Lemma~\ref{S4toS3T1T2T3}, $P(v_1,v_2)$ is internally disjoint from $P(a,b) \cup P(a_4,b_4)$; if it meets $P(s_0,v_0)$, we cut it there from $v_1$. Then
    \[
        \begin{aligned}
            S_K ={}& T \cup P[s_0,v_0] \cup S_0[s_0,s_2] \cup S'_4, \\
            T_K ={}& S_2 \cup S_3 \cup T'_4, \\
            U_K ={}& S_0[s_0,s_1] \cup P(v_1,v_2)
        \end{aligned}
    \]
    form a subgraph of type $B$ (Figure~\ref{fig:crossed-configurations}(b)).
    In particular, an $H$-path with an end in $V(T) \setminus \{a_1,a_2,a_3\}$ is internally disjoint from every $H$-path with an end in $V(S_2(s_1,a_2))$.

    \emph{(C4) Suppose that some $H$-path $P(v_3,v_4)$ goes from $S_2(s_1,a_2)$ to $S_4(s_2,a_4)$. Then there are no disjoint $H$-paths $P(c,a_1)$ and $P(a_4,d)$ with $c \in V(S_4(s_2,a_4))$ and $d \in V(S_1(s_2,a_1))$.}
    Suppose otherwise. By (C2), (C3) and Lemmas~\ref{S4toS3T1T2T3} and~\ref{crossedS1toS2S3}, the paths $P(c,a_1)$, $P(a_4,d)$, $P(v_3,v_4)$ and $P(s_0,v_0)$ are pairwise internally disjoint. Let $x$ be the vertex of $\{v_4,c\}$ closer to $s_2$ on $S_4$. Then
    \[
        \begin{aligned}
            S_K ={}& S_2 \cup S_3 \cup P[v_3,v_4] \cup S_4[x,a_4] \cup P[c,a_1], \\
            T_K ={}& T \cup P[s_0,v_0] \cup S_0[s_0,s_2] \cup S_1[s_2,d] \cup P[d,a_4]
        \end{aligned}
    \]
    form a subgraph of type $A$.

    \medskip
    \emph{Step 3.} \emph{There is no $H$-path from $S_2(s_1,a_2)$ to $S_3(s_1,a_3]$.}
    Suppose otherwise, and let $P(v_1,v_2)$ be as in (M2). By Lemma~\ref{S4toS3T1T2T3} and (C3), it is internally disjoint from $P(a,b)$, $P(a_4,b_4)$ and $P(s_0,v_0)$.

    \emph{(C5) There is no $H$-path $C = P(c,c')$ with $c \in V(S_3[s_1,v_2))$ and $c' \in V(S_3(v_2,a_3])$.}
    Such a path is internally disjoint from $P(v_1,v_2)$ by (M2), and from $P(a,b) \cup P(a_4,b_4) \cup P(s_0,v_0)$ by Lemmas~\ref{S4toS3T1T2T3} and~\ref{S0min} applied at $c'$. Replacing $S_3$ by $S_3[s_1,c] \cup C \cup S_3[c',a_3]$ gives a member of $\mathcal{H}$ in which $P(a,b)$ is still an $H$-path, so (M1) still holds, and in which $P(v_1,v_2) \cup S_3[v_2,c']$ ends closer to $a_3$, contrary to (M2).

    \emph{Step 3.1.} \emph{There is no $H$-path from $S_2(s_1,a_2)$ to $S_4(s_2,a)$.}
    Suppose that some $H$-path $P(v_3,v_4)$ goes from $S_2(s_1,a_2)$ to $S_4(s_2,a)$. By Lemmas~\ref{S4toS3T1T2T3} and~\ref{crossedS1toS2S3} and (C3), it is internally disjoint from $P(v_1,v_2)$, $P(s_0,v_0)$ and $S'_4$.
    Let $x$ be the vertex of $\{v_1,b\}$ closer to $s_1$ on $S_2$, and let
    \[
        T_K = S_2[x,a_2] \cup P[v_1,v_2] \cup S_3[v_2,a_3] \cup T'_4, \qquad \Lambda = S_1[a_1,s_2] \cup S_0[s_2,s_1] \cup S_3[s_1,v_2].
    \]

    \emph{(C6) There is no $H$-path $P(v_5,v_6)$ with $v_5 \in V(S_0 \cup S_3[s_1,v_2)) \setminus \{s_0\}$ and $v_6 \in V(T) \setminus \{a_1,a_2,a_3\}$. If, moreover, $s_0 \neq s_2$ or some $H$-path $P(v_7,v_8)$ goes from $S_4(s_2,a_4)$ to $S_1(s_2,a_1)$, then the same holds with $v_5 \in V(S_1(s_2,a_1))$.}
    By (C3) and Lemma~\ref{S4toS3T1T2T3}, $P(v_5,v_6)$ is internally disjoint from $P(a,b)$, $P(a_4,b_4)$, $P(v_1,v_2)$ and $P(v_3,v_4)$.

    Suppose that $P(v_5,v_6)$ exists, and, if $v_5 \in V(S_1(s_2,a_1))$ and $s_0 = s_2$, let $P(v_7,v_8)$ be as above; only in this case do we say that $P(v_7,v_8)$ is given. By Lemmas~\ref{S4toS3T1T2T3} and~\ref{crossedS1toS2S3} and (C3), $P(v_7,v_8)$, if given, is internally disjoint from $T'_4$, $P(v_1,v_2)$, $P(v_3,v_4)$, $P(v_5,v_6)$ and $P(s_0,v_0)$. Both $s_0$ and $v_5$ lie on the path $\Lambda$, which goes from $a_1$ to $v_2$; let $p$ be the one closer to $a_1$ on $\Lambda$ and $q$ the other, and let $P(p,v_p)$ and $P(q,v_q)$ be the corresponding paths among $P(s_0,v_0)$ and $P(v_5,v_6)$. Let
    \[
        \widehat{S}_4 =
        \begin{cases}
            P[v_8,v_7] \cup S_4[v_7,a_4] & \text{if } s_0 = s_2,\ a = a_4 \text{ and } P(v_7,v_8) \text{ is given},\\
            S'_4 & \text{otherwise},
        \end{cases}
    \]
    let $S_K$ be the minimal subtree of $S_1 \cup \Lambda[a_1,p] \cup \widehat{S}_4 \cup P[p,v_p] \cup T[v_p,a_2]$ containing $a_1,a_2,a_4$, and let $W$ be $P[q,v_q] \cup T[v_q,a_2]$ cut at its first vertex on $P[p,v_p] \cup T[v_p,a_2]$. Then
    \[
        U_K = W \cup \Lambda[q,v_2], \qquad V_K = V_1 \cup V_2 \cup V_3,
    \]
    where, with $D = S_1 \cup S_2 \cup S_4 \cup P[v_3,v_4]$, together with $P[v_7,v_8]$ if it is given, $V_1$ is a shortest path in $D - V(U_K) - \{a_1,a_2,a_4\}$ from $V(S_K) \setminus \{a_1,a_2,a_4\}$ to $V(T_K) \setminus \{a_2,a_3,a_4\}$; it is nonempty since $S_K$ and $T_K$ meet only in roots, $V_2$ is the subpath of $\Lambda[p,q]$ from its last vertex on $S_K$ to $q$, and $V_3 = S_2[x,s_1] \cup \Lambda[s_1,q]$ (if $q$ precedes $s_1$ on $\Lambda$, then $\Lambda[s_1,q] \subseteq U_K$, and $V_3$ is shortened to $S_2[x,s_1]$ by the convention following the definition of the types), form with $T_K$ a subgraph of type $D$ (Figure~\ref{fig:crossed-configurations}(c)). The path $V_1$ exists. Its admissible ends lie in $V(S_K) \cap V(D)$, which is contained in $S_1 \cup \widehat{S}_4$, so on the legs of $S_K$ at $a_1$ and $a_4$ and not beyond the attachment of $U_K$ on the leg at $a_2$; and its other end lies in $V(T_K) \cap V(D)$, which is contained in $V(S_2[x,a_2] \cup T'_4)$, since $P[v_3,v_4]$ may meet $T'_4$ in interior vertices of $P(a,b)$. Thus this end lies on the legs of $T_K$ at $a_2$ or $a_4$, or on the portion of its leg at $a_3$ strictly before the attachment $u_3=v_2$; in each case it is an admissible end for a connector of type $D$.

    Let $y$ be the neighbor of $a_1$ on $S_1$. As $a_1$ is a leaf of $S_K$ whose only neighbor in the graph defining $S_K$ is $y$, we have $y \in V(S_K) \setminus \{a_1,a_2,a_4\}$. Let $Q_S$ be the subpath of $S_2$ from $v_3$ to its first vertex on $S_2[x,a_2]$; it avoids $a_2$, and $S_2(s_1,x)$ meets none of $S_K$, $T_K$, $U_K$.

    If $s_2 \notin V(U_K)$, consider the walk $S_1[y,s_2] \cup S_4[s_2,v_4] \cup P[v_4,v_3] \cup Q_S$. It avoids $a_1$, $a_2$ and $a_4$, since $v_4 \neq a_4$. It also avoids $U_K$: since $q$ occurs no earlier than $s_0$ on $\Lambda$, it occurs at or after $s_2$. If $q=s_2$, then $s_2\in V(U_K)$, contrary to the present assumption. Thus $q$ occurs strictly after $s_2$, and $\Lambda[q,v_2]$ avoids $S_1\cup S_4$. Moreover, $W$ consists of an $H$-path internally disjoint from $P(v_3,v_4)$ followed by a path in $T$.

    Otherwise $q=s_0=s_2$, so $v_5\in V(S_1(s_2,a_1))$ and $P(v_7,v_8)$ is given. Since $v_5\in V(S_1(s_2,a_1))$, the neighbor $y$ of $a_1$ on $S_1$ is different from $s_2$. Consider then the walk $S_1[y,v_8] \cup P[v_8,v_7] \cup S_4[v_7,v_4] \cup P[v_4,v_3] \cup Q_S$, where $S_4[v_7,v_4]$ is the subpath of $S_4$ between $v_7$ and $v_4$. Since $y,v_8\in V(S_1(s_2,a_1))$ and $v_7\neq a_4$, it avoids $s_2$ and the roots, and it avoids the rest of $U_K$ as before. In both cases the walk lies in $D-V(U_K)-\{a_1,a_2,a_4\}$, starts on $S_K$ and ends on $T_K$, so it contains a path from $V(S_K)\setminus\{a_1,a_2,a_4\}$ to $V(T_K)\setminus\{a_2,a_3,a_4\}$ in this graph.

    We may moreover choose $V_1$ so that it meets $V_2$ at most in a common initial segment and $V_3$ at most in a common final segment ending at $x$. Let $Q_2$ be the subpath $V_2\cap S_1$ (possibly empty), and let $Q_3=S_2[x,s_1]$. Apart from $s_1$, $V_2\cup V_3$ meets $D$ only in $Q_2\cup Q_3$: the remaining parts lie in $S_0\cup S_3$. If $s_1\notin V(U_K)$, then $s_1$ has degree one in $D$ and belongs to neither $S_K$ nor $T_K$, so a simple path between these two trees in $D$ avoids $s_1$; if $s_1\in V(U_K)$, it is deleted already. The paths $Q_2$ and $Q_3$ are disjoint, and meet $S_K\cup T_K$ only at the first vertex of $Q_2$, when nonempty, and at $x$, respectively. If $V_1$ meets $Q_3$, replace its final part from its first vertex $w'$ on $Q_3$ by $Q_3[w',x]$. Then, if it meets $Q_2$, replace its initial part up to its last vertex $w$ on $Q_2$ by the subpath of $Q_2$ from its first vertex to $w$. These replacements stay outside $U_K$: any vertex of $Q_2$ or $Q_3$ in $U_K$ is at its end opposite the relevant tree. Since $Q_2$ and $Q_3$ are disjoint, the second replacement preserves the first. The resulting connectors have the required intersections. This proves (C6).

    We now show that some $H$-path $P(v_5,v_6)$ goes from $S_1(s_2,a_1)$ to $V(T) \setminus \{a_1,a_2,a_3\}$ and that the hypothesis of the second part of (C6) holds; this contradicts (C6).
    Let $Z = \{v_2,a_1,a_2,s_0\}$, $X = V(S_0 \cup S_2 \cup S_3[s_1,v_2) \cup S_4) \setminus Z$ and $Y = V(S_3(v_2,a_3] \cup T) \setminus Z$.
    By (C2), (C3), (M2), (C5), Lemmas~\ref{S0min} and~\ref{S4toS3T1T2T3} and the first part of (C6), no $H$-path goes from $X$ to $Y$, and by Lemma~\ref{crossedS1toS2S3}, every $H$-path from $S_1(s_2,a_1)$ to $Y$ ends in $V(T) \setminus \{a_1,a_2\}$. If $s_0 \neq s_2$, then $S_1(s_2,a_1)$ is attached to $X$ at $s_2$, and since $Z$ does not separate $X \cup V(S_1(s_2,a_1))$ from $Y$, there is such a path $P(v_5,v_6)$. If $s_0 = s_2$, then $S_1(s_2,a_1)$ is attached to the rest of $H$ only at $s_2,a_1 \in Z$, and by Lemmas~\ref{S0min} and~\ref{crossedS1toS2S3} every $H$-path from $X$ to $S_1(s_2,a_1)$ starts in $S_4(s_2,a_4]$; since $Z$ separates neither $X$ from $Y \cup V(S_1(s_2,a_1))$ nor $X \cup V(S_1(s_2,a_1))$ from $Y$, there are both such a path $P(v_7,v_8)$ and a path $P(v_5,v_6)$ as above. In both cases $v_6 \neq a_3$ by (C1).

    \emph{(C7) If $s_0 = s_2$, some $H$-path goes from $S_4(s_2,a_4)$ to $S_1(s_2,a_1)$.}
    Suppose not, so that $P(v_7,v_8) = P(a_4,v_8)$. Let $Z' = \{v_2,a_2,s_2,a_4\}$ and $Y' = V(S_1 \cup S_3(v_2,a_3] \cup T) \setminus Z'$. By (C2), (C3), (M2) and Lemmas~\ref{S1toS2} and~\ref{crossedS1toS2S3}, no $H$-path goes from $S_2(s_1,a_2)$ to $Y'$; by (C5), the first part of (C6) and Lemmas~\ref{S0min} and~\ref{crossedS1toS2S3}, none goes from $V(S_0 \cup S_3[s_1,v_2)) \setminus \{s_2\}$ to $Y'$; and by Lemma~\ref{S4toS3T1T2T3} and our assumption, every $H$-path from $S_4(s_2,a_4)$ to $Y'$ ends at $a_1$. As $Z'$ is not a separator, there is such a path $P(c,a_1)$.
    If $P(c,a_1)$ meets $P(a_4,v_8)$, following it from $c$ to its first vertex on $P(a_4,v_8)$ and then $P(a_4,v_8)$ to $v_8$ gives an $H$-path from $c$ to $S_1(s_2,a_1)$, contrary to our assumption. Otherwise $P(c,a_1)$ and $P(a_4,v_8)$ contradict (C4). This proves (C7), and hence Step~3.1.

    \emph{Step 3.2.} \emph{Conclusion of Step 3.}
    By (C2), (C3), Step~3.1 and (M2), every $H$-path from $S_2(s_1,a_2)$ leaving $S_2$ ends in $S_3(s_1,v_2] \cup \{a\}$, so, as $\{s_1,v_2,a_2,a\}$ is not a separator, there is an $H$-path $P(v_3,v_4)$ from $S_2(s_1,a_2)$ to $S_3(s_1,v_2)$, and then an $H$-path $P(v_5,v_6)$ from $S_3(s_1,v_2)$ to $V(H) \setminus (V(S_2 \cup S_3[s_1,v_2]) \cup \{a\})$.
    By Lemmas~\ref{S0min}, \ref{crossedS1toS2S3}, \ref{S4toS3T1T2T3} and (C5), $v_6 \in V(T) \setminus \{a_1,a_2,a_3\}$.
    By Lemma~\ref{S4toS3T1T2T3} and (C3), the paths $P(a,b)$, $P(v_1,v_2)$, $P(v_3,v_4)$, $P(v_5,v_6)$ satisfy the hypotheses of Lemma~\ref{S3toT}, a contradiction. This proves Step~3.

    \medskip
    \emph{Step 4.} \emph{Conclusion.}
    By (C2), (C3) and Step~3, every $H$-path leaving $S_2(s_1,a_2)$ ends in $S_4(s_2,a]$. As $\{s_1,a_2,a\}$ is not a separator, some path of $G - \{s_1,a_2,a\}$ leaves $S_2(s_1,a_2) \cup P(a,b)$; preceded if necessary by a subpath of $P(a,b)$ starting at $b$, it gives an $H$-path from $S_2(s_1,a_2)$ to $S_4(s_2,a)$. Let $P(v_3,v_4)$ be as in (M3); by Lemmas~\ref{S4toS3T1T2T3} and~\ref{crossedS1toS2S3}, it is internally disjoint from $P(s_0,v_0)$ and $S'_4$.

    Let $Z = \{s_1,s_2,a_1,a_2\}$. By (C2), (C3), Step~3 and Lemmas~\ref{S0min} and~\ref{S4toS3T1T2T3}, the only $H$-paths from $S_2(s_1,a_2) \cup S_4(s_2,a_4]$ to $V(H) \setminus (V(S_2 \cup S_4) \cup Z)$ go from $S_4$ to $S_1(s_2,a_1)$; since $Z$ is not a separator, such a path exists, so the hypothesis of (C1) holds. Let $W$ consist of $V(S_1(s_2,a_1) \cup S_2(s_1,a_2) \cup S_4(s_2,a_4])$ and the internal vertices of all $H$-paths with both ends in this set. Since $Z$ does not separate $W$ from $V(T) \setminus \{a_1,a_2,a_3\}$, a shortest path $R$ of $G - Z$ from $W$ to $V(H) \setminus (W \cup Z)$ has its internal vertices outside $W \cup V(H)$; if its first vertex is an internal vertex of an $H$-path with both ends in $W$, we precede $R$ by the subpath of that $H$-path from one of its ends to this vertex. This gives an $H$-path from $V(S_1(s_2,a_1) \cup S_2(s_1,a_2) \cup S_4(s_2,a_4])$ to $V(H) \setminus (W \cup Z) = V(S_0 \cup S_3 \cup T) \setminus Z$. By the above it starts in $S_1(s_2,a_1)$; by Lemmas~\ref{S0min} and~\ref{crossedS1toS2S3} it ends in $V(T) \setminus \{a_1,a_2\}$, and by (C1) not at $a_3$. This is the required $H$-path $P(v_5,v_6)$.

    \emph{(C8) There is an $H$-path $P(v_7,v_8)$ with $v_7 \in V(S_4(v_4,a_4))$ and $v_8 \in V(S_1(s_2,a_1))$.}
    We use two separations. In both, $X$ consists of vertices of $S_2(s_1,a_2) \cup S_4(v_4,a_4]$ and of internal vertices of some of the paths $P(a,b)$, $P(a_4,b_4)$, $P(v_3,v_4)$; in (a) and (b) below, a path is included only if it has an end in $S_2$, which holds for $P(a,b)$ and $P(v_3,v_4)$, and for $P(a_4,b_4)$ exactly when $b_4 \in V(S_2)$, and $Z$ is a set of four vertices with $N_{H \cup P(a,b) \cup P(a_4,b_4) \cup P(v_3,v_4)}(X) \subseteq Z$. As $Z$ is not a separator, $G - Z$ has a path from $X$ to $V(H) \setminus (X \cup Z)$; let $R$ be a shortest one, from $x_0 \in X$ to $c' \in V(H) \setminus (X \cup Z)$. Its internal vertices lie outside $X \cup V(H)$. If $x_0 \in V(H)$, let $c = x_0$ and $P(c,c') = R$; otherwise $x_0$ is an internal vertex of one of the three paths above, and $P(c,c')$ is $R$ preceded by the subpath of that path from its end $c$ in $S_2$ to $x_0$, which is internally disjoint from $H$ and meets $R$ only in $x_0$. In both cases $P(c,c')$ is an $H$-path with $c \in V(S_2(s_1,a_2) \cup S_4(v_4,a_4])$ and $c' \notin Z$. Two outcomes are excluded in both cases:
    \begin{itemize}
        \item[(i)] $c \in V(S_2)$: then $c' \in V(S_4(s_2,v_4))$ by (C2), (C3) and Step~3, contrary to (M3);
        \item[(ii)] $c \in V(S_4(v_4,a_4])$, $c' \in V(S_4[s_2,v_4))$ and $P(c,c')$ is internally disjoint from $P(a_4,b_4)$: by Lemma~\ref{S4toS3T1T2T3}, it is also disjoint from $P(s_0,v_0)$ (otherwise, following it from $c$ to its first vertex on $P(s_0,v_0)$ and then $P(s_0,v_0)$ to $v_0$ gives an $H$-path excluded by that lemma), so replacing $S_4$ by $S_4[s_2,c'] \cup P[c',c] \cup S_4[c,a_4]$ gives a member of $\mathcal{H}$: the new leg meets the rest of $H$ only in $s_2$ and $a_4$, since $P(c,c')$ is an $H$-path, and $P(s_0,v_0)$ and $P(a_4,b_4)$ remain $H$-paths, since they avoid the interior of $P(c,c')$ and have no end in the removed subpath $S_4(c',c)$. If $a \in V(S_4(c',c))$, then $P[b,a] \cup S_4[a,c]$ is an $H$-path of the new frame contradicting (M1). Otherwise $P(a,b)$ is still an $H$-path, (M1) has the same value, and $P(v_3,v_4) \cup S_4[v_4,c']$ is an $H$-path of the new frame from $S_2(s_1,a_2)$ to $S_4[s_2,c']$; this contradicts Lemma~\ref{S0min} if $c' = s_2$. Otherwise, the new frame has no $H$-path $R'$ from $S_2(s_1,a_2)$ to $S_3(s_1,a_3]$. Indeed, $R'$ is internally disjoint from the new frame, hence from $P(c,c')$; if it meets the removed subpath $S_4(c',c)$, its final segment from its last vertex there is an $H$-path of $H^*$ from $S_4$ to $S_3(s_1,a_3]$, excluded by Lemma~\ref{S4toS3T1T2T3}, and otherwise $R'$ itself is an $H$-path of $H^*$ from $S_2(s_1,a_2)$ to $S_3(s_1,a_3]$, excluded by Step~3. So (M2) is unchanged, and the new frame contradicts (M3).
    \end{itemize}
    Moreover, if $c \in V(S_4(v_4,a_4])$, then $c' \in V(S_1(s_2,a_1]) \cup V(S_4[s_2,v_4))$ by Lemmas~\ref{S4toS3T1T2T3} and~\ref{S0min}.

    \emph{(a) Some $H$-path goes from $S_4(v_4,a_4]$ to $S_1(s_2,a_1)$.} If $P(a_4,b_4)$ ends in $S_1$, it is such a path. Otherwise $b_4 \in V(S_2)$, so $a = a_4$ by (M1). Let $X$ consist of $V(S_2(s_1,a_2) \cup S_4(v_4,a_4])$ and the internal vertices of $P(a,b)$, $P(a_4,b_4)$ and $P(v_3,v_4)$; all these paths end in $X \cup \{s_1,a_2,v_4\}$. As $\{s_1,a_2,v_4,a_1\}$ is not a separator, we obtain $P(c,c')$ as above. By (i), $c \in V(S_4(v_4,a_4])$, so $P(c,c')$ is internally disjoint from $X$, in particular from $P(a_4,b_4)$; by (ii), and since $a_1$ lies in the separator, $c' \in V(S_1(s_2,a_1))$.

    \emph{(b) Some such path starts in $S_4(v_4,a_4)$.} Suppose that every such path starts at $a_4$, and let $P(a_4,v_8)$ be one of them. Let $X$ consist of $V(S_2(s_1,a_2) \cup S_4(v_4,a_4))$ and the internal vertices of $P(a,b)$, of $P(v_3,v_4)$ and, if it ends in $S_2$, of $P(a_4,b_4)$; all these paths end in $X \cup \{s_1,a_2,v_4,a_4\}$. As $\{s_1,a_2,v_4,a_4\}$ is not a separator, we obtain $P(c,c')$ as above. By (i), $c \in V(S_4(v_4,a_4))$. If $P(c,c')$ meets $P(a_4,b_4)$, then $P(a_4,b_4)$ ends in $S_1$, and following $P(c,c')$ from $c$ to its first vertex on $P(a_4,b_4)$ and then $P(a_4,b_4)$ to $b_4$ gives an $H$-path from $c$ to $S_1(s_2,a_1)$, contrary to our assumption. So, by (ii), $c' \in V(S_1(s_2,a_1])$; since $c \neq a_4$, the path $P(c,c')$ does not start at $a_4$, so by our assumption $c' \notin V(S_1(s_2,a_1))$, that is, $c' = a_1$. If $P(c,a_1)$ meets $P(a_4,v_8)$, following it from $c$ to its first vertex on $P(a_4,v_8)$ and then $P(a_4,v_8)$ to $v_8$ gives an $H$-path from $c$ to $S_1(s_2,a_1)$, contrary to our assumption. Otherwise $P(c,a_1)$ and $P(a_4,v_8)$ contradict (C4). This proves (C8).

    By Lemmas~\ref{S4toS3T1T2T3} and~\ref{crossedS1toS2S3} and (C3), the paths $P(v_5,v_6)$, $P(v_7,v_8)$, $P(s_0,v_0)$, $T'_4$ and $P(v_3,v_4)$ are pairwise internally disjoint, except that $P(v_3,v_4)$ may meet $T'_4$ and $P(s_0,v_0)$ may meet $P(v_5,v_6)$; moreover $P(v_5,v_6)$ is internally disjoint from $P(a_4,b_4)$. Let $T_K = S_2 \cup S_3 \cup T'_4$, let $P_3$ be $P[v_4,v_3]$ cut at its first vertex on $T_K$, let
    \[
        \widehat{S}_4 =
        \begin{cases}
            P[v_8,v_7] \cup S_4[v_7,a_4] & \text{if } a = a_4,\\
            S'_4 & \text{otherwise},
        \end{cases}
    \]
    let $S_K$ be the minimal subtree of $S_1 \cup \widehat{S}_4 \cup P[v_5,v_6] \cup T[v_6,a_2]$ containing $a_1,a_2,a_4$, and let $W_0$ be $P'_0$ cut at its first vertex on $P[v_5,v_6] \cup T[v_6,a_2]$. Let $V_2$ be the path from $s_0$ to $S_K$ in $S_0[s_0,s_2] \cup S_1$, let $V_3 = S_0[s_0,s_2] \cup S_4[s_2,v_4] \cup P_3$, and let $V_1$ be a shortest path from $S_K$ to $T_K$ in $(S_1 - s_2) \cup S_4[v_4,a_4) \cup P[v_7,v_8] \cup P_3$; it exists since $v_7 \in V(S_4(v_4,a_4))$. Then
    \[
        U_K = S_0[s_0,s_1] \cup W_0, \qquad V_K = V_1 \cup V_2 \cup V_3
    \]
    form with $S_K$ and $T_K$ a subgraph of type $D$ (Figure~\ref{fig:crossed-configurations}(d)); we choose $V_1$, as in (C6), so that it meets $V_2$ at most in a common initial segment and $V_3$ at most in a common final segment: $V_2$ and $V_3$ meet the graph in which $V_1$ is chosen only in the subpath of $V_2$ on $S_1 - s_2$, which starts on $S_K$, and in $P_3$, which ends on $T_K$ (together with its first vertex $v_4$), and the same two replacements apply. This contradiction completes the proof.

    \begin{figure}[htbp]
        \centering
        \begin{minipage}{0.48\textwidth}
            \centering
            \compactsubfiginput{figures/crossed4.tex}\par
            \textbf{(a)} type $B$
        \end{minipage}\hfill
        \begin{minipage}{0.48\textwidth}
            \centering
            \compactsubfiginput{figures/crossed2.tex}\par
            \textbf{(b)} type $B$
        \end{minipage}
        \par\medskip
        \begin{minipage}{0.48\textwidth}
            \centering
            \compactsubfiginput{figures/crossed3.tex}\par
            \textbf{(c)} type $D$
        \end{minipage}\hfill
        \begin{minipage}{0.48\textwidth}
            \centering
            \compactsubfiginput{figures/crossed5.tex}\par
            \textbf{(d)} Step 4, type $D$
        \end{minipage}
        \caption{The configurations in Lemma~\ref{crossed}.}
        \label{fig:crossed-configurations}
        \vspace{0.25em}
    \end{figure}
\end{proof}

By Lemma~\ref{classical}, $H^*$ is either reduced or crossed; in the reduced case, after interchanging $a_1$ and $a_2$ if necessary, $b_4 \in V(S_1(s_1,a_1))$.

The next lemmas rule out additional $H$-paths that would produce one of the forbidden types.
In Lemmas~\ref{S1pairs}, \ref{crossS1T} and~\ref{S1toT}, $H^*$ is reduced, or crossed with $b_4 \in V(S_1(s_2,a_1))$. Thus $S_1$ and $S_4$ start at $s_2$, and $S_2$ and $S_3$ at $s_1$; in the reduced case, $S_0$ is trivial and $s_0 = s_1 = s_2$. We use the following notation. Given an $H$-path $P(a,b)$ with $a \in V(S_4(s_2,a_4])$ and $b \in V(S_1(s_2,a_1))$, let $c = b_4$ if $b_4 \in V(S_1)$ and $c = b$ otherwise, and let
\[
    (A,B) =
    \begin{cases}
        (P[b_4,a_4],\ S_4) & \text{if } c = b_4,\\
        (P[b,a] \cup S_4[a,a_4],\ S_2[s_1,b_4] \cup P[b_4,a_4]) & \text{otherwise}.
    \end{cases}
\]
Thus $A$ joins $c$ to $a_4$, $B$ joins $s_2$ to $a_4$ (in the second case $H^*$ is reduced, so $s_1 = s_2$), and $A \cap B = \{a_4\}$: if $P(a,b)$ and $P(a_4,b_4)$ met in an inner vertex $p$, then $P(a_4,b_4)$ and $P[b_4,p] \cup P[p,b]$ would contradict Lemma~\ref{S1toS2}. (The same argument shows that $P(a,b)$ and $P(a_4,b_4)$ meet at most in $a_4$ in Claim~A of Lemma~\ref{notreduced12}.) Moreover, by Lemma~\ref{S4toS3T1T2T3}, every $H$-path with an end in $V(S_3 \cup T) \setminus \{s_1,a_1,a_2\}$ is internally disjoint from $A \cup B$.

\begin{lemma}\label{S1pairs}
    Suppose that $H^*$ is reduced, or crossed with $b_4 \in V(S_1(s_2,a_1))$. Let $m \in \{1,2\}$, with $m = 1$ if $H^*$ is crossed, and suppose that there is an $H$-path $P(a,b)$ with $a \in V(S_4(s_2,a_4])$ and $b \in V(S_m(s_{i_m},a_m))$.
    \begin{itemize}
        \item[(a)] There are no two disjoint $H$-paths $P(v_1,v_2)$ and $P(v_3,v_4)$ such that
    \begin{enumerate}
        \item $v_1 \in V(S_m(s_{i_m},a_m))$ and $v_2 \in V(S_3(s_1,a_3] \cup T_3(t,a_3])$;
        \item $v_3 \in V(S_m(s_{i_m},a_m))$ and $v_4 \in V(T_1 \cup T_2 \cup T_3) \setminus \{a_1,a_2,a_3\}$.
    \end{enumerate}
    Each of the two paths may meet $P(s_0,v_0)$.
        \item[(b)] Let $X = V(S_3 \cup T) \setminus \{s_1,a_1,a_2\}$, and suppose that there are two disjoint $H$-paths from $S_m(s_{i_m},a_m)$ to $X$. Then either every $H$-path from $S_m(s_{i_m},a_m)$ to $X$ ends in $T_m[t,a_m)$, or every such path ends in $S_3(s_1,a_3]$; in the latter case, the two given paths avoid $P(s_0,v_0)$, and for every $H$-path from $S_m(s_{i_m},a_m)$ to $X$, some $H$-path from $S_m(s_{i_m},a_m)$ to the same end avoids $P(s_0,v_0)$.
    \end{itemize}
\end{lemma}

\begin{proof}
    \emph{Proof of (a).}
    If $m = 2$, then $H^*$ is reduced, and we exchange the labels $a_1$ and $a_2$, together with the indices $1$ and $2$ of the legs; since all four legs of $S$ start at $s_1 = s_2$, this preserves all our assumptions except possibly that $b_4 \in V(S_1)$, which is not used below.
    Hence we may assume that $m = 1$.

    Let $c$, $A$ and $B$ be as above. Then $P(v_1,v_2)$, $P(v_3,v_4)$ and $P(s_0,v_0)$ are internally disjoint from $A \cup B$, and by Lemma~\ref{S1toS2}, $v_4 \notin V(T_2(t,a_2))$.

    The types do not involve $P(s_0,v_0)$, so we may replace it by any $H$-path from $s_0$ to $V(T) \setminus \{a_1,a_2,a_3\}$ contained in $P(s_0,v_0) \cup P(v_1,v_2) \cup P(v_3,v_4)$. Let $z$ be the first vertex of $P(s_0,v_0)$ on $P(v_1,v_2) \cup P(v_3,v_4)$, if any.
    If $z \in V(P(v_3,v_4))$, we replace $P(s_0,v_0)$ by $P[s_0,z] \cup P[z,v_4]$; if $z \in V(P(v_1,v_2))$ and $v_2 \in V(T_3(t,a_3))$, by $P[s_0,z] \cup P[z,v_2]$.
    If $z \in V(P(v_1,v_2))$ and $v_2 \in V(S_3(s_1,a_3])$, let $Q_0 = P[s_0,z')$ if $P(s_0,v_0)$ meets $P(v_3,v_4)$, where $z'$ is its first vertex on $P(v_3,v_4)$, and $Q_0 = P(s_0,v_0)$ otherwise. In the first case we replace $P(s_0,v_0)$ by $P[s_0,z'] \cup P[z',v_4]$; in both cases we replace $P(v_1,v_2)$ by $P[v_1,y] \cup P[y,y'] \cup P[y',v_2]$, where $y$ and $y'$ are the first and the last vertex of $P(v_1,v_2)$ on $Q_0$.
    The new $P(v_1,v_2)$ is a path: $P[v_1,y]$ and $P[y',v_2]$ are disjoint subpaths of the old one meeting $Q_0$ only in $y$ and $y'$, and $P[y,y'] \subseteq Q_0$. It meets the new $P(s_0,v_0)$ exactly in $P[y,y']$, since the rest of the new $P(s_0,v_0)$, if any, lies in $P(v_3,v_4)$. Hence $P(v_1,v_2)$ and $P(v_3,v_4)$ are still disjoint, and each of them meets $P(s_0,v_0)$ in at most one subpath; for $P(v_3,v_4)$, and for $P(v_1,v_2)$ when $v_2 \in V(T_3(t,a_3))$, this subpath ends at $v_0$.

    Suppose first that $P(v_1,v_2)$ is disjoint from $P(s_0,v_0)$, and that $v_4, v_0 \notin V(T_3[v_2,a_3])$ if $v_2 \in V(T_3)$. Let $x$ be the vertex of $\{v_1,c\}$ closer to $s_2$ on $S_1$, let
    \[
        L =
        \begin{cases}
            S_3[v_2,a_3] & \text{if } v_2 \in V(S_3),\\
            \emptyset & \text{otherwise},
        \end{cases}
    \]
    and let
    \[
        S_K = T \cup P[s_0,v_0] \cup S_0[s_0,s_2] \cup B, \qquad T_K = S_1[x,a_1] \cup P[v_1,v_2] \cup L \cup A.
    \]
    Let $U_1$ be the subpath of $S_0[s_0,s_1] \cup S_3$ from $s_0$ to its first vertex on $T_K \cup \{a_3\}$, let $Q$ be the subpath of $P[v_3,v_4]$ from $v_3$ to its first vertex on $S_K$, and let
    \[
        U_2 =
        \begin{cases}
            Q & \text{if } v_3 \in V(S_1[x,a_1]),\\
            S_1[x,v_3] \cup Q & \text{otherwise}.
        \end{cases}
    \]
    Then, after interchanging the labels $a_1$ and $a_2$, $K = S_K \cup T_K \cup U_1 \cup U_2$ is of type $B$, as in Figure~\ref{fig:s1pairs}(a). Here $a'_3 = a_3$ if $v_2 \in V(S_3)$ and $a'_3 = v_2$ otherwise; in the latter case $U_2$ avoids $T_3[v_2,a_3]$, since $v_4 \notin V(T_3[v_2,a_3])$ and $P(v_3,v_4)$ meets $P(s_0,v_0)$ at most in a subpath ending at $v_0$.
    Symmetrically, if $v_2 \in V(T_3)$, $v_4 \in V(T_3(v_2,a_3))$, $P(v_3,v_4)$ is disjoint from $P(s_0,v_0)$ and $v_0 \notin V(T_3[v_4,a_3])$, the same construction with $P(v_1,v_2)$ and $P(v_3,v_4)$ interchanged gives a subgraph of type $B$.

    Suppose now that $v_2 \in V(S_3(s_1,a_3])$. By the above, $P(v_1,v_2)$ meets $P(s_0,v_0)$, in $P[y,y']$. Let $u$ be the vertex of $\{y,y'\}$ farther from $s_0$ on $P(s_0,v_0)$, let $x$ be the vertex of $\{v_3,c\}$ closer to $s_2$ on $S_1$, and let $R$ be the subpath of $P[u,v_0] \cup T[v_0,a_2]$ from $u$ to its first vertex on $P[v_3,v_4] \cup T[v_4,a_2]$.
    \[
        \begin{aligned}
            S_K ={}& S_1[x,a_1] \cup P[v_3,v_4] \cup T[v_4,a_2] \cup A, \\
            T_K ={}& S_0 \cup S_2 \cup S_3 \cup B, \\
            U_K ={}& R \cup P[u,v_2], \\
            V_K ={}& S_1[s_2,x] \cup P[v_1,u] \cup P[s_0,y'].
        \end{aligned}
    \]
    Then $K = S_K \cup T_K \cup U_K \cup V_K$ contains a subgraph of type $D$, as in Figure~\ref{fig:s1pairs}(b). Indeed, use $S_1[x,s_2]$ for the $S_K$--$T_K$ connector. For the $S_K$--$U_K$ connector, use $P[v_1,u]$, preceded by $S_1[x,v_1]$ if $v_1$ precedes $x$; for the $T_K$--$U_K$ connector, use $P[s_0,y']$, cutting at the first vertex of $U_K$. Trim the first of these two routes at its last vertex of $S_K$. If $u=y'$, the two connectors to $U_K$ share only a final subpath ending at $u$; if $u=y\ne y'$, they are disjoint. Any other common nontrivial segment lies on $S_1$ and starts at the common end in $S_K$. No vertex lies on all three connectors.

    Finally, suppose that $v_2 \in V(T_3(t,a_3))$. By the above, letting $P(y_1,y_2) = P(v_1,v_2)$ and $P(x_1,x_2) = P(v_3,v_4)$ if $v_4 \in V(T_3(v_2,a_3))$, and $P(y_1,y_2) = P(v_3,v_4)$ and $P(x_1,x_2) = P(v_1,v_2)$ otherwise, $P(y_1,y_2)$ is disjoint from $P(s_0,v_0)$, $v_0 \in V(T_3[x_2,a_3))$ and $y_2 \notin V(T_3[x_2,a_3])$.
    Let $u_2 = y_2$ if $y_2 \in V(T_3)$ and $u_2 = t$ otherwise, and let $x$ be the vertex of $\{y_1,c\}$ closer to $s_2$ on $S_1$.
    \[
        \begin{aligned}
            S_K ={}& S_1[x,a_1] \cup P[y_1,y_2] \cup T[y_2,a_2] \cup A, \\
            T_K ={}& S_0 \cup S_2 \cup S_3 \cup B, \\
            U_K ={}& T_3[u_2,a_3], \\
            V_K ={}& S_1[s_2,x] \cup P[x_1,x_2] \cup P[s_0,v_0].
        \end{aligned}
    \]
    Then $K = S_K \cup T_K \cup U_K \cup V_K$ is of type $D$, with $u = v_0$, as in Figure~\ref{fig:s1pairs}(c); if $P(x_1,x_2)$ meets $P(s_0,v_0)$, then $x_2 = v_0 = u$ by the above, and the two paths share a final subpath ending at $u$, as allowed in type $D$.

    In all cases we obtain one of the forbidden types, a contradiction.

    \begin{figure}[htbp]
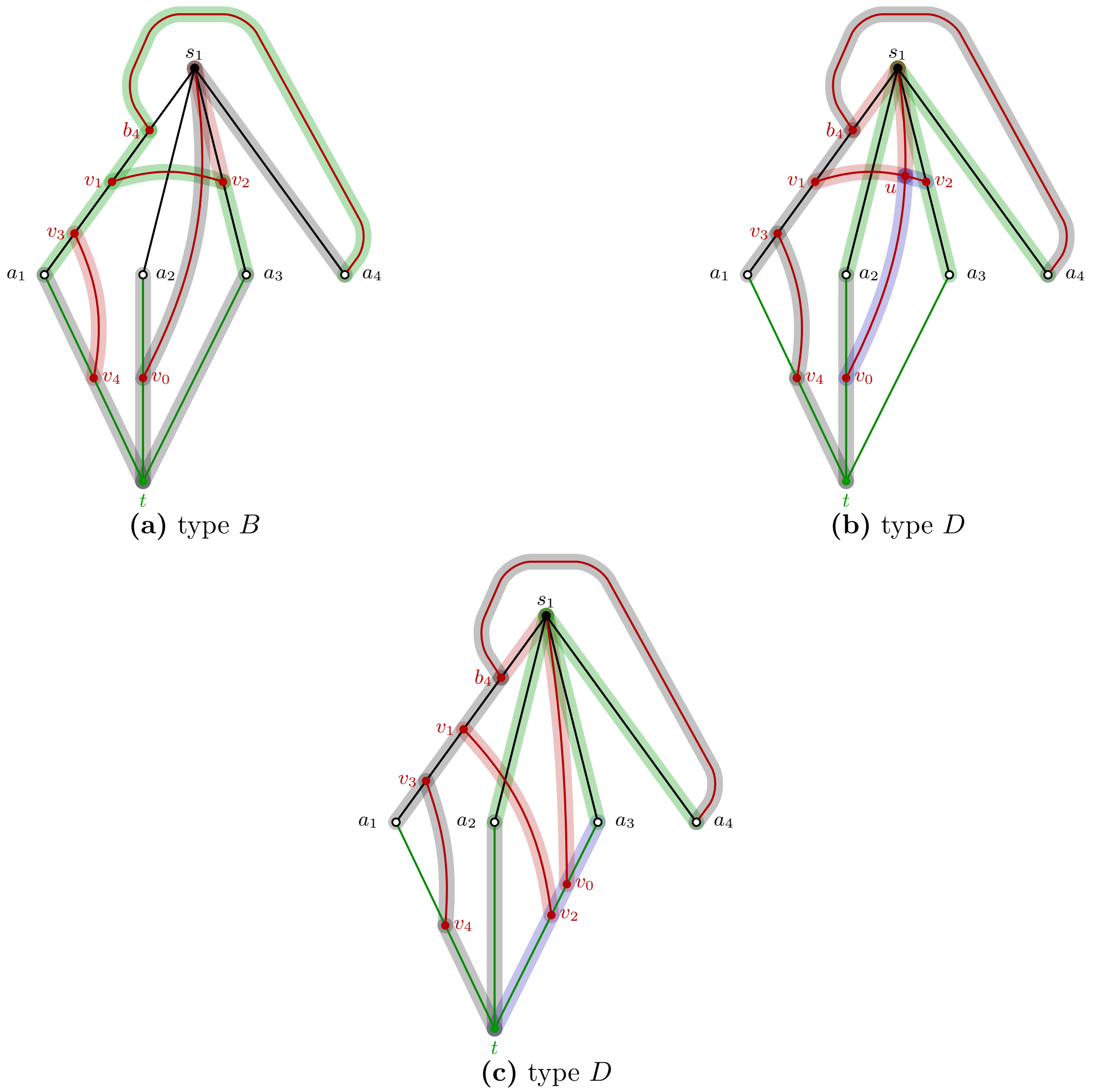

        \centering
        \begin{minipage}{0.48\textwidth}
            \centering
            \compactsubfiginput{figures/S1toS3.tex}\par
            \textbf{(a)} type $B$
        \end{minipage}\hfill
        \begin{minipage}{0.48\textwidth}
            \centering
            \compactsubfiginput{figures/S1toS3n2.tex}\par
            \textbf{(b)} type $D$
        \end{minipage}

        \vspace{0.4em}

        \begin{minipage}{0.48\textwidth}
            \centering
            \compactsubfiginput{figures/S1toT3n2.tex}\par
            \textbf{(c)} type $D$
        \end{minipage}
        \caption{The configurations in Lemma~\ref{S1pairs}.}
        \label{fig:s1pairs}
    \end{figure}

    \emph{Proof of (b).}
    Let $P(v_1,v_2)$ and $P(v_3,v_4)$ be the two given paths. By Lemma~\ref{S1toS2}, $v_2,v_4 \notin V(T_n(t,a_n])$, where $n \neq m$, and by~(a), $v_2, v_4 \in V(S_3(s_1,a_3])$ or $v_2, v_4 \in V(T_m[t,a_m))$.
    In the first case, both paths avoid $P(s_0,v_0)$: otherwise, let $z$ be the vertex of $P(s_0,v_0)$ on $P(v_1,v_2) \cup P(v_3,v_4)$ closest to $v_0$, say on $P(v_3,v_4)$; then $P[v_3,z] \cup P[z,v_0]$ and $P(v_1,v_2)$ contradict~(a).
    Now let $P(v_5,v_6)$ be an $H$-path from $S_m(s_{i_m},a_m)$ to $X$. If it meets $P(v_1,v_2) \cup P(v_3,v_4)$, let $y$ be its vertex on $P(v_1,v_2) \cup P(v_3,v_4)$ closest to $v_6$, say on $P(v_1,v_2)$, and let $R = P[v_1,y] \cup P[y,v_6]$; otherwise let $R = P(v_5,v_6)$. Then $R$ is an $H$-path from $S_m(s_{i_m},a_m)$ to $v_6$ disjoint from $P(v_3,v_4)$ or from $P(v_1,v_2)$, so $R$ and this path are two disjoint $H$-paths from $S_m(s_{i_m},a_m)$ to $X$; applying to them~(a) and the previous paragraph, which holds for any two such paths, shows that $v_6$ lies on the same side as $v_2$ and $v_4$, and that $R$ avoids $P(s_0,v_0)$ in the first case.
\end{proof}

\begin{lemma}\label{crossS1T}
    Suppose that $H^*$ is reduced, or crossed with $b_4 \in V(S_1(s_2,a_1))$, and that there is an $H$-path $P(a,b)$ with $a \in V(S_4(s_2,a_4])$ and $b \in V(S_1(s_2,a_1))$, and let $c$ be as above. There are no two disjoint $H$-paths $P(v_1,v_2)$ and $P(v_3,v_4)$ such that
    \begin{enumerate}
        \item $v_1 \in V(S_1(s_2,a_1))$ and $v_2 \in V(T_1[t,a_1))$;
        \item $v_3 \in V(T_1(v_2,a_1))$ and $v_4 \in V(S_0 \cup S_1[s_2,x) \cup S_2[s_1,a_2) \cup S_3[s_1,a_3))$, where $x$ is the vertex of $\{v_1,c\}$ closer to $s_2$ on $S_1$.
    \end{enumerate}
\end{lemma}

\begin{proof}
    Suppose, for a contradiction, that such paths exist. Both are internally disjoint from $A \cup B$; moreover $v_2, v_3 \notin V(A \cup B)$, $v_4 \notin V(A)$, $v_1 \notin V(B)$, and $v_1 \in V(A)$ only if $v_1 = c$. Let
    \[
        L =
        \begin{cases}
            S_1[s_2,v_4] \cup P[v_4,v_3] \cup T_1[v_3,a_1] & \text{if } v_4 \in V(S_1),\\
            P[v_4,v_3] \cup T_1[v_3,a_1] & \text{otherwise},
        \end{cases}
    \]
    and
    \[
        \begin{aligned}
            S_K ={}& S_0 \cup S_2 \cup S_3 \cup B \cup L, \\
            T_K ={}& T_1[t,v_2] \cup T_2 \cup T_3 \cup P[v_2,v_1] \cup S_1[x,a_1] \cup A.
        \end{aligned}
    \]
    Then $K = S_K \cup T_K$ is of type $A$, as in Figure~\ref{fig:cross-s1-t}(a), a contradiction; note that $T_1[t,v_2]$ and $T_1[v_3,a_1]$ are disjoint since $v_3 \in V(T_1(v_2,a_1))$, and that $S_1[s_2,v_4]$ and $S_1[x,a_1]$ are disjoint since $v_4 \in V(S_1[s_2,x))$.

    \begin{figure}[htbp]
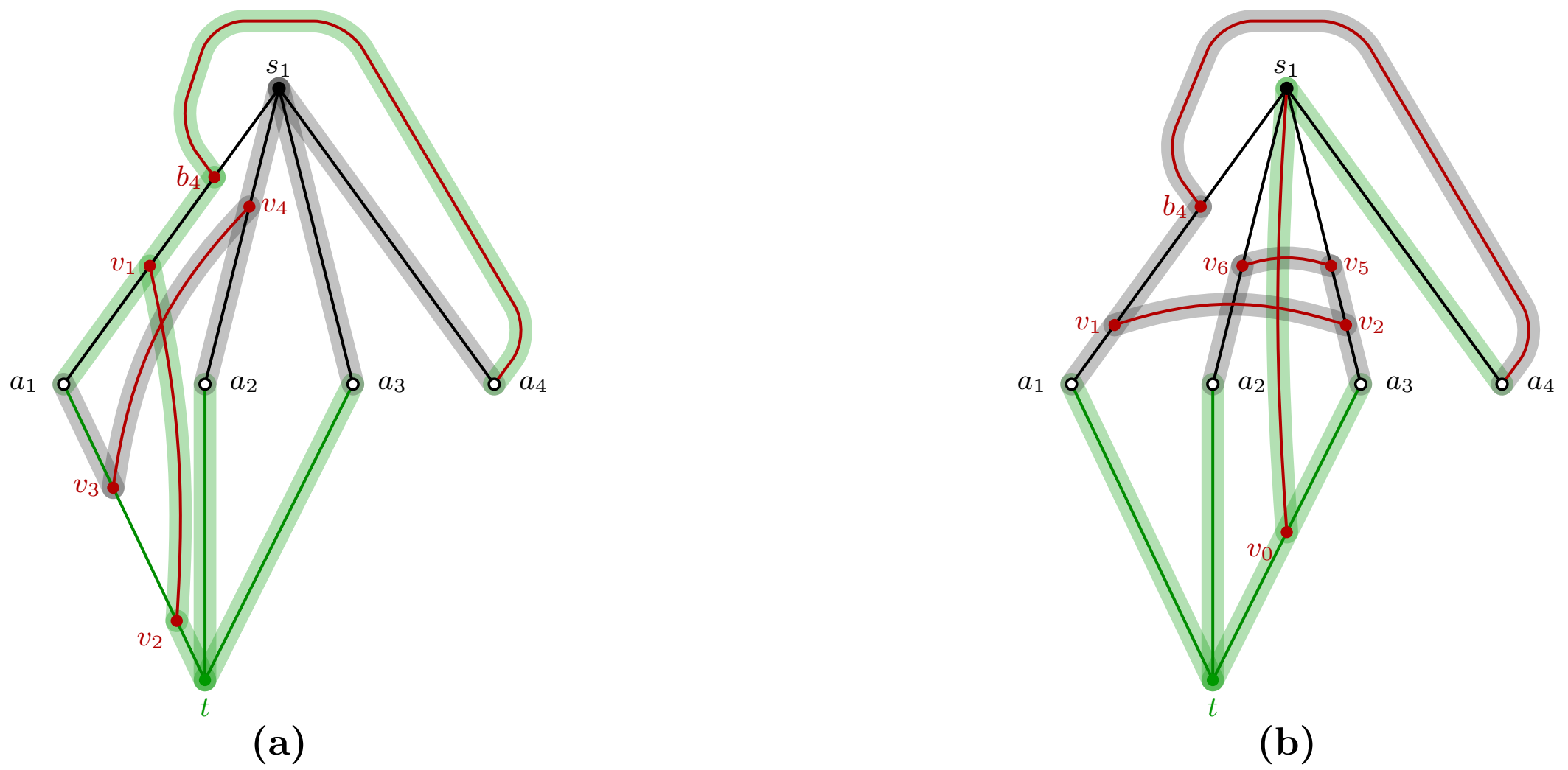

        \centering
        \begin{minipage}{0.48\textwidth}
            \centering
            \compactsubfiginput{figures/crossS1T.tex}\par
            \textbf{(a)}
        \end{minipage}\hfill
        \begin{minipage}{0.48\textwidth}
            \centering
            \compactsubfiginput{figures/allSi.tex}\par
            \textbf{(b)}
        \end{minipage}
        \caption{(a) The type-$A$ configuration in Lemma~\ref{crossS1T}, drawn for $v_4 \in V(S_2(s_1,a_2))$ and $x = c = b_4$. (b) The type-$A$ configuration in the proof of Theorem~\ref{structure}.}
        \label{fig:cross-s1-t}
    \end{figure}
\end{proof}

   \begin{lemma}\label{notreduced12}
    Suppose that $H^*$ is reduced and that $b_4 \in V(S_1(s_1,a_1))$.
    There is no $H$-path $P(a,b)$ with $a \in V(S_4(s_1,a_4])$ and $b \in V(S_2(s_1,a_2))$.
\end{lemma}

\begin{proof}
    Suppose, for a contradiction, that such a path exists, and let $T'_4 = P[b,a] \cup S_4[a,a_4]$. Let $\mathcal{H}$ be the set of augmented frames obtained from $H^*$ by replacing $T$, $P(s_0,v_0)$ and $S_3$ so that $P(a_4,b_4)$ and $P(a,b)$ remain $H$-paths; they have the same $S_0$ as $H^*$, hence are minimal. Claims~A and~B below hold for every member of $\mathcal{H}$.

    \smallskip
    \emph{Claim A. There are no two $H$-paths $P(v_1,v_2)$ and $P(v_3,v_4)$ with $v_1 \in V(S_1(s_1,a_1))$, $v_3 \in V(S_2(s_1,a_2))$ and $v_2,v_4 \in V(S_3 \cup T) \setminus \{s_1,a_1,a_2\}$.}
    Suppose otherwise. By Lemma~\ref{S4toS3T1T2T3}, both paths are internally disjoint from $P(a_4,b_4)$ and $P(a,b)$, and by Lemma~\ref{S1toS2} they meet at most in $v_2 = v_4$; for the same reason they do not both meet $P(s_0,v_0)$. If one of them ends in $S_3$ and meets $P(s_0,v_0)$, we cut it at $P(s_0,v_0)$ and follow $P(s_0,v_0)$ to $v_0$.

    Since $a_3 \in V(S_3)$, we distinguish the following cases.
    If $v_2, v_4 \in V(T) \setminus \{a_3\}$, let
    \[
        \begin{aligned}
            S_K ={}& T \cup P[v_2,v_1] \cup S_1[v_1,b_4] \cup P[b_4,a_4], \\
            T_K ={}& S_2 \cup S_3 \cup T'_4, \\
            U_K ={}& R \cup P[v_3,v_4],
        \end{aligned}
    \]
    where $R$ is the subpath of $S_1$ from $s_1$ to its first vertex on $S_K$; then $K$ is of type $B$, as in Figure~\ref{fig:nr12}(a), the hub of $T_K$ being $b$.
    If $v_2 \in V(T) \setminus \{a_3\}$ and $v_4 \in V(S_3)$, let $P'$ be the subpath of $P(s_0,v_0) \cup T[v_0,a_2]$ from $s_1$ to $P[v_1,v_2] \cup T[v_2,a_2]$, and
    \[
        \begin{aligned}
            S_K ={}& T[v_2,a_2] \cup P[v_2,v_1] \cup S_1[v_1,a_1] \cup S_1[b_4,a_1] \cup P[b_4,a_4], \\
            T_K ={}& S_2[v_3,a_2] \cup S_2[b,a_2] \cup P[v_3,v_4] \cup S_3[v_4,a_3] \cup T'_4, \\
            U_K ={}& P' \cup S_3[s_1,v_4] \cup S_4[s_1,a],
        \end{aligned}
    \]
    and let $V_K$ consist of the subpaths of $S_1$ and $S_2$ from $s_1$ to their first vertices on $S_K$ and on $T_K$, respectively; then $K$ is of type $E$, as in Figure~\ref{fig:nr12}(b). Here $S_1[v_1,a_1] \cup S_1[b_4,a_1]$ is the longer of these two final segments of $S_1$, and similarly for $S_2[v_3,a_2] \cup S_2[b,a_2]$ and in the constructions below; writing the union avoids specifying which of $v_1,b_4$, respectively $v_3,b$, is closer to $s_1$. If $v_2 \in V(S_3)$ and $v_4 \in V(T) \setminus \{a_3\}$, let $P''$ be the subpath of $P(s_0,v_0) \cup T[v_0,a_1]$ from $s_1$ to $P[v_3,v_4] \cup T[v_4,a_1]$, and
    \[
        \begin{aligned}
            S_K ={}& T[v_4,a_1] \cup P[v_4,v_3] \cup S_2[v_3,a_2] \cup S_2[b,a_2] \cup T'_4, \\
            T_K ={}& S_1[v_1,a_1] \cup S_1[b_4,a_1] \cup P[v_1,v_2] \cup S_3[v_2,a_3] \cup P[b_4,a_4], \\
            U_K ={}& P'' \cup S_3[s_1,v_2] \cup S_4[s_1,a],
        \end{aligned}
    \]
    and let $V_K$ consist of the subpaths of $S_2$ and $S_1$ from $s_1$ to their first vertices on $S_K$ and on $T_K$, respectively; then $K$ is of type $E$ after interchanging $a_1$ and $a_2$, with $u_4 = a$ now on the leg of $S_K$ to $a_4$.
    If $v_2, v_4 \in V(S_3)$, let
    \[
        \begin{aligned}
            S_K ={}& T \cup P[s_0,v_0] \cup S_4, \\
            T_K ={}& S_1[v_1,a_1] \cup S_1[b_4,a_1] \cup P[v_1,v_2] \cup S_3[v_2,a_3] \cup P[b_4,a_4],
        \end{aligned}
    \]
    let $U_1$ be the path $P[a,b] \cup S_2[b,v_3] \cup P[v_3,v_4]$ followed by the subpath of $S_3[v_4,a_3]$ from $v_4$ to its first vertex on $T_K$, and let $U_2$ be the subpath of $S_1$ from $s_1$ to its first vertex on $T_K$; then $K$ is of type $B$ after interchanging $a_1$ and $a_2$, as in Figure~\ref{fig:nr12}(c), since $P(a,b)$ is disjoint from $P(s_0,v_0)$ by Lemma~\ref{S4toS3T1T2T3}.

    \smallskip
    \emph{Claim B. There is an $H$-path $P(v_1,v_2)$ with $v_1 \in V(S_1(s_1,a_1) \cup S_2(s_1,a_2))$ and $v_2 \in V(S_3 \cup T) \setminus \{s_1,a_1,a_2\}$.}
    As $\{s_1,a_1,a_2,a_4\}$ is not a $4$-separator, there is an $H$-path from $S_1(s_1,a_1) \cup S_2(s_1,a_2) \cup S_4(s_1,a_4)$ to $V(S_3 \cup T) \setminus \{s_1,a_1,a_2\}$, and by Lemma~\ref{S4toS3T1T2T3} it does not start in $S_4$.

    \smallskip
    \emph{Claim C. There is no $H$-path $P(v_1,v_2)$ with $v_1 \in V(S_1(s_1,a_1))$ and $v_2 \in V(S_3 \cup T) \setminus \{s_1,a_1,a_2\}$.}
    Suppose otherwise.
    In every member of $\mathcal{H}$ having such a path, every $H$-path from $S_2(s_1,a_2)$ or $S_4(s_1,a_4]$ ends in $S_1 \cup S_2 \cup S_4$, by Claim~A and Lemmas~\ref{S1toS2} and~\ref{S4toS3T1T2T3}; hence every path from $a_4$ to $V(S_3 \cup T) \setminus \{s_1,a_1,a_2\}$ in $G - \{s_1,a_1,a_2\}$ meets $S_1(s_1,a_1)$. As $\{s_1,a_1,a_2,w\}$ is not a $4$-separator for any vertex $w$, Menger's theorem gives two disjoint such $H$-paths. We apply Lemma~\ref{S1pairs} with $m = 1$, the path $P(a_4,b_4)$ playing the role of $P(a,b)$.

    First suppose that some member of $\mathcal{H}$ has such a path $P(v_1,v_2)$ with $v_2 \in V(S_3(s_1,a_3])$ disjoint from $P(s_0,v_0)$. Choose the member and the path with $|E(S_3[v_2,a_3])|$ minimal.
    By Lemma~\ref{S1pairs}(b), applied to two disjoint paths as above, and since $P(v_1,v_2)$ ends in $S_3$, these two paths end in $S_3(s_1,a_3]$ and avoid $P(s_0,v_0)$, hence end in $S_3(s_1,v_2]$ by the choice of $P(v_1,v_2)$; so one of them, $P(v_3,v_4)$, has $v_4 \in V(S_3(s_1,v_2))$.
    As $\{s_1,a_1,a_2,v_2\}$ is not a $4$-separator, there is an $H$-path $P(v_5,v_6)$ with $v_5 \in V(S_1 \cup S_2 \cup S_3(s_1,v_2) \cup S_4) \setminus \{s_1,a_1,a_2\}$ and $v_6 \in V(S_3(v_2,a_3] \cup T) \setminus \{a_1,a_2,v_2\}$.
    By Claim~A and Lemma~\ref{S4toS3T1T2T3}, $v_5 \notin V(S_2 \cup S_4)$, and by Lemma~\ref{S1pairs}(b) and the choice of $P(v_1,v_2)$, $v_5 \notin V(S_1)$. Hence $v_5 \in V(S_3(s_1,v_2))$, and $P(v_5,v_6)$ is disjoint from $P(v_1,v_2)$, $P(v_3,v_4)$, $P(a_4,b_4)$ and $P(a,b)$: otherwise, cutting it at the first common vertex from $v_6$ gives an $H$-path from $S_1$ or $S_4$ to $v_6$.
    If $P(v_5,v_6)$ meets $P(s_0,v_0)$, let $p$ be its first vertex on $P(s_0,v_0)$ from $v_5$; then $P[v_5,p] \cup P[p,v_0]$ is an $H$-path from $S_3(s_1,v_2)$ to $T \setminus \{a_1,a_2,a_3\}$, disjoint from $P(v_1,v_2)$, $P(v_3,v_4)$, $P(a_4,b_4)$ and $P(a,b)$, and we conclude as in the last case below. So we may assume that $P(v_5,v_6)$ avoids $P(s_0,v_0)$.
    If $v_6 \in V(S_3)$, then replacing $S_3$ by $S_3[s_1,v_5] \cup P[v_5,v_6] \cup S_3[v_6,a_3]$ and $P(v_1,v_2)$ by $P[v_1,v_2] \cup S_3[v_2,v_6]$ contradicts the choice.
    Otherwise $P(a_4,b_4)$, $P(v_1,v_2)$, $P(v_3,v_4)$ and $P(v_5,v_6)$ contradict Lemma~\ref{S3toT}, applied after interchanging $a_1$ and $a_2$.

    Hence, by Lemma~\ref{S1pairs}(b), in every member of $\mathcal{H}$, every $H$-path $P(v_1,v_2)$ with $v_1 \in V(S_1(s_1,a_1))$ and $v_2 \in V(S_3 \cup T) \setminus \{s_1,a_1,a_2\}$ has $v_2 \in V(T_1[t,a_1))$.
    Choose a member and such a path $P(q,q')$ with $|E(T_1[t,q'])|$ minimal.
    As $\{s_1,a_1,a_2,q'\}$ is not a $4$-separator, there is an $H$-path $P(v_1,v_2)$ with $v_1 \in V(S_1 \cup S_2 \cup S_4 \cup T_1(q',a_1)) \setminus \{s_1,a_1,a_2\}$ and $v_2 \in V(S_3 \cup T_1[t,q') \cup T_2 \cup T_3) \setminus \{s_1,a_1,a_2\}$.
    By the above, the choice of $q'$, Claim~A and Lemma~\ref{S4toS3T1T2T3}, $v_1 \in V(T_1(q',a_1))$. Moreover, $P(v_1,v_2)$ is disjoint from $P(q,q')$, since otherwise $P[q,y] \cup P[y,v_2]$, where $y$ is the last common vertex, contradicts the above; and $P(v_1,v_2)$ and $P(q,q')$ are disjoint from $P(a_4,b_4)$ and $P(a,b)$ by Lemma~\ref{S4toS3T1T2T3}.

    If $v_2 \in V(T_j(t,a_j])$ for some $j \in \{2,3\}$, let $k$ be the other element of $\{2,3\}$, and replace $T$ by the tripod
    \[
        T' = T_1 \cup P[v_1,v_2] \cup T_j[v_2,a_j] \cup T_k,
    \]
    with center $v_1$ and legs $T_1[v_1,a_1]$, $P[v_1,v_2] \cup T_j[v_2,a_j]$ and $T_1[v_1,t] \cup T_k$. Let
    \[
        P' =
        \begin{cases}
            P[s_0,v_0] \cup T_j[v_0,t] & \text{if } v_0 \in V(T_j(t,v_2)),\\
            P[s_0,v_0] & \text{otherwise},
        \end{cases}
    \]
    and replace $P(s_0,v_0)$ by the subpath of $P'$ from $s_0$ to its first vertex on $T'$, which is not a root. This gives a member of $\mathcal{H}$ in which $P(q,q')$ is an $H$-path from $S_1(s_1,a_1)$ to $q'$, an interior vertex of the leg of $T'$ to $a_k$; this contradicts the above.
    If $v_2 \in V(S_3(s_1,a_3))$, then $P(q,q')$ and $P(v_1,v_2)$ contradict Lemma~\ref{crossS1T}, applied to $P(a_4,b_4)$.
    If $v_2 \in V(T_1[t,q'))$, delete $T_1(v_2,v_1)$, add $P(v_1,v_2)$, and replace $P(s_0,v_0)$, or $P(s_0,v_0) \cup T_1[v_0,v_2]$ if $v_0 \in V(T_1(v_2,v_1))$, by its subpath from $s_0$ to its first vertex on the new tripod; this vertex is not a root. In the resulting member of $\mathcal{H}$, $P[q,q'] \cup T_1[q',v_2]$ is an $H$-path from $S_1(s_1,a_1)$ to $v_2$, contrary to the choice of $q'$.

    \smallskip
    \emph{Claim D. There is no $H$-path $P(v_1,v_2)$ with $v_1 \in V(S_2(s_1,a_2))$ and $v_2 \in V(S_3 \cup T) \setminus \{s_1,a_1,a_2\}$.}
    The proof of Claim~C applies after exchanging the roles of $a_1, S_1, T_1, b_4, P(a_4,b_4)$ and of $a_2, S_2, T_2, b, P(a,b)$: Claim~C replaces Claim~A, Lemma~\ref{S1pairs} is applied with $m = 2$ to $P(a,b)$, and Lemma~\ref{S3toT} is applied as stated. With $P(q,q')$, $q'$ and $P(v_1,v_2)$ as in Claim~C, now $q \in V(S_2(s_1,a_2))$, $q' \in V(T_2[t,a_2))$ and $v_1 \in V(T_2(q',a_2))$.

    If $v_2 \in V(T_j(t,a_j])$ for some $j \in \{1,3\}$, we conclude as in Claim~C, with the tripod $T_2 \cup P[v_1,v_2] \cup T_j[v_2,a_j] \cup T_k$ centered at $v_1$, where $k$ is the other element of $\{1,3\}$.
    If $v_2 \in V(S_3(s_1,a_3))$, then, after interchanging $a_1$ and $a_2$, $P(q,q')$ and $P(v_1,v_2)$ contradict Lemma~\ref{crossS1T}, applied to $P(a,b)$.
    If $v_2 \in V(T_2[t,q'))$, we conclude as in Claim~C.

    \smallskip
    Claims~B, C and~D give a contradiction.

    \begin{figure}[htbp]
        \centering
        \begin{minipage}{0.48\textwidth}
            \centering
            \compactsubfiginput{figures/reduced_type_b.tex}\par
            \textbf{(a)} Claim~A, type $B$
        \end{minipage}\hfill
        \begin{minipage}{0.48\textwidth}
            \centering
            \compactsubfiginput{figures/caseE.tex}\par
            \textbf{(b)} Claim~A, type $E$
        \end{minipage}

        \vspace{0.4em}

        \begin{minipage}{0.48\textwidth}
            \centering
            \compactsubfiginput{figures/bothS3.tex}\par
            \textbf{(c)} Claim~A, type $B$
        \end{minipage}
        \caption{The configurations in Lemma~\ref{notreduced12}.}
        \label{fig:nr12}
    \end{figure}
\end{proof}

\begin{lemma}\label{S1toT}
    Suppose that $H^*$ is reduced or crossed, and that $b_4 \in V(S_1(s_2,a_1))$. There is no $H$-path $P(v_1,v_2)$ with $v_1 \in V(S_1(s_2,a_1))$ and $v_2 \in V(T) \setminus \{a_1,a_2,a_3\}$.
\end{lemma}

\begin{proof}
    Suppose, for a contradiction, that such a path exists. Let $\mathcal{H}$ be the set of augmented frames obtained from $H^*$ by replacing $T$ and $P(s_0,v_0)$ so that $P(a_4,b_4)$ remains an $H$-path; they have the same $S_0$ as $H^*$, hence are minimal.
    By Lemma~\ref{S1toS2} and, if $H^*$ is crossed, Lemma~\ref{S0min}, every $H$-path from $S_1(s_2,a_1)$ to $V(S_0 \cup S_2 \cup S_3 \cup T) \setminus \{s_2,a_1,a_2\}$ ends in $V(S_3 \cup T) \setminus \{s_1,a_1,a_2\}$. By Lemma~\ref{notreduced12} if $H^*$ is reduced, and by Lemmas~\ref{crossed} and~\ref{S0min} if $H^*$ is crossed, together with Lemma~\ref{S4toS3T1T2T3}, every $H$-path from $S_4(s_2,a_4]$ ends in $S_1 \cup S_4 \cup \{a_2\}$.

    Consider a member of $\mathcal{H}$ having such a path $P(v_1,v_2)$. Every path from $a_4$ to $V(S_0 \cup S_2 \cup S_3 \cup T) \setminus \{s_2,a_1,a_2\}$ in $G - \{s_2,a_1,a_2\}$ meets $S_1(s_2,a_1)$; as $\{s_2,a_1,a_2,w\}$ is not a $4$-separator for any vertex $w$, Menger's theorem gives two disjoint $H$-paths from $S_1(s_2,a_1)$ to $V(S_3 \cup T) \setminus \{s_1,a_1,a_2\}$.
    By Lemma~\ref{S1pairs}(b), applied with $m = 1$ to $P(a_4,b_4)$ and to these two paths, and since $P(v_1,v_2)$ ends in $T \setminus \{a_3\}$, every $H$-path from $S_1(s_2,a_1)$ to $V(S_3 \cup T) \setminus \{s_1,a_1,a_2\}$ ends in $T_1[t,a_1)$.

    Choose a member and such a path $P(q,q')$ with $|E(T_1[t,q'])|$ minimal.
    Since $\{s_2,a_1,a_2,q'\}$ is not a $4$-separator, there is an $H$-path $P(v_1,v_2)$ with $v_1 \in V(S_1 \cup S_4 \cup T_1(q',a_1)) \setminus \{s_2,a_1\}$ and $v_2 \in V(S_0 \cup S_2 \cup S_3 \cup T_1[t,q') \cup T_2 \cup T_3) \setminus \{s_2,a_1,a_2\}$.
    By the above and the choice of $q'$, $v_1 \in V(T_1(q',a_1))$, and we conclude as in Claim~C of Lemma~\ref{notreduced12}: if $v_2 \in V(T_2 \cup T_3) \setminus \{t\}$, by the same re-routing of $T$; if $v_2 \in V(S_0 \cup S_2 \cup S_3) \setminus V(T)$, since $P(q,q')$ and $P(v_1,v_2)$ contradict Lemma~\ref{crossS1T}, applied to $P(a_4,b_4)$; and if $v_2 \in V(T_1[t,q'))$, by the choice of $q'$.
\end{proof}

\begin{proof}[Proof of Theorem~\ref{structure}]

    First, $H^*$ is not crossed. Indeed, suppose otherwise. By Lemma~\ref{crossed}, there is no $H$-path between $S_4(s_2,a_4]$ and $S_2(s_1,a_2)$; in particular $b_4 \in V(S_1(s_2,a_1))$; by Lemma~\ref{S4toS3T1T2T3}, there is no $H$-path between $S_4$ and $V(S_3 \cup T_1 \cup T_2 \cup T_3) \setminus \{s_1,a_1,a_2\}$; and by Lemma~\ref{S0min}, there is no $H$-path from $S_1(s_2,a_1) \cup S_4(s_2,a_4]$ to $V(S_0) \setminus \{s_2\}$. Similarly, by Lemma~\ref{crossedS1toS2S3}, there is no $H$-path between $S_1(s_2,a_1)$ and $V(S_2 \cup S_3) \setminus \{s_1,a_2,a_3\}$; by Lemma~\ref{S1toT}, there is no $H$-path between $S_1(s_2,a_1)$ and $V(T_1 \cup T_2 \cup T_3) \setminus \{a_1,a_2,a_3\}$. Therefore $\{s_2,a_1,a_2,a_3\}$ would be a $4$-separator.

    Hence, by Lemmas~\ref{classical} and~\ref{notreduced12}, the augmented frame $H^*$ is reduced and, up to exchanging $a_1$ and $a_2$, $b_4 \in V(S_1(s_2,a_1))$; in particular, $s_1 = s_2$. Hence, by the definition of an augmented frame in the reduced case, $s_0 = s_1 = s_2$.

    Let $\mathcal{H}$ be the set of augmented frames obtained from $H^*$ by replacing $S_3$ so that $P(s_0,v_0)$ and $P(a_4,b_4)$ remain $H$-paths; they have the same $S_0$ as $H^*$, hence are minimal.
    By Lemmas~\ref{notreduced12} and~\ref{S4toS3T1T2T3}, every $H$-path from $S_4(s_1,a_4]$ ends in $S_1 \cup S_4 \cup \{a_2\}$. By Lemmas~\ref{S1toS2} and~\ref{S1toT}, every $H$-path from $S_1(s_1,a_1)$ to $V(S_2 \cup S_3 \cup T) \setminus \{s_1,a_1,a_2\}$ ends in $S_3(s_1,a_3]$ and is disjoint from $P(s_0,v_0)$ (otherwise, cutting it at $P(s_0,v_0)$ and following $P(s_0,v_0)$ to $v_0$ contradicts Lemma~\ref{S1toT}).
    As $\{s_1,a_1,a_2,a_4\}$ is not a $4$-separator, there is such a path in every member of $\mathcal{H}$.
    Choose a member and such a path $P(v_1,v_2)$ with $|E(S_3[v_2,a_3])|$ minimal. The paths $P(v_1,v_2)$ and $P(a_4,b_4)$ are internally disjoint, and their only possible common vertex is $v_1 = b_4$. Indeed, any other common vertex would yield an $H$-path from $a_4$ to $S_3(s_1,a_3]$, contrary to Lemma~\ref{S4toS3T1T2T3}: follow $P(a_4,b_4)$ from $a_4$ to its first vertex on $P(v_1,v_2)$ other than $v_1$, and then $P(v_1,v_2)$ to $v_2$.

    Every path from $a_4$ to $V(S_2 \cup S_3 \cup T) \setminus \{s_1,a_1,a_2\}$ in $G - \{s_1,a_1,a_2\}$ meets $S_1(s_1,a_1)$. As $\{s_1,a_1,a_2,w\}$ is not a $4$-separator for any vertex $w$, Menger's theorem gives two disjoint $H$-paths from $S_1(s_1,a_1)$ to $V(S_2 \cup S_3 \cup T) \setminus \{s_1,a_1,a_2\}$. By the above and the choice of $P(v_1,v_2)$, they end in $S_3(s_1,v_2]$, so one of them, $P(v_3,v_4)$, has $v_4 \in V(S_3(s_1,v_2))$.

    Since $\{s_1,a_1,a_2,v_2\}$ is not a $4$-separator, there is an $H$-path $P(v_5,v_6)$ with $v_5 \in V(S_1 \cup S_3(s_1,v_2) \cup S_4) \setminus \{s_1,a_1\}$ and $v_6 \in V(S_2 \cup S_3(v_2,a_3] \cup T) \setminus \{s_1,a_1,a_2,v_2\}$.
    By the above and the choice of $P(v_1,v_2)$, $v_5 \in V(S_3(s_1,v_2))$. Moreover, $P(v_5,v_6)$ is disjoint from $P(v_1,v_2)$, $P(v_3,v_4)$ and $P(a_4,b_4)$: otherwise, cutting it at the first common vertex from $v_6$ gives an $H$-path from $S_1$ or $S_4$ to $v_6$.

    By Lemma~\ref{S3toT}, applied after interchanging $a_1$ and $a_2$ to $P(a_4,b_4)$, $P(v_1,v_2)$, $P(v_3,v_4)$ and $P(v_5,v_6)$, $v_6 \notin V(T)$, and for the same reason (applied to $P[v_5,p] \cup P[p,v_0]$, where $p$ is the first vertex of $P(v_5,v_6)$ on $P(s_0,v_0)$) $P(v_5,v_6)$ is disjoint from $P(s_0,v_0)$.
    If $v_6 \in V(S_3)$, then replacing $S_3$ by $S_3[s_1,v_5] \cup P[v_5,v_6] \cup S_3[v_6,a_3]$ and $P(v_1,v_2)$ by $P[v_1,v_2] \cup S_3[v_2,v_6]$ contradicts the choice.
    Hence $v_6 \in V(S_2(s_1,a_2])$. Let $x$ be the vertex of $\{v_1,b_4\}$ closer to $s_1$ on $S_1$, and
    \[
    \begin{aligned}
        S_K ={}& S_1[x,a_1] \cup P[v_1,v_2] \cup S_3[v_5,a_3] \cup P[v_5,v_6] \cup S_2[v_6,a_2] \cup P[b_4,a_4], \\
        T_K ={}& T_1 \cup T_2 \cup T_3 \cup P[s_0,v_0] \cup S_4.
    \end{aligned}
    \]
    Then $K = S_K \cup T_K$ is of type $A$, as in Figure~\ref{fig:cross-s1-t}(b), a contradiction.

\end{proof}

\section{Finding a $K_6^-$-Minor}\label{sec:k6minus}

In this section we derive the $K_6^-$-minor statement. The $K_4^-$ case splits into the framed case, where the frame structure theorem applies, and the no-frame case, which is handled by a planar argument. The remaining $K_4^-$-free case is handled by induction.

\begin{theorem}\label{thm:k6minus}
    Let $G$ be a non-planar $5$-connected graph.
    Then $G$ contains a $K_6^-$-minor.
\end{theorem}

To prove this, we first consider the case where $G$ contains a $K_4^-$. We label its vertices as $a_1,a_2,a_3,a_4$ so that the omitted edge of the chosen $K_4^-$ subgraph is $a_3a_4$; this edge may still belong to $G$.

\subsection{The framed $K_4^-$ case}

\begin{lemma}[The framed $K_4^-$ case]\label{K4minusgivesK6minus}
    Let $G$ be a $5$-connected graph, and let $a_1,a_2,a_3,a_4 \in V(G)$ be distinct.
    Suppose that $a_1a_2, a_1a_3, a_1a_4, a_2a_3, a_2a_4 \in E(G)$, and that $G-a_3$ has no planar embedding with $a_1,a_2,a_4$ on the outer face or $G-a_4$ has no planar embedding with $a_1,a_2,a_3$ on the outer face.
    Then $G$ contains a $K_6^-$-minor.
\end{lemma}

\begin{proof}
    As the hypotheses are symmetric in $a_3$ and $a_4$, Lemma~\ref{lem:legless-tripod-to-frame} gives, after possibly exchanging $a_3$ and $a_4$, a pointed frame
    $(H,a_3,a_4)$ in $G$, rooted at $\{a_1,a_2,a_3,a_4\}$.

    By Theorem~\ref{structure}, $G$ contains a subgraph $K$ of one of types
    $A$, $B$, $C$, $D$, or $E$. Here $a_3a_4$ is the omitted edge of the chosen $K_4^-$ subgraph.
    We now show that each type yields a $K_6^-$-minor. The corresponding
    branch-set models are collected in Figure~\ref{fig:k6-branch-sets}.

    \begin{itemize}
        \item If $K = S \cup T$ is of type $A$, then since $\{a_1,a_2,a_3,a_4\}$ is not a $4$-separator, there is a $K$-path $P(a,b)$ with $a \in V(S) \setminus \{a_1,a_2,a_3,a_4\}$ and $b \in V(T) \setminus \{a_1,a_2,a_3,a_4\}$. Then
        \[
        \begin{aligned}
            \{&(S_0 \cup S_1 \cup S_2 \cup S_3 \cup S_4) \setminus \{a_1,a_2,a_3,a_4\}, \\
            &(T_0 \cup T_1 \cup T_2 \cup T_3 \cup T_4 \cup P(a,b]) \setminus \{a_1,a_2,a_3,a_4\}, a_1, a_2, a_3, a_4\}
        \end{aligned}
        \]
        is a $K_6^-$-minor (see Figure~\ref{fig:k6-branch-sets}(a)).

        \item If $K = S \cup T \cup U$ is of type $B$, then
        \[
        \begin{aligned}
            \{&S_0 \cup S_1[s_{i_1},a_1) \cup S_2[s_{i_2},a_2)\cup S_3[s_{i_3},a'_3) \cup S_4[s_{i_4},u_1) \cup U_2[u_2,v_2), \\
            &T_2[t,a_2) \cup T_3[t,w) \cup T_4[t,a_4), a_1, a_2, T_3[w,a'_3] \cup S_3[a'_3,a_3] \cup U_1[v_1,u_1), S_4[u_1,a_4]\}
        \end{aligned}
        \]
        is a $K_6^-$-minor, with $w$ as in the definition (see Figure~\ref{fig:k6-branch-sets}(b)).

        \item If $K = S \cup T \cup U$ is of type $C$, we choose it with $|E(S_2[a'_2,a_2])|$ minimal. Since $\{a_1,a'_2,a_3\}$ is not a $3$-separator, there is a $K$-path $P(a,b)$ with $a \in V(T) \setminus \{a_1,a'_2,a_3\}$ and $b \in V(S \cup U) \setminus \{a_1,a'_2,a_3\}$.
        If $b \in V(S_0 \cup S_1 \cup S_2[s_{i_2},a'_2) \cup S_3[s_{i_3},u) \cup S_4) \setminus \{a_1\}$, then
        \[
        \begin{aligned}
            \{&S_0 \cup S_1[s_{i_1},a_1) \cup S_2[s_{i_2},a'_2) \cup S_3[s_{i_3},u) \cup S_4[s_{i_4},v), \\
            &(T \cup P[a,b)) \setminus \{a_1,a'_2,a_3\}, a_1, S_2[a'_2,a_2], S_3[u,a_3] \cup U[u,v), S_4[v,a_4]\}
        \end{aligned}
        \]
        is a $K_6^-$-minor (see Figure~\ref{fig:k6-branch-sets}(c)).
        If $b \in V(U) \setminus \{u,v\}$, then
        \[
        \begin{aligned}
            \{&S_0 \cup S_1[s_{i_1},a_1) \cup S_2[s_{i_2},a'_2) \cup S_3[s_{i_3},u) \cup S_4[s_{i_4},v), \\
            &(T \cup P[a,b)) \setminus \{a_1,a'_2,a_3\}, a_1, S_2[a'_2,a_2], S_3[u,a_3] \cup U[u,b), U[b,v] \cup S_4(v,a_4]\}
        \end{aligned}
        \]
        is a $K_6^-$-minor, whose only missing adjacency is between the first two sets (see Figure~\ref{fig:k6-branch-sets}(d)).
        If $b \in V(S_2(a'_2,a_2])$, then replacing $T$ by the subtree of $T \cup P(a,b)$ spanning $a_1$, $b$ and $a_3$ contradicts the choice of $K$.
        Hence $b \in V(S_3[u,a_3))$, and we choose $K$ and $P(a,b)$ with, in addition, $|E(S_3[u,b])|$ minimal.

        Since $\{a_1,a_2,a'_2,b\}$ is not a $4$-separator, there is a $K$-path $P(e,f)$ with $e \in V(T \cup S_3(b,a_3]) \setminus \{a_1,a'_2\}$ and $f \in V(S \cup U) \setminus (V(S_3(b,a_3]) \cup \{a_1,a_2,a'_2,b\})$.
        By the above and the choice of $P(a,b)$, $e \in V(S_3(b,a_3])$, and $P(e,f)$ is disjoint from $P(a,b)$: otherwise, with $y$ the first vertex of $P(a,b)$ from $a$ on $P(e,f)$, $P[a,y] \cup P[y,f]$ is excluded in the same way.
        If $f \in V(U - u) \cup V(S_4[v,a_4])$, then $S \cup T \cup P(e,f) \cup U'$, where $U' = U[f,v]$ if $f \in V(U - u)$ and $U' = \emptyset$ if $f \in V(S_4[v,a_4])$, is of type $C$ with $e$ in the role of $u$, and $b \in V(S_3(s_{i_3},e))$, so the first case applies.
        If $f \in V(S_3[u,b))$, then replacing $S_3$ by $S_3[s_{i_3},f] \cup P[f,e] \cup S_3[e,a_3]$ gives a subgraph of type $C$ in which $P[a,b] \cup S_3[b,f]$ contradicts the choice of $P(a,b)$.
        If $f \in V(S_3(s_{i_3},u))$, then replacing $S_3$ by $S_3[s_{i_3},f] \cup P[f,e] \cup S_3[e,a_3]$ and $U$ by $S_3[f,u] \cup U$ gives a subgraph of type $C$ in which $P[a,b] \cup S_3[b,u]$ ends at an interior vertex of the new $U$, so the second case applies.
        Therefore $f \in V(S_0 \cup S_1 \cup S_2 \cup S_4[s_{i_4},v)) \setminus \{a_1,a_2,a'_2\}$.
        Then
        \[
        \begin{aligned}
            \{&S_0 \cup S_1[s_{i_1},a_1) \cup S_2[s_{i_2},a_2)\cup S_3[s_{i_3},u) \cup S_4[s_{i_4},v) \cup P[f,e), \\
            &T_1[t,a_1) \cup T_2[t,a'_2) \cup P[a,b) \cup T_3[t,a_3), a_1, a_2, S_3(b,a_3], S_3[u,b] \cup U \cup S_4[v,a_4]\}
        \end{aligned}
        \]
        is a $K_6^-$-minor, whose only possibly missing adjacency is between the second set and $a_2$ if $a'_2 \neq a_2$, and between the first two sets if $a'_2 = a_2$ (see Figure~\ref{fig:k6-branch-sets}(e)).

        \item If $K = S \cup T \cup U \cup V$ is of type $D$, we add the common subpath of $V_1$ and $V_2$, if any, to the first set below, that of $V_1$ and $V_3$ to the second, and that of $V_2$ and $V_3$ to the third, and we shorten $V_1$, $V_2$, $V_3$ accordingly. Each shortened path still joins the two sets it is meant to join, and two shortened paths now meet at most in a common end, which lies in the set receiving it. Hence the parts of $V_1$, $V_2$, $V_3$ outside the three enlarged sets are pairwise disjoint, which is all we use below.

        Then
        \[
        \begin{aligned}
            \{&S_2[s,u_2) \cup S_1[s,a_1) \cup S_4[s,a_4)
                \cup V_1[v_1,v'_1) \cup V_2[v_2,v'_2), \\
            &T_2[t,a_2) \cup T_3[t,u_3) \cup T_4[t,a_4)
                \cup V_3[v_3,v'_3), \\
            &U_2[u,u_2) \cup U_3[u,u_3),
                S_2[u_2,a_2], \\
            &\{a_1\} \cup T_3[u_3,a_3], a_4\}
        \end{aligned}
        \]
        is a $K_6^-$-minor (see Figure~\ref{fig:k6-branch-sets}(f)).

        \item If $K = S \cup T \cup U \cup V$ is of type $E$, then, with the convention of the definition for $S_r[s,u_r)$ and $T_r[t,u_r)$, let
        \[
            \begin{aligned}
                S' ={}& S_1[s,a_1) \cup S_2[s,u_2) \cup S_4[s,u_4), \\
                T' ={}& T_2[t,u_2) \cup T_3[t,u_3) \cup T_4[t,u_4), \\
                U' ={}& U - \{u_2,u_3,u_4\},
            \end{aligned}
        \]
        and add the interior of each of $V_1$, $V_2$ to one of the two sets it joins. For $r \in \{2,3,4\}$, let $A_r$ be the subpath of $S_r$ or $T_r$ from $u_r$ to $a_r$. Then
        \[
            \{S',\ T',\ U',\ A_2,\ \{a_1\} \cup A_3,\ A_4\}
        \]
        is a $K_6^-$-minor, whose only possibly missing adjacency is the pair among $S'$, $T'$, $U'$ not joined by $V_1$ or $V_2$ (see Figure~\ref{fig:k6-branch-sets}(g)). Indeed, these six sets are connected and pairwise disjoint; $S'$ is adjacent to $a_1$, $A_2$ and $A_4$ along $S_1$, $S_2$, $S_4$; $T'$ is adjacent to $A_2$, $A_3$, $A_4$ along $T_2$, $T_3$, $T_4$; $U'$ is adjacent to $u_2$, $u_3$, $u_4$; the edge $a_1a_3$ makes $\{a_1\}\cup A_3$ connected; and the edges $a_1a_2$, $a_2a_4$, $a_1a_4$ make $A_2$, $\{a_1\}\cup A_3$, $A_4$ pairwise adjacent.

    \end{itemize}

    \begin{figure}[tbp]
        \centering
\begin{minipage}[t]{0.32\linewidth}\centering\adjustbox{max width=\linewidth,max height=4.2cm}{%
    \begingroup\setlength{\textwidth}{\linewidth}\newcommand{\typeAFrame}[3]{%
\begin{tikzpicture}[scale=0.73, line cap=round, line join=round]
    \coordinate (s) at (-2.15,4.4);
    \coordinate (sp) at (2.15,4.4);
    \coordinate (a1) at (-3.2,0);
    \coordinate (a2) at (-1.1,0);
    \coordinate (a3) at (1.1,0);
    \coordinate (a4) at (3.2,0);

    \coordinate (t) at (-2.15,-4.4);
    \coordinate (tp) at (2.15,-4.4);

    \coordinate (a) at ($(s)!0.5!(sp)$);
    \coordinate (b) at ($(t)!0.5!(tp)$);

    \draw[Spath] (s)--(sp);
    \draw[Spath] (s)--(a1);
    \draw[Spath] (s)--(a2);
    \draw[Spath] (sp)--(a3);
    \draw[Spath] (sp)--(a4);

    \draw[Hpath] (a) to [bend left=6](b);

    \draw[rootedge] (a1) to [bend left=36](a2);
    \draw[rootedge] (a1) to [bend left=36](a3);
    \draw[rootedge] (a1) to [bend left=36](a4);
    \draw[rootedge] (a2) to [bend left=36](a3);
    \draw[rootedge] (a2) to [bend left=36](a4);

    \draw[Tpath] (t)--(tp);
    \draw[Tpath] (a1)--(t);
    \draw[Tpath] (a2)--(t);
    \draw[Tpath] (a3)--(tp);
    \draw[Tpath] (a4)--(tp);

\tikzset{
        setAHL/.style={draw=green!60!black, opacity=0.30, line width=7pt, line cap=round, line join=round},
        setBHL/.style={draw=black, opacity=0.24, line width=7pt, line cap=round, line join=round},
        setCHL/.style={draw=red!70!black, opacity=0.24, line width=7pt, line cap=round, line join=round},
        baseLbl/.style={font=\small, fill=white, fill opacity=0.88, text opacity=1, inner sep=1.5pt, rounded corners=1pt},
        mainLbl/.style={baseLbl, text=black},
        mainPlainLbl/.style={font=\small, text=black, inner sep=1.5pt},
        redLbl/.style={baseLbl, text=red!70!black},
        redPlainLbl/.style={font=\small, text=red!70!black, inner sep=1.5pt},
        greenLbl/.style={baseLbl, text=green!60!black}
    }
    #1

    \draw[setCHL] (s) -- (sp);
    \draw[setCHL] (t) -- (tp);

    \draw[setCHL] (a1) circle (0.13);
    \draw[setCHL] (a2) circle (0.13);
    \draw[setCHL] (a3) circle (0.13);
    \draw[setCHL] (a4) circle (0.13);

    \foreach \p in {a1,a2,a3,a4}
        \fill[white] (\p) circle (0.08);
    \foreach \p in {s,sp}
        \fill[black] (\p) circle (0.08);
    \foreach \p in {t,tp}
        \fill[green!60!black] (\p) circle (0.09);
    \foreach \p in {s,sp,a1,a2,a3,a4}
        \draw[black, thick] (\p) circle (0.08);
    \fill[red!70!black] (a) circle (0.09);
    \fill[red!70!black] (b) circle (0.09);

    #2

    \node[mainLbl, above=4pt] at (s) {$s_1$};
    \node[mainLbl, above=4pt] at (sp) {$s_2$};
    \node[mainLbl, left=6pt] at (a1) {$a_1$};
    \node[mainLbl, left=6pt] at (a2) {$a_2$};
    \node[mainLbl, right=6pt] at (a3) {$a_3$};
    \node[mainLbl, right=6pt] at (a4) {$a_4$};
    \node[mainLbl, text=green!60!black, below=4pt] at (t) {$t_1$};
    \node[mainLbl, text=green!60!black, below=4pt] at (tp) {$t_2$};
    \node[redPlainLbl, above=4pt] at (a) {$a$};
    \node[redPlainLbl, below=4pt] at (b) {$b$};

    #3
\end{tikzpicture}%
}

\typeAFrame{}{}{}\endgroup}\par\small\textbf{(a)} type $A$\end{minipage}\hfill\begin{minipage}[t]{0.32\linewidth}\centering\adjustbox{max width=\linewidth,max height=4.2cm}{%
    \begingroup\setlength{\textwidth}{\linewidth}\begin{tikzpicture}[scale=0.73, line cap=round, line join=round]
        \coordinate (s) at (-2.15,4.4);
        \coordinate (sp) at (2.15,4.4);
        \coordinate (a1) at (-3.2,0);
        \coordinate (a2) at (-1.1,0);
        \coordinate (a3) at (1.1,0);
        \coordinate (a4) at (3.2,0);

        \coordinate (t) at (2.15,-4.4);
        \coordinate (ap3) at ($(sp)!0.5!(a3)$);

        \coordinate (u1) at ($(sp)!0.65!(a4)$);
        \coordinate (u2) at ($(s)!0.30!(sp)$);
        \coordinate (v2) at ($(t)!0.8!(a2)$);

        \draw[rootedge] (a1) to [bend left=36](a2);
        \draw[rootedge] (a1) to [bend left=36](a3);
        \draw[rootedge] (a1) to [bend left=36](a4);
        \draw[rootedge] (a2) to [bend left=36](a3);
        \draw[rootedge] (a2) to [bend left=36](a4);

        \draw[Spath] (s)--(sp);
        \draw[Spath] (s)--(a1);
        \draw[Spath] (s)--(a2);
        \draw[Spath] (sp)--(a3);
        \draw[Spath] (sp)--(a4);

        \draw[Tpath] (a2)--(t);
        \draw[Tpath] (ap3) to[bend left=18] node[pos=0.4, coordinate] (v1) {} (t);
        \draw[Tpath] (a4)--(t);

        \draw[Hpath] (u1) to[bend right=18] (v1);
        \draw[Hpath] (u2) to[bend left=12] (v2);

        \tikzset{
            setAHL/.style={draw=green!60!black, opacity=0.30, line width=7pt, line cap=round, line join=round},
            setBHL/.style={draw=black, opacity=0.24, line width=7pt, line cap=round, line join=round},
            setCHL/.style={draw=red!70!black, opacity=0.24, line width=7pt, line cap=round, line join=round},
            baseLbl/.style={font=\small, fill=white, fill opacity=0.88, text opacity=1, inner sep=1.5pt, rounded corners=1pt},
            mainLbl/.style={baseLbl, text=black},
            mainPlainLbl/.style={font=\small, text=black, inner sep=1.5pt},
            redLbl/.style={baseLbl, text=red!70!black},
            redPlainLbl/.style={font=\small, text=red!70!black, inner sep=1.5pt},
            greenLbl/.style={baseLbl, text=green!60!black}
        }

        \draw[setCHL] (a1) circle (0.13);
        \draw[setCHL] (a2) circle (0.13);
        \draw[setCHL] (s) -- (sp);
        \draw[setCHL] (t) -- (v2);
        \draw[setCHL] (u1) -- (a4);
        \draw[setCHL] (ap3) to[bend left=18] (v1);
        \draw[setCHL] (ap3) -- (a3);

        \foreach \p in {s,sp,a1,a2,a3,a4,t}
            \fill[white] (\p) circle (0.08);
        \foreach \p in {s,sp}
            \fill[black] (\p) circle (0.08);
        \fill[green!60!black] (t) circle (0.09);
        \fill[green!60!black] (ap3) circle (0.09);

        \foreach \p in {s,sp,a1,a2,a3,a4}
            \draw[black, thick] (\p) circle (0.08);
        \foreach \p in {u1,v1}
            \fill[red!70!black] (\p) circle (0.09);
        \foreach \p in {u2,v2}
            \fill[red!70!black] (\p) circle (0.09);

        \node[mainLbl, above=4pt] at (s) {$s_1$};
        \node[mainLbl, above=4pt] at (sp) {$s_2$};
        \node[mainLbl, left=6pt] at (a1) {$a_1$};
        \node[mainLbl, left=6pt] at (a2) {$a_2$};
        \node[mainLbl, left=6pt] at (a3) {$a_3$};
        \node[mainLbl, right=6pt] at (a4) {$a_4$};
        \node[mainLbl, text=green!60!black, below=4pt] at (t) {$t$};
        \node[mainLbl, left=4pt] at (ap3) {$a'_3$};

        \node[mainLbl, right=4pt] at (u1) {$u_1$};
        \node[mainLbl, left=4pt] at (v1) {$v_1$};

\end{tikzpicture}\endgroup}\par\small\textbf{(b)} type $B$\end{minipage}\hfill\begin{minipage}[t]{0.32\linewidth}\centering\adjustbox{max width=\linewidth,max height=4.2cm}{%
    \begingroup\setlength{\textwidth}{\linewidth}\begin{tikzpicture}[scale=0.73, line cap=round, line join=round]
    \coordinate (s) at (-2.15,4.4);
    \coordinate (sp) at (2.15,4.4);
    \coordinate (a1) at (-3.2,0);
    \coordinate (a2) at (-1.1,0);
    \coordinate (a3) at (1.1,0);
    \coordinate (a4) at (3.2,0);
    \coordinate (t) at (-1.2,-3.0);

    \coordinate (u3) at ($(sp)!0.5!(a3)$);
    \coordinate (ap2) at ($(s)!0.5!(a2)$);

    \coordinate (a) at ($(t)!0.6!(a3)$);
    \coordinate (b) at ($(s)!0.25!(sp)$);

    \coordinate (u4) at ($(sp)!0.5!(a4)$);

    \draw[Spath] (s)--(sp);
    \draw[Spath] (s)--(a1);
    \draw[Spath] (s)--(a2);
    \draw[Spath] (sp)--(a3);
    \draw[Spath] (sp)--(a4);

    \draw[Tpath] (a1)--(t);
    \draw[Tpath] (ap2) to[bend right = 18](t);
    \draw[Tpath] (a3)--(t);

    \draw[Hpath] (u4) to[bend right=18] (u3);
    \draw[Hpath] (a) to[bend right=18] (b);

    \draw[rootedge] (a1) to [bend left=36](a2);
    \draw[rootedge] (a1) to [bend left=36](a3);
    \draw[rootedge] (a1) to [bend left=36](a4);
    \draw[rootedge] (a2) to [bend left=36](a3);
    \draw[rootedge] (a2) to [bend left=36](a4);

        \tikzset{
            setAHL/.style={draw=green!60!black, opacity=0.30, line width=7pt, line cap=round, line join=round},
            setBHL/.style={draw=black, opacity=0.24, line width=7pt, line cap=round, line join=round},
            setCHL/.style={draw=red!70!black, opacity=0.24, line width=7pt, line cap=round, line join=round},
            baseLbl/.style={font=\small, fill=white, fill opacity=0.88, text opacity=1, inner sep=1.5pt, rounded corners=1pt},
            mainLbl/.style={baseLbl, text=black},
            mainPlainLbl/.style={font=\small, text=black, inner sep=1.5pt},
            redLbl/.style={baseLbl, text=red!70!black},
            redPlainLbl/.style={font=\small, text=red!70!black, inner sep=1.5pt},
            greenLbl/.style={baseLbl, text=green!60!black}
        }

        \draw[setCHL] (a1) circle (0.13);
        \draw[setCHL] (u4) -- (a4);
        \draw[setCHL] (s) -- (sp);
        \draw[setCHL] (t) -- (a);
        \draw[setCHL] (ap2) -- (a2);
        \draw[setCHL] (u3) -- (a3);

    \foreach \p in {a1,a2,a3,a4}
        \fill[white] (\p) circle (0.08);
    \foreach \p in {s,sp}
        \fill[black] (\p) circle (0.08);
    \fill[green!60!black] (t) circle (0.09);
    \foreach \p in {s,sp,a1,a2,a3,a4}
        \draw[black, thick] (\p) circle (0.08);
    \fill[red!70!black] (u3) circle (0.09);
    \fill[red!70!black] (a) circle (0.09);
    \fill[red!70!black] (b) circle (0.09);

    \fill[green!60!black] (ap2) circle (0.09);

    \fill[red!70!black] (u4) circle (0.09);

    \node[mainLbl, above=4pt] at (s) {$s_1$};
    \node[mainLbl, above=4pt] at (sp) {$s_2$};
    \node[mainLbl, left=6pt] at (a1) {$a_1$};
    \node[mainLbl, right=6pt] at (a2) {$a_2$};
    \node[mainLbl, right=6pt] at (a3) {$a_3$};
    \node[mainLbl, right=6pt] at (a4) {$a_4$};
    \node[mainLbl, text=green!60!black, below=4pt] at (t) {$t$};
    \node[mainLbl, right=4pt] at (ap2) {$a'_2$};
    \node[redLbl, right=4pt] at (a) {$a$};
    \node[redLbl, above=4pt] at (b) {$b$};

    \node[redLbl, left=4pt] at (u3) {$u$};
    \node[redLbl, right=4pt] at (u4) {$v$};

\end{tikzpicture}
\par\smallskip
{}\endgroup}\par\small\textbf{(c)} type $C$, initial case\end{minipage}

\vspace{0.8em}

\begin{minipage}[t]{0.32\linewidth}\centering\adjustbox{max width=\linewidth,max height=4.2cm}{%
    \begingroup\setlength{\textwidth}{\linewidth}\begin{tikzpicture}[scale=0.73, line cap=round, line join=round]
    \coordinate (s) at (-2.15,4.4);
    \coordinate (sp) at (2.15,4.4);
    \coordinate (a1) at (-3.2,0);
    \coordinate (a2) at (-1.1,0);
    \coordinate (a3) at (1.1,0);
    \coordinate (a4) at (3.2,0);
    \coordinate (t) at (-1.2,-3.0);

    \coordinate (u3) at ($(sp)!0.5!(a3)$);
    \coordinate (ap2) at ($(s)!0.5!(a2)$);

    \coordinate (a) at ($(t)!0.2!(a3)$);

    \coordinate (u4) at ($(sp)!0.5!(a4)$);

    \draw[Spath] (s)--(sp);
    \draw[Spath] (s)--(a1);
    \draw[Spath] (s)--(a2);
    \draw[Spath] (sp)--(a3);
    \draw[Spath] (sp)--(a4);

    \draw[Tpath] (a1)--(t);
    \draw[Tpath] (ap2) to[bend right = 18](t);
    \draw[Tpath] (a3)--(t);

    \draw[Hpath] (u4) to[bend right=18] node[pos=0.5, coordinate] (b) {} (u3);
    \draw[Hpath] (a) to[bend right=36] (b);

    \draw[rootedge] (a1) to [bend left=36](a2);
    \draw[rootedge] (a1) to [bend left=36](a3);
    \draw[rootedge] (a1) to [bend left=36](a4);
    \draw[rootedge] (a2) to [bend left=36](a3);
    \draw[rootedge] (a2) to [bend left=36](a4);

        \tikzset{
            setAHL/.style={draw=green!60!black, opacity=0.30, line width=7pt, line cap=round, line join=round},
            setBHL/.style={draw=black, opacity=0.24, line width=7pt, line cap=round, line join=round},
            setCHL/.style={draw=red!70!black, opacity=0.24, line width=7pt, line cap=round, line join=round},
            baseLbl/.style={font=\small, fill=white, fill opacity=0.88, text opacity=1, inner sep=1.5pt, rounded corners=1pt},
            mainLbl/.style={baseLbl, text=black},
            mainPlainLbl/.style={font=\small, text=black, inner sep=1.5pt},
            redLbl/.style={baseLbl, text=red!70!black},
            redPlainLbl/.style={font=\small, text=red!70!black, inner sep=1.5pt},
            greenLbl/.style={baseLbl, text=green!60!black}
        }

        \draw[setCHL] (a1) circle (0.13);
        \draw[setCHL] (u4) -- (a4);
        \draw[setCHL] (s) -- (sp);
        \draw[setCHL] (t) -- (a);
        \draw[setCHL] (ap2) -- (a2);
        \draw[setCHL] (u3) -- (a3);
        \draw[setCHL] (u4) to[bend right=18] (b);

    \foreach \p in {a1,a2,a3,a4}
        \fill[white] (\p) circle (0.08);
    \foreach \p in {s,sp}
        \fill[black] (\p) circle (0.08);
    \fill[green!60!black] (t) circle (0.09);
    \foreach \p in {s,sp,a1,a2,a3,a4}
        \draw[black, thick] (\p) circle (0.08);
    \fill[red!70!black] (u3) circle (0.09);
    \fill[red!70!black] (a) circle (0.09);
    \fill[red!70!black] (b) circle (0.09);

    \fill[green!60!black] (ap2) circle (0.09);

    \fill[red!70!black] (u4) circle (0.09);

    \node[mainLbl, above=4pt] at (s) {$s_1$};
    \node[mainLbl, above=4pt] at (sp) {$s_2$};
    \node[mainLbl, left=6pt] at (a1) {$a_1$};
    \node[mainLbl, right=6pt] at (a2) {$a_2$};
    \node[mainLbl, left=6pt] at (a3) {$a_3$};
    \node[mainLbl, right=6pt] at (a4) {$a_4$};
    \node[mainLbl, text=green!60!black, below=4pt] at (t) {$t$};
    \node[mainLbl, right=4pt] at (ap2) {$a'_2$};
    \node[redLbl, right=4pt] at (a) {$a$};
    \node[redLbl, above=4pt] at (b) {$b$};

    \node[redLbl, left=4pt] at (u3) {$u$};
    \node[redLbl, right=4pt] at (u4) {$v$};

\end{tikzpicture}
\par\smallskip
{}\endgroup}\par\small\textbf{(d)} type $C$, $b \in V(U)$\end{minipage}\hfill\begin{minipage}[t]{0.32\linewidth}\centering\adjustbox{max width=\linewidth,max height=4.2cm}{%
    \begingroup\setlength{\textwidth}{\linewidth}\begin{tikzpicture}[scale=0.73, line cap=round, line join=round]
    \coordinate (s) at (-2.15,4.4);
    \coordinate (sp) at (2.15,4.4);
    \coordinate (a1) at (-3.2,0);
    \coordinate (a2) at (-1.1,0);
    \coordinate (a3) at (1.1,0);
    \coordinate (a4) at (3.2,0);
    \coordinate (t) at (-1.2,-3.0);

    \coordinate (u3) at ($(sp)!0.3!(a3)$);
    \coordinate (ap2) at ($(s)!0.6!(a2)$);

    \coordinate (a) at ($(t)!0.6!(a3)$);
    \coordinate (b) at ($(sp)!0.45!(a3)$);
    \coordinate (e) at ($(sp)!0.6!(a3)$);
    \coordinate (f) at ($(s)!0.4!(a2)$);
    
    \coordinate (u4) at ($(sp)!0.3!(a4)$);

    \draw[Spath] (s)--(sp);
    \draw[Spath] (s)--(a1);
    \draw[Spath] (s)--(a2);
    \draw[Spath] (sp)--(a3);
    \draw[Spath] (sp)--(a4);

    \draw[Tpath] (a1)--(t);
    \draw[Tpath] (ap2) to[bend right = 18](t);
    \draw[Tpath] (a3)--(t);

    \draw[Hpath] (u4) to[bend right=18] (u3);
    \draw[Hpath] (a) to[bend left=18] (b);
    \draw[Hpath] (e) to[bend left=18] (f);

    \draw[rootedge] (a1) to [bend left=36](a2);
    \draw[rootedge] (a1) to [bend left=36](a3);
    \draw[rootedge] (a1) to [bend left=36](a4);
    \draw[rootedge] (a2) to [bend left=36](a3);
    \draw[rootedge] (a2) to [bend left=36](a4);

        \tikzset{
            setAHL/.style={draw=green!60!black, opacity=0.30, line width=7pt, line cap=round, line join=round},
            setBHL/.style={draw=black, opacity=0.24, line width=7pt, line cap=round, line join=round},
            setCHL/.style={draw=red!70!black, opacity=0.24, line width=7pt, line cap=round, line join=round},
            baseLbl/.style={font=\small, fill=white, fill opacity=0.88, text opacity=1, inner sep=1.5pt, rounded corners=1pt},
            mainLbl/.style={baseLbl, text=black},
            mainPlainLbl/.style={font=\small, text=black, inner sep=1.5pt},
            redLbl/.style={baseLbl, text=red!70!black},
            redPlainLbl/.style={font=\small, text=red!70!black, inner sep=1.5pt},
            greenLbl/.style={baseLbl, text=green!60!black}
        }

        \draw[setCHL] (a1) circle (0.13);
        \draw[setCHL] (a2) circle (0.13);
        \draw[setCHL] (s) -- (sp);
        \draw[setCHL] (s) -- (ap2);
        \draw[setCHL] (t) -- (a);
        \draw[setCHL] (b) -- (u3);
        \draw[setCHL] (u4) -- (a4);
        \draw[setCHL] (u3) to[bend left=18] (u4);
        \draw[setCHL] (e) -- (a3);

    \foreach \p in {a1,a2,a3,a4}
        \fill[white] (\p) circle (0.08);
    \foreach \p in {s,sp}
        \fill[black] (\p) circle (0.08);
    \fill[green!60!black] (t) circle (0.09);
    \foreach \p in {s,sp,a1,a2,a3,a4}
        \draw[black, thick] (\p) circle (0.08);
    \fill[red!70!black] (u3) circle (0.09);
    \fill[red!70!black] (a) circle (0.09);
    \fill[red!70!black] (b) circle (0.09);
    \fill[red!70!black] (e) circle (0.09);
    \fill[red!70!black] (f) circle (0.09);

    \fill[green!60!black] (ap2) circle (0.09);

    \fill[red!70!black] (u4) circle (0.09);

    \node[mainLbl, above=4pt] at (s) {$s_1$};
    \node[mainLbl, above=4pt] at (sp) {$s_2$};
    \node[mainLbl, left=6pt] at (a1) {$a_1$};
    \node[mainLbl, right=6pt] at (a2) {$a_2$};
    \node[mainLbl, right=6pt] at (a3) {$a_3$};
    \node[mainLbl, right=6pt] at (a4) {$a_4$};
    \node[mainLbl, text=green!60!black, below=4pt] at (t) {$t$};
    \node[mainLbl, right=4pt] at (ap2) {$a'_2$};
    \node[redLbl, right=4pt] at (a) {$a$};
    \node[redLbl, right=4pt] at (b) {$b$};

    \node[redLbl, right=4pt] at (e) {$e$};
    \node[redLbl, left=4pt] at (f) {$f$};

    \node[redLbl, left=4pt] at (u3) {$u$};
    \node[redLbl, right=4pt] at (u4) {$v$};

\end{tikzpicture}
\par\smallskip
{}\endgroup}\par\small\textbf{(e)} type $C$, final case\end{minipage}\hfill\begin{minipage}[t]{0.32\linewidth}\centering\adjustbox{max width=\linewidth,max height=4.2cm}{%
    \begingroup\setlength{\textwidth}{\linewidth}\begin{tikzpicture}[scale=0.88, line cap=round, line join=round]
    \coordinate (s) at (0,4.6);
    \coordinate (a1) at (-3.2,0);
    \coordinate (a2) at (-0.9,0);
    \coordinate (a3) at (1.5,0);
    \coordinate (a4) at (3.8,0);
    \coordinate (t) at (2.65,-4.4);
    \coordinate (u) at (-3.2,-4.4);

    \coordinate (ui) at ($(s)!0.60!(a2)$);
    \coordinate (u3) at ($(t)!0.83!(a3)$);

    \coordinate (vs) at ($(s)!0.33!(a4)$);
    \coordinate (vsp) at ($(s)!0.24!(a2)$);
    \coordinate (vt) at ($(t)!0.26!(a4)$);
    \coordinate (vtp) at ($(t)!0.24!(a2)$);
    \coordinate (vu) at ($(u)!0.36!(ui)$);
    \coordinate (vup) at ($(u)!0.32!(u3)$);

        \draw[rootedge] (a1) to [bend left=36](a2);
    \draw[rootedge] (a1) to [bend left=36](a3);
    \draw[rootedge] (a1) to [bend left=36](a4);
    \draw[rootedge] (a2) to [bend left=36](a3);
    \draw[rootedge] (a2) to [bend left=36](a4);

        \tikzset{
            setAHL/.style={draw=green!60!black, opacity=0.30, line width=7pt, line cap=round, line join=round},
            setBHL/.style={draw=black, opacity=0.24, line width=7pt, line cap=round, line join=round},
            setCHL/.style={draw=red!70!black, opacity=0.24, line width=7pt, line cap=round, line join=round},
            baseLbl/.style={font=\large, fill=white, fill opacity=0.88, text opacity=1, inner sep=1.5pt, rounded corners=1pt},
            mainLbl/.style={baseLbl, text=black},
            mainPlainLbl/.style={font=\large, text=black, inner sep=1.5pt},
            redLbl/.style={baseLbl, text=red!70!black},
            redPlainLbl/.style={font=\large, text=red!70!black, inner sep=1.5pt},
            greenLbl/.style={baseLbl, text=green!60!black},
            blueLbl/.style={baseLbl, text=blue!70!black}
        }

    \draw[setCHL] (a4) circle (0.13);
    \draw[setCHL] (a1) circle (0.13);
    \draw[setCHL] (s) -- (vsp);
    \draw[setCHL] (s) -- (vs);
    \draw[setCHL] (t) -- (vt);
    \draw[setCHL] (t) -- (vtp);
    \draw[setCHL] (u) -- (vu);
    \draw[setCHL] (u) -- (vup);
    \draw[setCHL] (ui) -- (a2);
    \draw[setCHL] (u3) -- (a3);
    \draw[setCHL] (a1) to [bend left=36] (a3);

    \draw[Spath] (s)--(a2);
    \draw[Spath] (s)--(a1);
    \draw[Spath] (s)--(a4);

    \draw[Tpath] (t)--(a2);
    \draw[Tpath] (t)--(a3);
    \draw[Tpath] (t)--(a4);

    \draw[Upath] (u)--(ui);
    \draw[Upath] (u)--(u3);

    \draw[Hpath] (vs) to[bend left=20]  (vt);
    \draw[Hpath] (vsp) to[bend right=10]  (vu);
    \draw[Hpath] (vtp) to[bend left=12]  (vup);

    \foreach \p in {s,a1,a2,a3,a4,t}
        \fill[white] (\p) circle (0.08);
    \fill[black] (s) circle (0.08);
    \fill[green!60!black] (t) circle (0.09);
    \foreach \p in {s,a1,a2,a3,a4}
        \draw[black, thick] (\p) circle (0.08);
    \foreach \p in {u,ui,u3}
        \fill[blue!70!black] (\p) circle (0.09);
    \foreach \p in {vs,vsp,vt,vtp,vu,vup}
        \fill[red!70!black] (\p) circle (0.09);

    \node[mainLbl, above=4pt] at (s) {$s$};
    \node[mainLbl, right=7pt] at (a2) {$a_2$};
    \node[mainLbl, left=7pt] at (a1) {$a_1$};
    \node[mainLbl, right=7pt] at (a3) {$a_3$};
    \node[mainLbl, right=7pt] at (a4) {$a_4$};
    \node[mainLbl, text=green!60!black, below=4pt] at (t) {$t$};

    \node[blueLbl, below=3pt] at (u) {$u$};
    \node[blueLbl, right=3pt] at (ui) {$u_2$};
    \node[blueLbl, right=3pt] at (u3) {$u_3$};

    \node[redLbl, above right=2pt] at (vs) {$v_1$};
    \node[redLbl, right=2pt] at (vsp) {$v_2$};
    \node[redLbl, right=3pt] at (vt) {$v'_1$};
    \node[redLbl, below left=3pt] at (vtp) {$v_3$};
    \node[redLbl, above left=2pt] at (vu) {$v'_2$};
    \node[redLbl, below right = 4pt] at (vup) {$v'_3$};

\end{tikzpicture}\endgroup}\par\small\textbf{(f)} type $D$\end{minipage}

\vspace{0.8em}

\begin{minipage}[t]{0.32\linewidth}\centering\adjustbox{max width=\linewidth,max height=4.2cm}{%
    \begingroup\setlength{\textwidth}{\linewidth}\begin{tikzpicture}[scale=0.88, line cap=round, line join=round]
    \coordinate (s) at (0,4.6);
    \coordinate (aj) at (-3.2,0);
    \coordinate (ai) at (-0.9,0);
    \coordinate (a3) at (1.5,0);
    \coordinate (a4) at (3.8,0);
    \coordinate (t) at (2.65,-4.4);
    \coordinate (u) at (-3.2,-4.4);

    \coordinate (ui) at ($(t)!0.70!(ai)$);
    \coordinate (u4) at ($(t)!0.70!(a4)$);
    \coordinate (u3) at ($(t)!0.83!(a3)$);

    \coordinate (a) at ($(u)!0.5!(ui)$);
    \coordinate (b) at ($(s)!0.5!(ai)$);
    \coordinate (c) at ($(u)!0.5!(u4)$);
    \coordinate (d) at ($(t)!0.2!(ai)$);

    \draw[rootedge] (aj) to [bend left=36](ai);
    \draw[rootedge] (aj) to [bend left=36](a3);
    \draw[rootedge] (aj) to [bend left=36](a4);
    \draw[rootedge] (ai) to [bend left=36](a3);
    \draw[rootedge] (ai) to [bend left=36](a4);

    \tikzset{
            setAHL/.style={draw=green!60!black, opacity=0.30, line width=7pt, line cap=round, line join=round},
            setBHL/.style={draw=black, opacity=0.24, line width=7pt, line cap=round, line join=round},
            setCHL/.style={draw=red!70!black, opacity=0.24, line width=7pt, line cap=round, line join=round},
            baseLbl/.style={font=\large, fill=white, fill opacity=0.88, text opacity=1, inner sep=1.5pt, rounded corners=1pt},
            mainLbl/.style={baseLbl, text=black},
            mainPlainLbl/.style={font=\large, text=black, inner sep=1.5pt},
            redLbl/.style={baseLbl, text=red!70!black},
            redPlainLbl/.style={font=\large, text=red!70!black, inner sep=1.5pt},
            greenLbl/.style={baseLbl, text=green!60!black},
            blueLbl/.style={baseLbl, text=blue!70!black}
    }

    \draw[setCHL] (aj) to [bend left = 36] (a3);
    \draw[setCHL] (ai) -- (ui);
    \draw[setCHL] (a3) -- (u3);
    \draw[setCHL] (a4) -- (u4);
    \draw[setCHL] (s) -- (b);
    \draw[setCHL] (u) -- (a);
    \draw[setCHL] (u) -- (c);
    \draw[setCHL] (t) -- (d);

    \draw[Spath] (s)--(ai);
    \draw[Spath] (s)--(aj);
    \draw[Spath] (s)--(a4);

    \draw[Tpath] (t)--(ai);
    \draw[Tpath] (t)--(a3);
    \draw[Tpath] (t)--(a4);

    \draw[Upath] (u)--(ui);
    \draw[Upath] (u)--(u3);
    \draw[Upath] (u)--(u4);

    \draw[Hpath] (a) to [bend left=18] (b);
    \draw[Hpath] (c) to [bend right=18] (d);

    \foreach \p in {s,ai,aj,a3,a4,t}
        \fill[white] (\p) circle (0.08);
    \fill[black] (s) circle (0.08);
    \fill[green!60!black] (t) circle (0.09);
    \foreach \p in {s,ai,aj,a3,a4}
        \draw[black, thick] (\p) circle (0.08);
    \foreach \p in {u,ui,u3,u4}
        \fill[blue!70!black] (\p) circle (0.09);
    \foreach \p in {a,b,c,d}
        \fill[red!70!black] (\p) circle (0.09);

    \node[mainLbl, above=4pt] at (s) {$s$};
    \node[mainLbl, right=7pt] at (ai) {$a_2$};
    \node[mainLbl, left=7pt] at (aj) {$a_1$};
    \node[mainLbl, right=7pt] at (a3) {$a_3$};
    \node[mainLbl, right=7pt] at (a4) {$a_4$};
    \node[mainLbl, text=green!60!black, below=4pt] at (t) {$t$};

    \node[blueLbl, below=3pt] at (u) {$u$};
    \node[blueLbl, above right=3pt] at (ui) {$u_2$};
    \node[blueLbl, right=3pt] at (u3) {$u_3$};
    \node[blueLbl, right=4pt] at (u4) {$u_4$};

    \node[redLbl, above left=3pt] at (a) {$v'_1$};
    \node[redLbl, right=4pt] at (b) {$v_1$};
    \node[redLbl, below= 6pt] at (c) {$v'_2$};
    \node[redLbl, below left=3pt] at (d) {$v_2$};

\end{tikzpicture}\endgroup}\par\small\textbf{(g)} type $E$\end{minipage}
        \caption{Branch-set constructions for types $A$--$E$.}
        \label{fig:k6-branch-sets}
    \end{figure}
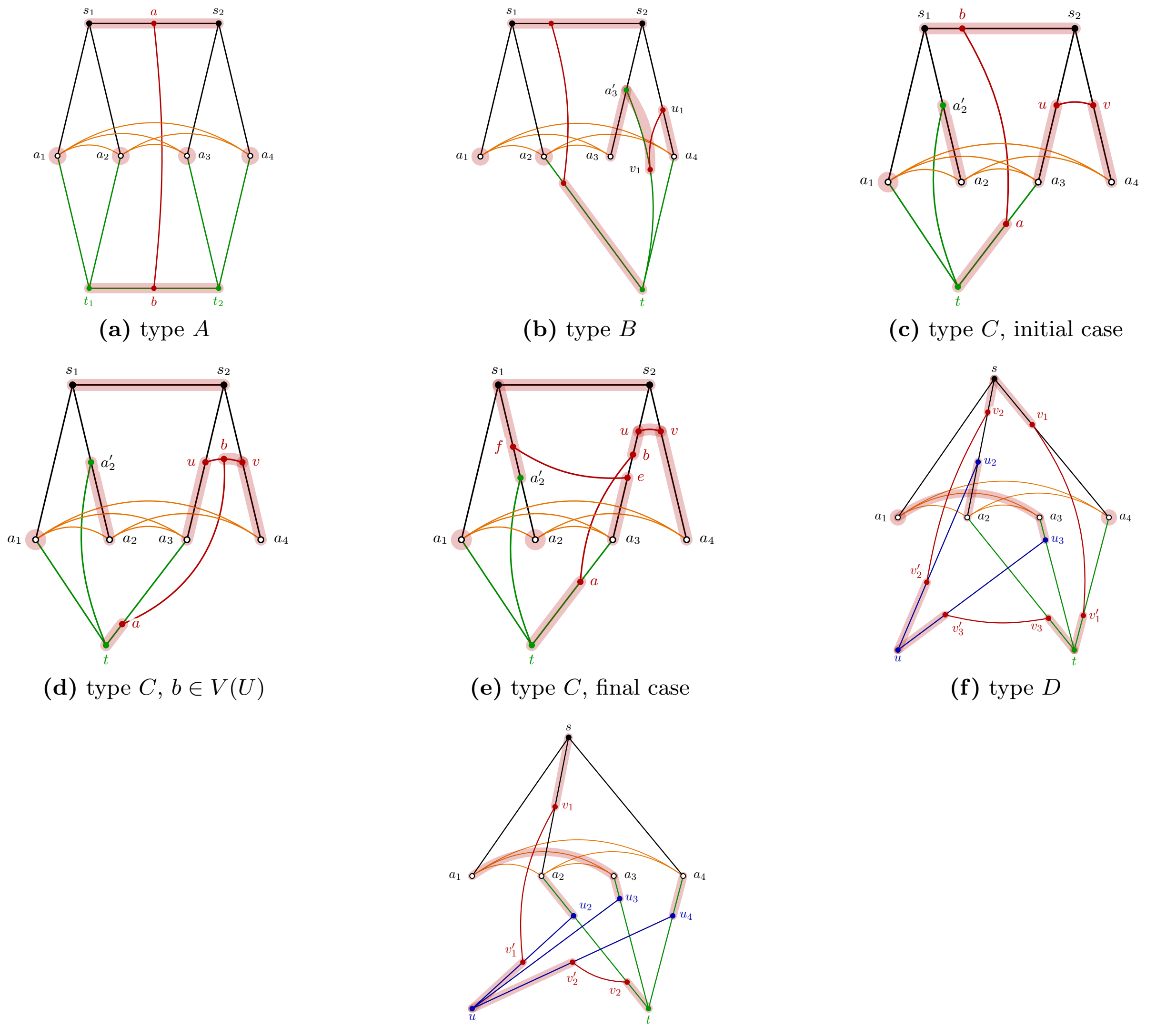

\end{proof}

\subsection{The no-frame $K_4^-$ case}

We now treat the complementary case where Lemma~\ref{lem:legless-tripod-to-frame} does not produce a frame. Thus both $G-a_3$ and $G-a_4$ have planar embeddings with the prescribed roots on the outer face.

\begin{lemma}[Whitney~\cite{Whitney1933}]\label{lem:whitney}
    Every $3$-connected planar graph has a unique planar embedding, up to
    reflection and the choice of the outer face. In particular its set of
    facial cycles is determined by the graph.
\end{lemma}

We also use the following standard form of Tutte's theorem on peripheral
cycles in $3$-connected planar graphs.

\begin{lemma}[Tutte's theorem~\cite{Tutte1963}]\label{lem:tutte-facial-cycles}
    Let $G$ be a $3$-connected plane graph, and let $C$ be a facial cycle of $G$.
    Then $C$ is induced and non-separating; equivalently, $G-V(C)$ is connected
    whenever it is nonempty.
\end{lemma}

\begin{remark}\label{rem:no-interleaving}
Let $C$ be a cycle. Two sets $A,B \subseteq V(C)$ \emph{interleave} on $C$ if there are no distinct vertices $p,q$ of $C$ such that $A$ and $B$ lie in the two paths of $C$ between $p$ and $q$, respectively; this holds in particular if there are $a,a' \in A$ and $b,b' \in B$, all distinct, in the cyclic order $a,b,a',b'$. If $C$ is drawn in the plane and $A,B$ are disjoint connected subgraphs drawn in the same closed region bounded by $C$ and disjoint from $C$, then $N(A)\cap V(C)$ and $N(B)\cap V(C)$ do not interleave.
\end{remark}

To treat a $K_4$, we use the following notation. For distinct vertices $x,y$ of a $5$-connected graph $G$ such that $G - y$ is planar, let $F^{y}_{x}$ denote the face of the plane graph $G - \{x,y\}$ whose boundary contains $N(x) \setminus \{y\}$. This face is well defined: $G - \{x,y\}$ is $3$-connected, hence uniquely embedded by Lemma~\ref{lem:whitney}; in an embedding of $G - y$ the vertex $x$ lies in such a face; and it is unique, since $|N(x) \setminus \{y\}| \ge 4$ while two faces of a $3$-connected plane graph share at most two vertices.

\begin{lemma}\label{lem:k4-pair}
    Let $G$ be a non-planar $5$-connected graph, and let $a_1,a_2,a_3,a_4$ induce a $K_4$ in $G$ with $G - a_i$ planar for every $i$. Then $F^{y}_{x} = F^{x}_{y}$ for some distinct $x,y \in \{a_1,a_2,a_3,a_4\}$.
\end{lemma}

\begin{proof}
    Suppose that $F_3 := F^{a_4}_{a_3}$ and $F_4 := F^{a_3}_{a_4}$ are distinct, and let $J := G - \{a_3,a_4\}$. By Lemma~\ref{lem:tutte-facial-cycles}, the boundaries of $F_3$ and $F_4$ are induced cycles of $J$; both contain $a_1$ and $a_2$, hence the edge $a_1a_2$. So $F_3$ and $F_4$ are the two faces of $J$ incident with $a_1a_2$. Drawing $a_3$ in $F_3$ and $a_4$ in $F_4$, deleting $a_1a_2$ and drawing $a_3a_4$ through the merged face gives a planar embedding $\mathcal{E}$ of $G - a_1a_2$.

    Let $J' := G - \{a_1,a_3\}$. Since $G - a_1 \subseteq G - a_1a_2$, the embedding $\mathcal{E}$ restricts to an embedding of $G - a_1$, so $a_3$ lies in $F^{a_1}_{a_3}$. In $\mathcal{E}$, the vertex $a_1$ lies in a face of $J'$ whose boundary contains $N(a_1) \setminus \{a_2,a_3\}$, a set of at least three vertices which is also contained in the boundary of $F^{a_3}_{a_1}$; hence this face is $F^{a_3}_{a_1}$. Finally, $a_1a_3$ is an edge of $G - a_1a_2$, so $a_1$ and $a_3$ lie in the same face of $J'$. Thus $F^{a_3}_{a_1} = F^{a_1}_{a_3}$.
\end{proof}

\begin{lemma}[Apex-pair setup]\label{lem:apex-pair}
Let $G$ be a non-planar $5$-connected graph, and let $x,y$ be vertices of $G$ such that $G-x$ and $G-y$ are planar, and $xy \notin E(G)$ or $F^{y}_{x} = F^{x}_{y}$. Suppose that $x,y$ have two common neighbors $p,q$ with $pq \in E(G)$. Let $N'(x) := N(x) \setminus \{y\}$ and $N'(y) := N(y) \setminus \{x\}$. Then, for $J := G - \{x,y\}$, there is a facial cycle $C$ of $J$ such that
\[
   N'(x) \cup N'(y) \subseteq V(C), \qquad pq \in E(C),
\]
$C$ is induced in $J$, the sets $N'(x)$ and $N'(y)$ interleave on $C$, and, writing $X := V(J) \setminus V(C)$, the graph $J[X]$ is connected and every vertex of $C$ has a neighbor in $X$.
\end{lemma}

\begin{proof}
Since $G$ is $5$-connected, $J=G-\{x,y\}$ is $3$-connected; being a subgraph of the
planar graph $G-x$, it is planar, hence uniquely embedded by
Lemma~\ref{lem:whitney}. Fix that embedding; as $J$ is $2$-connected, every face of $J$
is bounded by a cycle.

In the planar graphs $J+x=G-y$ and $J+y=G-x$, the vertex $x$, respectively $y$, is drawn
inside the face $F_x := F^y_x$, respectively $F_y := F^x_y$, of that embedding of $J$. If $F_x\ne F_y$, then $xy \notin E(G)$ by hypothesis, and $x$ and $y$ could be drawn simultaneously, in $F_x$ and
$F_y$ respectively, giving a planar embedding of $G$, a contradiction. Hence $F_x=F_y=:F$. Set $C:=\partial F$, a cycle of $J$.

Every edge from $x$ to $J$ is drawn inside $F$ and ends on its boundary, so
$N'(x)\subseteq V(C)$; likewise $N'(y)\subseteq V(C)$. Since $J$ is $3$-connected and $C$
is a facial cycle of $J$, Lemma~\ref{lem:tutte-facial-cycles} shows that $C$ is induced
in $J$ and that $J-V(C)$ is connected whenever it is nonempty. As
$p,q\in N'(x)\cap N'(y)\subseteq V(C)$ and $pq\in E(J)$ with $C$ induced, we get
$pq\in E(C)$.

Suppose $N'(x)$ and $N'(y)$ did not interleave on $C$. Then $x$ and $y$, drawn inside $F$,
would be disjoint connected subgraphs in one region bounded by $C$ whose attachments do
not interleave, so they could be redrawn side by side inside $F$ without crossing, in a common face of $J+x+y$, and joined by the edge $xy$ if $xy \in E(G)$; this
would give a planar embedding of $G$, a contradiction. Hence $N'(x)$ and $N'(y)$ interleave
on $C$.

Finally, let $c\in V(C)$. As $C$ is induced, $c$ has exactly two neighbors on $C$; it
has at most two in $\{x,y\}$; and $\deg_G(c)\ge5$. So $c$ has a neighbor in $X$. In
particular $X\ne\emptyset$, and $J[X]=J-V(C)$ is connected.
\end{proof}

\begin{lemma}[The no-frame case]\label{lem:no-frame-k4minus}
Let $G$ be a non-planar $5$-connected graph, and let $a_1,a_2,a_3,a_4 \in V(G)$ be distinct vertices
such that
\[
    a_1a_2,\ a_1a_3,\ a_1a_4,\ a_2a_3,\ a_2a_4 \in E(G).
\]
Suppose that $G-a_3$ has a planar embedding with $a_1,a_2,a_4$ on the outer face, and that $G-a_4$ has a planar embedding with $a_1,a_2,a_3$ on the outer face. Suppose moreover that $a_3a_4 \notin E(G)$, or that $F^{a_4}_{a_3} = F^{a_3}_{a_4}$ and any four vertices of $G$ spanning at least five edges span a $K_4$.
Then $G$ contains a $K_6^-$-minor.
\end{lemma}

\begin{proof}
Suppose, for a contradiction, that $G$ has no $K_6^-$-minor.

\emph{Claim 1. Every $K_4^-$ of $G$ contains both $a_3$ and $a_4$.}

By symmetry between $a_3$ and $a_4$, it suffices to prove that $G-a_4$ is
$K_4^-$-free. Suppose $D=\{d_1,d_2,d_3,d_4\}\subseteq V(G-a_4)$ spans a $K_4^-$, with
$d_1d_2\in E(G)$ and $d_1,d_2\in N(d_3)\cap N(d_4)$. As $G$
has no $K_6^-$-minor, Lemma~\ref{K4minusgivesK6minus} gives that $G-d_3$ and $G-d_4$ are
both planar.

Apply Lemma~\ref{lem:apex-pair} to $(a_3,a_4)$ and $(a_1,a_2)$. We obtain a uniquely embedded graph
$J=G-\{a_3,a_4\}$ and an induced facial cycle $C$ bounding the face $F$ into which both
apexes are inserted. Moreover, $N'(a_3)\cup N'(a_4)\subseteq V(C)$,
$a_1a_2\in E(C)$, the sets $N'(a_3)$ and $N'(a_4)$ interleave on $C$, and, for
$X:=V(J)\setminus V(C)$, the graph $J[X]$ is connected and every vertex of $C$ has a
neighbor in $X$. Note $|V(C)|\ge4$, since $|N'(a_3)|\ge4$. Two embeddings occur below and are kept apart: that of $J$, in which $F$
contains no vertex, so all of $X$ lies in the closed disk $D^\ast$ bounded by $C$
opposite to $F$; and one of $G-d$ for the vertex $d$ chosen below. Adjacency statements
are facts about $G$, established in either and used freely.

\emph{Neither $d_3$ nor $d_4$ lies in $X$.} Suppose $d\in\{d_3,d_4\}\cap X$, so $G-d$ is
planar.

If $X=\{d\}$, then every $c\in V(C)$ has $N(c)\subseteq\{c^-,c^+,d,a_3,a_4\}$, so
$\deg_G(c)\ge5$ forces $c\in N(a_3)\cap N(a_4)$; and $|V(C)|\ge \deg_G(d) \ge 5$,
since $N(d) \subseteq V(C)$. Choose $w_1,\dots,w_5\in N(d)$ in their cyclic order on $C$. Then
\[
\{\,C[w_1,w_2),\; C[w_2,w_3)\cup a_3w_2,\; C[w_3,w_4)\cup a_4w_3,\;
   C[w_4,w_5),\; C[w_5,w_1),\; \{d\}\,\}
\]
is a $K_6^-$-minor, as in Figure~\ref{fig:no-frame-k4minus}(a), a contradiction.

So $|X|\ge2$, and $X\setminus\{d\}\ne\emptyset$ as $J[X]$ is connected; fix a component
$A$ of $X\setminus\{d\}$. Since $a_3,a_4$ have no neighbors in $X$, we get
$N_G(A)=\{d\}\cup S_A$ with $S_A:=N(A)\cap V(C)$; as $\{d\}\cup S_A$ separates $A$ from
$a_3$, $5$-connectivity gives $|S_A|\ge4$.

If $S_A=V(C)$, embed $G-d$: then $A$ lies in one region bounded by $C$ and attaches to
every vertex of $C$, so by Remark~\ref{rem:no-interleaving} no apex can lie in that
region, each having two non-consecutive neighbors on $C$. Both apexes therefore lie in
the other region, where their interleaved attachments again contradict
Remark~\ref{rem:no-interleaving}. Hence $S_A\ne V(C)$. Call an \emph{$A$-arc} a subpath $C[u_1,u_4]$ of $C$ containing $S_A$, with $u_1,u_4\in S_A$ and $V(C) \setminus V(C[u_1,u_4]) \neq \emptyset$. $A$-arcs exist: for $c\in V(C)\setminus S_A$, let $u_1$ and $u_4$ be the first vertices of $S_A$ after and before $c$ in the cyclic order of $C$; then $C[u_1,u_4]$, the subpath avoiding $c$, is an $A$-arc. Fix an $A$-arc $C[u_1,u_4]$; we choose it more specifically below, and the argument up to that point holds for any $A$-arc. Put $I:=V(C(u_1,u_4))$; as $|S_A|\ge4$, at
least two vertices of $S_A$ lie strictly inside, so $I\ne\emptyset$.

Let $W$ be the union of $A$, of $I$, and of every component of $X\setminus\{d\}$ having
a neighbor in $I$. We claim
\[
   N_G(W)\subseteq\{u_1,u_4,d,a_3,a_4\}. \tag{$\dagger$}
\]
Neither $A$ nor $I$ has a neighbor on $C$ outside $C[u_1,u_4]$, and every neighbor of
$W$ in $X$ is $d$ or lies in a component meeting $I$, hence in $W$; so it suffices that
each component $B$ of $W$ other than $A$ satisfy $N(B)\cap V(C)\subseteq V(C[u_1,u_4])$. In
the embedding of $J$, such a $B$ is connected, lies in $X\subseteq D^\ast$ and is
disjoint from $A$; a neighbor $c'\notin V(C[u_1,u_4])$ of $B$, together with a neighbor
$b\in I$ of $B$, would give attachments in the cyclic order $u_1,b,u_4,c'$,
contradicting Remark~\ref{rem:no-interleaving}. Finally $d\in X$, so $d\notin I$ and
$W\cap V(C)=I$.

Thus, the set $V(C)\setminus V(C[u_1,u_4])$ is nonempty and disjoint from
$W\cup\{u_1,u_4,d,a_3,a_4\}$, so $N_G(W)$ separates it from $W$; $5$-connectivity and
$(\dagger)$ give $N_G(W)=\{u_1,u_4,d,a_3,a_4\}$, whence $I\cap N(a_3)\ne\emptyset$ and
$I\cap N(a_4)\ne\emptyset$.

Embed $G-d$. As $N'(a_3)\cup N'(a_4)\subseteq V(C)$ and $d\notin V(C)$, the interleaving
persists, so by Remark~\ref{rem:no-interleaving} the apexes lie on opposite sides of
$C$; this is already a contradiction if $a_3a_4 \in E(G)$. For consecutive $w_1,w_2\in N(a_4)$ let the \emph{$a_4$-sector} be the closed
region bounded by $a_4w_1$, $C[w_1,w_2]$, $a_4w_2$, and define $a_3$-sectors
symmetrically. Being connected and disjoint from $V(C)\cup\{a_3,a_4\}$, the set $A$ lies
on one side of $C$ and crosses no spoke, hence lies in a single sector, say the
$a_4$-sector on $w_1,w_2$. Its boundary meets $C$ only in $C[w_1,w_2]$, so
$S_A\subseteq C[w_1,w_2]$. Choosing $C[u_1,u_4]$ as the shortest subpath of $C[w_1,w_2]$ containing $S_A$, which is an $A$-arc (as $a_4$ has neighbors in $C(w_2,w_1)$, it is not all of $C$), we get $I\subseteq C(w_1,w_2)$, which
contains no neighbor of $a_4$ --- contradicting $I\cap N(a_4)\ne\emptyset$. The
$a_3$-case is symmetric. Hence $G-d$ is non-planar, a contradiction.

Therefore, $d_3,d_4\notin X$. Since $d_3,d_4\ne a_4$, each lies in
$\{a_3\}\cup V(C)$. If $d_3=a_3$, then $d_1,d_2\in N'(a_3)\subseteq V(C)$, so
$\{d_1,d_2,d_4\}\subseteq V(C)$ spans a triangle, contradicting that $C$ is
induced and has length at least four. The same argument applies if $d_4=a_3$.
Hence, $d_3,d_4\in V(C)$.
This holds for every choice of the pair $d_3,d_4$ in a $K_4^-$; so $D$ does not span a $K_4$, since otherwise its four vertices would lie in $V(C)$ and span a triangle in $C$. Hence $d_3d_4 \notin E(G)$, and $D$ is a set of four vertices spanning five edges but not a $K_4$.

If $a_3a_4 \in E(G)$, then, by the hypothesis of the lemma, any four vertices of $G$ spanning at least five edges span a $K_4$; this contradicts the existence of $D$ and proves Claim~1 in this case. So we may assume that $a_3a_4 \notin E(G)$; in particular $N'(a_k) = N(a_k)$ for $k \in \{3,4\}$. As $a_1a_2 \in E(G)$, we may also assume that $d_3 \notin \{a_1,a_2\}$.

We use the following fact.
\begin{itemize}
    \item[(K)] \emph{Let $c \in V(C) \setminus \{a_1,a_2\}$ be such that $G - c$ is planar. Then $C - \{c,a_1,a_2\}$ has exactly two components $\Lambda_3(c)$ and $\Lambda_4(c)$, and, for $k \in \{3,4\}$, all neighbors of $a_k$ in $C - \{c,a_1,a_2\}$ lie in $\Lambda_k(c)$, which contains at least two of them.}
\end{itemize}
Indeed, if a component $\Lambda$ of $C - \{c,a_1,a_2\}$ contained neighbors of both $a_3$ and $a_4$, then $\{a_3\}$, $\{a_4\}$, $X$ and $\{a_1\}$, $\{a_2\}$, $\Lambda$ would be the branch sets of a $K_{3,3}$-minor in $G - c$, as in Figure~\ref{fig:no-frame-k4minus}(b). As $N(a_k) \subseteq V(C)$ and $|N(a_k)| \ge 5$, each $a_k$ has at least two neighbors in $C - \{c,a_1,a_2\}$, and (K) follows.
By Lemma~\ref{K4minusgivesK6minus}, $G - d_3$ and $G - d_4$ are planar, so (K) applies to $d_3$, and to $d_4$ if $d_4 \notin \{a_1,a_2\}$. In particular, $d_3$ is adjacent on $C$ to neither $a_1$ nor $a_2$.

Suppose that $a_3 \in \{d_1,d_2\}$, say $d_1 = a_3$. Then $d_2 \in N(a_3) \subseteq V(C)$, and as $C$ is induced, $d_3$ and $d_4$ are the two neighbors of $d_2$ on $C$; thus $d_2 \notin \{a_1,a_2\}$ and $d_2 \in V(\Lambda_3(d_3))$. If $d_4 \in \{a_1,a_2\}$, then $\Lambda_3(d_3) = \{d_2\}$ contains only one neighbor of $a_3$, contrary to (K). Otherwise, the component of $C - \{d_4,a_1,a_2\}$ containing $d_2$ also contains $d_3$ and $\Lambda_4(d_3)$, hence neighbors of both $a_3$ and $a_4$, contrary to (K) applied to $d_4$.

Hence $d_1,d_2 \in V(C) \cup X$, and they do not both lie in $V(C)$: otherwise $d_1$ has the three neighbors $d_2,d_3,d_4$ on $C$, contrary to $C$ being induced.
Let $\Lambda$ be the component of $C - \{d_3,d_4\}$ containing neither $a_1$ nor $a_2$; it exists since $d_3d_4 \notin E(G)$ and $a_1a_2 \in E(C)$. We show that $\Lambda$ contains neighbors of $a_3$ and of $a_4$. Then
\[
   \{a_3\},\quad\{a_4\},\quad X \qquad\text{and}\qquad
   \{a_1\},\quad\{a_2\},\quad \Lambda
\]
are the branch sets of a $K_{3,3}$-minor in $G-d_3$, contradicting the
planarity of $G-d_3$; see Figure~\ref{fig:no-frame-k4minus}(b).

First suppose that $d_1 \in V(C)$, say. Then $d_3d_1d_4$ is a subpath of $C$ and $d_2 \in X$. In the embedding of $J$, the triangles $d_1d_2d_3$ and $d_1d_2d_4$ bound faces: the vertices on their side not containing $F$ lie in $X$, so they have no neighbors outside the triangle, and there are none by $5$-connectivity. Hence the neighbors of $d_1$ in $J$ are $d_2$, $d_3$ and $d_4$, and as $\deg_G(d_1) \ge 5$, $d_1 \in N(a_3) \cap N(a_4)$. Moreover $d_1 \notin \{a_1,a_2\}$, since otherwise $d_3$ or $d_4$, a neighbor of $d_1$ on $C$, would be the other vertex of $\{a_1,a_2\}$, contrary to $d_3 \notin \{a_1,a_2\}$ and to the fact that $d_3$ is not adjacent on $C$ to $a_1, a_2$. Hence $\Lambda = \{d_1\}$, as required.

Now suppose that $d_1,d_2 \in X$. The internally disjoint paths $d_3d_1d_4$ and $d_3d_2d_4$ lie in the disk $D^\ast$; let $d$ be the one of $d_1,d_2$ such that the region of $D^\ast$ bounded by $d_3 d d_4$ and $C[d_3,d_4] \supseteq \Lambda$ avoids the other path, and let $W$ consist of the vertices of $\Lambda$ and of the vertices of $X$ in this region other than $d$. Then $W \neq \emptyset$, $N_J(W) \subseteq \{d,d_3,d_4\}$, and only vertices of $\Lambda$ can be adjacent to $a_3$ or $a_4$. As $G$ is $5$-connected, $\Lambda$ contains neighbors of both $a_3$ and $a_4$. This proves the claim.

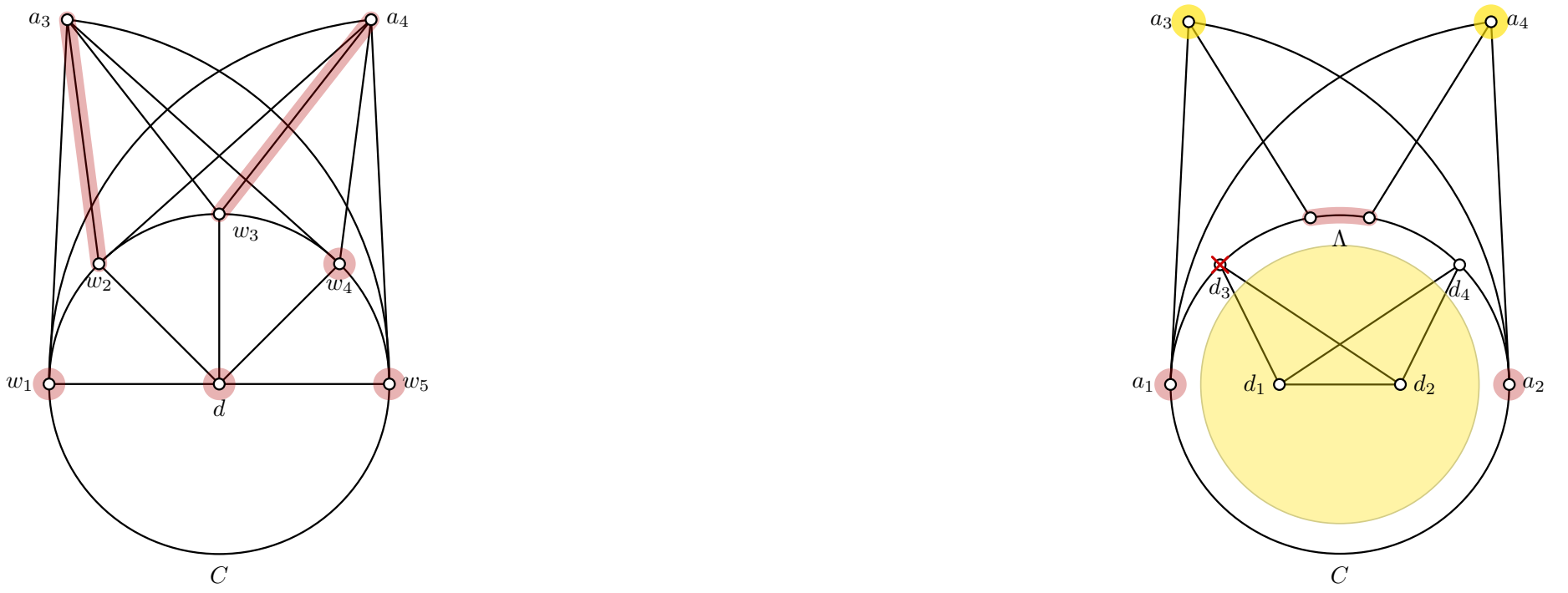
\begin{figure}[htbp]
    \centering
        \begin{minipage}{0.49\textwidth}
            \centering
            \adjustbox{max width=0.8\linewidth,max height=4.6cm}{%
    \begingroup\setlength{\textwidth}{\linewidth}\begin{tikzpicture}[scale=0.9, line cap=round, line join=round]
    \def\rin{2.8}
    \def\rout{4.2}

    \coordinate (O) at (0,0);
    \coordinate (u) at (-2.5,6);
    \coordinate (v) at (2.5,6);

    \coordinate (w1) at (180:\rin);
    \coordinate (w2) at (135:\rin);
    \coordinate (w3) at (90:\rin);
    \coordinate (w4) at (45:\rin);
    \coordinate (w5) at (0:\rin);

    \tikzset{
        setAHL/.style={draw=red!70!black, opacity=0.30, line width=7pt, line cap=round, line join=round},
    }

    \draw[thick] (O) circle (\rin);

    \draw[thick] (u) -- (w1);
    \draw[thick] (u) -- (w2);
    \draw[thick] (u) -- (w3);
    \draw[thick] (u) -- (w4);
    \draw[thick] (u) to [bend left = 42] (w5);

    \draw[thick] (v) to [bend right = 42] (w1);
    \draw[thick] (v) -- (w2);
    \draw[thick] (v) -- (w3);
    \draw[thick] (v) -- (w4);
    \draw[thick] (v) -- (w5);

    \draw[thick] (O) -- (w1);
    \draw[thick] (O) -- (w2);
    \draw[thick] (O) -- (w3);
    \draw[thick] (O) -- (w4);
    \draw[thick] (O) -- (w5);

    \draw[setAHL] (u) -- (w2);
    \draw[setAHL] (v) -- (w3);
    \draw[setAHL] (O) circle (0.13);
    \draw[setAHL] (w1) circle (0.13);
    \draw[setAHL] (w4) circle (0.13);
    \draw[setAHL] (w5) circle (0.13);

    \foreach \p in {O,u,v,w1,w2,w3,w4,w5} {
        \fill[white] (\p) circle (0.09);
        \draw[thick] (\p) circle (0.09);
    }

    \node[left=3pt] at (u) {$a_3$};
    \node[right=3pt] at (v) {$a_4$};
    \node[below=3pt] at (O) {$d$};

    \node[left=3pt] at (w1) {$w_1$};
    \node[below=2pt] at (w2) {$w_2$};
    \node[below right=2pt] at (w3) {$w_3$};
    \node[below=3pt] at (w4) {$w_4$};
    \node[right=2pt] at (w5) {$w_5$};

    \node[fill=white, inner sep=1pt] at (0,-3.15) {$C$};
\end{tikzpicture}\endgroup}\par
            \textbf{(a)} $K_6^-$-minor model
        \end{minipage}\hfill
        \begin{minipage}{0.49\textwidth}
            \centering
            \adjustbox{max width=0.8\linewidth,max height=4.6cm}{%
    \begingroup\setlength{\textwidth}{\linewidth}\begin{tikzpicture}[scale=0.9, line cap=round, line join=round]
    \def\rin{2.8}
    \def\rout{4.2}

    \coordinate (O) at (0,0);
    \coordinate (a) at (-2.5,6);
    \coordinate (b) at (2.5,6);

    \coordinate (c) at (-1,0);
    \coordinate (d) at (1,0);

    \coordinate (w1) at (180:\rin);
    \coordinate (w2) at (135:\rin);
    \coordinate (u) at (80:\rin);
    \coordinate (v) at (100:\rin);
    \coordinate (w4) at (45:\rin);
    \coordinate (w5) at (0:\rin);

    \tikzset{
        setAHL/.style={draw=red!70!black, opacity=0.30, line width=7pt, line cap=round, line join=round},
        setBHL/.style={draw=yellow!85!orange, opacity=0.65, line width=8pt, line cap=round, line join=round},
        setBDisk/.style={fill=yellow!85!orange, fill opacity=0.35, draw=yellow!70!black, draw opacity=0.70, line width=0.6pt},
    }

    \draw[thick] (O) circle (\rin);

    \draw[thick] (a) -- (w1);
    \draw[thick] (a) -- (v);
    \draw[thick] (a) to [bend left=40] (w5);

    \draw[thick] (b) to [bend right = 40] (w1);
    \draw[thick] (b) -- (u);
    \draw[thick] (b) -- (w5);

    \draw[thick] (c) -- (d);
    \draw[thick] (c) -- (w2);
    \draw[thick] (c) -- (w4);

    \draw[thick] (d) -- (w2);
    \draw[thick] (d) -- (w4);

    \path[setBDisk] (O) circle (2.30);
    \draw[setAHL] (w1) circle (0.13);
    \draw[setBHL] (a) circle (0.13);
    \draw[setBHL] (b) circle (0.13);
    \draw[setAHL] (w5) circle (0.13);
    \draw[setAHL] (u) arc[start angle=80, end angle=100, radius=\rin];

    \foreach \p in {a,b,c,d,u,v,w1,w2,w4,w5} {
        \fill[white] (\p) circle (0.09);
        \draw[thick] (\p) circle (0.09);
    }

    \draw[red!80!black, line width=1.1pt] ($(w2)+(-0.13,-0.13)$) -- ($(w2)+(0.13,0.13)$);
    \draw[red!80!black, line width=1.1pt] ($(w2)+(-0.13,0.13)$) -- ($(w2)+(0.13,-0.13)$);

    \node[left=3pt] at (a) {$a_3$};
    \node[right=3pt] at (b) {$a_4$};

    \node[left=3pt] at (w1) {$a_1$};
    \node[below=2pt] at (w2) {$d_3$};
    \node[below =3pt] at (90:\rin) {$\Lambda$};

    \node[below=3pt] at (w4) {$d_4$};
    \node[right=2pt] at (w5) {$a_2$};

    \node[left=2pt] at (c) {$d_1$};
    \node[right=2pt] at (d) {$d_2$};

    \node[fill=white, inner sep=1pt] at (0,-3.15) {$C$};
\end{tikzpicture}\endgroup}\par
            \textbf{(b)} $K_{3,3}$-minor model
        \end{minipage}
        \caption{The two minor models used in the no-frame $K_4^-$ case.}
        \label{fig:no-frame-k4minus}
    \vspace{0.25em}
\end{figure}

\emph{Claim 2. The planar degree count gives a contradiction.}
Put $G_3=G-a_3$ and $G_4=G-a_4$. By Claim~1, both planar graphs $G_3$ and
$G_4$ are $K_4^-$-free, and they are $4$-connected since $G$ is $5$-connected.

We first record the planar estimate used twice below. Let $J$ be a $2$-connected $K_4^-$-free plane graph with
$N \ge 4$ vertices and $M$ edges, so that every face of $J$ is bounded by a cycle. If $f$ is the number of faces and $f_3$ the number of triangular faces,
then every non-triangular face has degree at least $4$, and every edge $xy$ is incident with at most one triangular face. Indeed, two triangular faces containing $xy$ cannot be bounded by the same triangle, since then $J$ would be this triangle, contrary to $N\ge4$; so they are bounded by distinct triangles $xyz$ and $xyw$, and $\{x,y,z,w\}$ spans a $K_4^-$. (For $J=K_3$ the estimate below fails.) Hence $3f_3\le M$. Summing face-degrees gives $2M\ge 3f_3+4(f-f_3)=4f-f_3$, so
\[
    f\le \frac{2M+f_3}{4}\le \frac{2M+M/3}{4}=\frac{7M}{12}.
\]
Euler's formula $2=N-M+f$ then yields $2\le N-\tfrac{5M}{12}$, i.e.,
\begin{equation}\label{eq:k4free-planar-bound}
    M\le \frac{12}{5}(N-2).
\end{equation}

Let
\[
    n_3=|N_G(a_3)\setminus\{a_1,a_2,a_4\}|,\qquad n_4=|N_G(a_4)\setminus\{a_1,a_2,a_3\}|,
\]
and $n=|V(G)|$. Since $G_3$ is $2$-connected and has $n-1\ge5$ vertices, \eqref{eq:k4free-planar-bound} gives
$2|E(G_3)|\le \tfrac{24}{5}(n-3)$. On the other hand $G$ is $5$-connected, so every vertex has
degree $\ge5$ in $G$. In $G_3$: $a_1,a_2$ have degree $\ge4$; $a_4$ has degree $n_4+2$; each of the $n_3$ vertices of $N_G(a_3)\setminus\{a_1,a_2,a_4\}$ has degree
$\ge4$; and every remaining vertex has degree $\ge5$. Thus
\[
    2|E(G_3)|\ge 4+4+(n_4+2)+4n_3+5(n-4-n_3)=5n-n_3+n_4-10.
\]
Combining,
\[
    n_3-n_4\ge \frac{n+22}{5}.
\]
The same argument applied to $G_4$ gives $2|E(G_4)|\le\tfrac{24}{5}(n-3)$ and
\[
    2|E(G_4)|\ge 4+4+(n_3+2)+4n_4+5(n-4-n_4)=5n+n_3-n_4-10,
\]
hence $n_4-n_3\ge\frac{n+22}{5}$. Adding the two inequalities gives $0\ge 2(n+22)/5$, a contradiction.
Therefore, $G$ contains a $K_6^-$-minor.
\end{proof}

\subsection{The $K_4^-$-free induction}

\begin{lemma}[Inner cycle]\label{lem:inner-cycle}
    Let $G$ be $5$-connected and $K_4^-$-free, let $uv\in E(G)$ with $G'=G/uv$
    planar, let $x$ be the contracted vertex, and fix a planar embedding of $G'$ in which $x$ is not incident with the outer face; such an embedding exists since $G'-x$ is $3$-connected and therefore has a face not containing the point where $x$ is drawn.
    Put $H=G'-x$. Then $H$ has a facial cycle $C_1$ such that:
    \begin{enumerate}
        \item[(i)] $N_{G'}(x)\subseteq V(C_1)$, and $x$ is the only vertex of $G'$ in the
            open disk $\Delta_1$ bounded by $C_1$;
        \item[(ii)] with $R:=V(G)\setminus(V(C_1)\cup\{u,v\})$, every vertex of $C_1$ has
            a neighbor in $R$;
        \item[(iii)] $G[R]$ is connected.
    \end{enumerate}
\end{lemma}

\begin{proof}
    Contracting an edge lowers connectivity by at most one, so $G'$ is
    $4$-connected and $H=G'-x$ is $3$-connected. Let $F_1$ be the face of $H$
    containing the point where $x$ was drawn, and $C_1$ its boundary; as $H$ is
    $2$-connected, $C_1$ is a cycle. Every edge at $x$ runs into $F_1$ and ends on
    its boundary, so $N_{G'}(x)\subseteq V(C_1)$, and $x$ is the only vertex of
    $G'$ in the open disk $\Delta_1\subseteq F_1$ bounded by $C_1$. This is (i).

    By Lemma~\ref{lem:tutte-facial-cycles}, $C_1$ is induced in $H$ (equivalently in $G$,
    since contracting $uv$ creates no edge between two vertices of $C_1$, only
    incidences at $x\notin V(C_1)$).

    From $N_{G'}(x)=\bigl(N_G(u)\cup N_G(v)\bigr)\setminus\{u,v\}\subseteq V(C_1)$
    we get $N_G(u)\cup N_G(v)\subseteq V(C_1)\cup\{u,v\}$, so $u,v$ have no neighbor
    in $R$; thus $G[R]=H[R]=H-V(C_1)$. For $w\in V(C_1)$, since $C_1$ is induced,
    the only neighbors of $w$ on $C_1$ are its two cycle-neighbors, and apart from
    these and possibly $u,v$ every neighbor of $w$ lies in $R$; as
    $\deg_G(w)\ge 5$, at least one of them lies in $R$. This proves (ii), and in
    particular $R\neq\emptyset$.

    Finally, $G[R]=H-V(C_1)$ is connected by the non-separating property and nonempty
    by (ii), giving (iii).
\end{proof}

\begin{lemma}\label{planaritypreserved}
    Let $G$ be a non-planar $5$-connected graph containing no $K_4^-$.
    If there exists an edge $uv \in E(G)$ such that $G/uv$ is planar, then $G$ contains a $K_6^-$-minor.
\end{lemma}

\begin{proof}
    Suppose that there is an edge $uv\in E(G)$ with $G/uv$ planar; fix a planar
    embedding as in Lemma~\ref{lem:inner-cycle} and let $x$ be the contracted vertex. By Lemma~\ref{lem:inner-cycle}
    there is a facial cycle $C_1$ of $G/uv-x$ with $N_{G/uv}(x)\subseteq V(C_1)$
    and, writing $R:=V(G)\setminus(V(C_1)\cup\{u,v\})$, every vertex of $C_1$ has a
    neighbor in $R$ and $G[R]$ is connected. Thus $G[R]$ is a connected branch set
    adjacent to every vertex of $C_1$; it will be one branch set of each model
    below, and it remains to exhibit five further pairwise-disjoint connected sets
    in $\{u,v\}\cup V(C_1)$, pairwise adjacent except for one pair.

    Since $G$ is non-planar, the neighbors of $u$ and of $v$ interleave on $C_1$;
    otherwise the embedding of $G/uv$ could be modified into a planar embedding of
    $G$. Moreover $N_G(u),N_G(v)\subseteq V(C_1)\cup\{u,v\}$ and
    $\deg_G(u),\deg_G(v)\ge 5$, so each of $u,v$ has at least four neighbors on
    $C_1$; and since $G$ is $K_4^-$-free, two common neighbors of $u,v$ would form
    a $K_4^-$ with $u,v$, so $u$ and $v$ share at most one neighbor.
    Hence, we may pick three neighbors of $u$ and three of $v$ on $C_1$, all six
    distinct, with the two triples interleaving. Write $A_u = (N_G(u) \cap V(C_1)) \setminus N_G(v)$ and $A_v = (N_G(v) \cap V(C_1)) \setminus N_G(u)$. If $A_u$ and $A_v$ contain $a,b,a',b'$ in this cyclic order, we take them and one more vertex of each set; otherwise these two sets occupy two disjoint subpaths of $C_1$, and, as the neighbors of $u$ and of $v$ interleave on $C_1$, their common neighbor $z$ lies strictly between two vertices of one of these subpaths, say of $A_u$; we then take $z$ and two vertices of $A_v$, and three vertices of $A_u$ including one on each side of $z$. Relabel them $w_1,\dots,w_6$ in
    cyclic order along $C_1$.

    Assign to each $w_i$ an \emph{owner} in $\{u,v\}$: $u$ if $w_i$ was chosen as one of the three neighbors of $u$, and $v$ otherwise (in the second case, $z$ is assigned to $v$). Each $w_i$ is adjacent to its owner, and the constructions below only use these edges; a further edge $zu$ only adds adjacencies between branch sets. The owner-sequence has three
    $u$'s and three $v$'s that interleave, i.e., the $u$'s are not three
    consecutive positions. Up to rotation, reflection, and exchanging $u$ and $v$,
    there are exactly two such patterns: $(u,v,u,v,u,v)$ and $(u,u,v,u,v,v)$.

    For $w_i,w_j$ in cyclic order let $C_1[w_i,w_j)$ denote the arc from $w_i$
    (inclusive) to $w_j$ (exclusive) in the direction $w_1\to w_2\to\cdots\to
    w_6\to w_1$. Since $C_1$ is induced, two of these arcs are adjacent only if they
    are consecutive.

    \smallskip
    \noindent\textbf{Case 1: owner-sequence $(u,v,u,v,u,v)$}, i.e.
    $uw_1,uw_3,uw_5,vw_2,vw_4,vw_6\in E(G)$. Then
    \[
        \{\,\{u\}\cup C_1[w_1,w_2),\;\{v\}\cup C_1[w_2,w_3),\;C_1[w_3,w_4),\;
        C_1[w_4,w_6),\;C_1[w_6,w_1),\;G[R]\,\}
    \]
    is a $K_6^-$-minor; see Figure~\ref{fig:planar-contraction-models}(a).

    \smallskip
    \noindent\textbf{Case 2: owner-sequence $(u,u,v,u,v,v)$}, i.e.
    $uw_1,uw_2,uw_4,vw_3,vw_5,vw_6\in E(G)$. Then
    \[
        \{\,\{u\}\cup C_1[w_4,w_5),\;\{v\}\cup C_1[w_3,w_4),\;G[R],\;
        C_1[w_6,w_2),\;C_1[w_2,w_3),\;C_1[w_5,w_6)\,\}
    \]
    is a $K_6^-$-minor; see Figure~\ref{fig:planar-contraction-models}(b).

    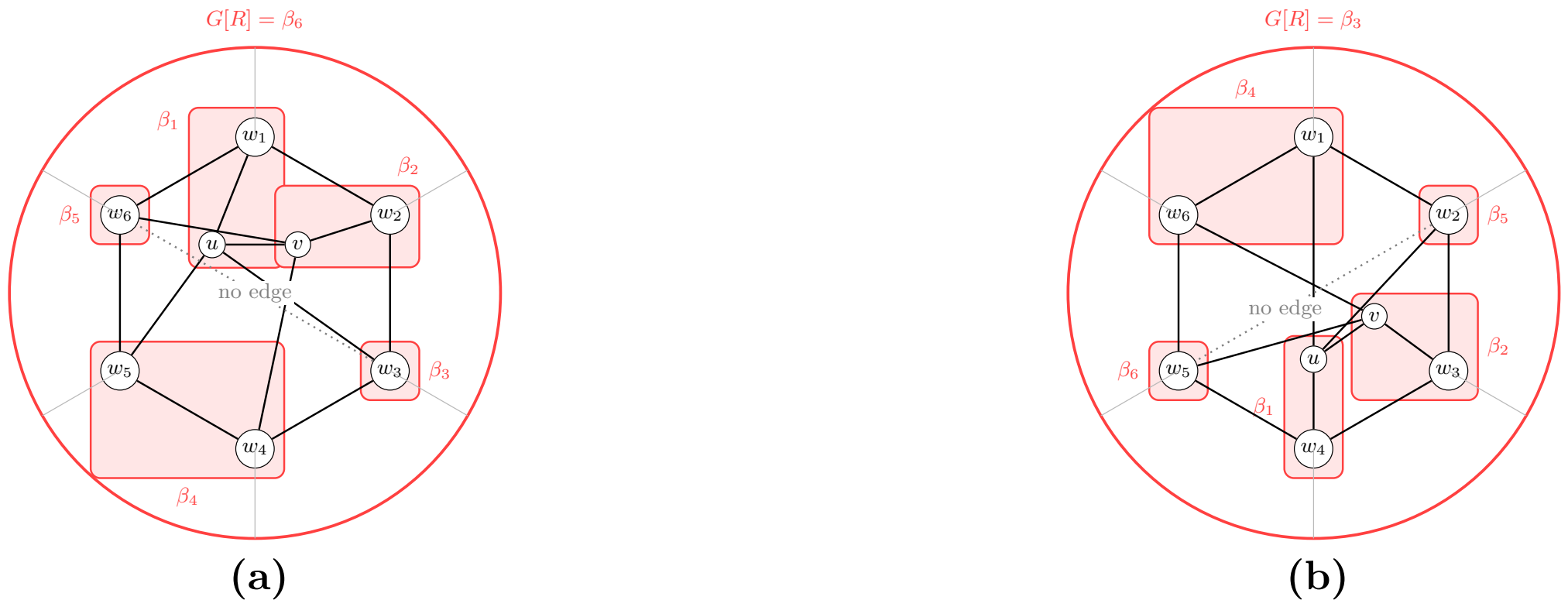
\begin{figure}[htbp]
        \centering
        \begin{minipage}{0.48\textwidth}
            \centering
            \adjustbox{max width=0.8\linewidth,max height=4.6cm}{%
    \begingroup\setlength{\textwidth}{\linewidth}\begin{tikzpicture}[scale=1.15,
  wn/.style={circle,draw,fill=white,inner sep=1.6pt,font=\small},
  ed/.style={thick},
  redge/.style={gray!55},
  bs/.style={rounded corners,fill=red!10,draw=red!75,thick,inner sep=4pt}]
  \node[wn](w1) at (90:2){$w_1$};  \node[wn](w2) at (30:2){$w_2$};
  \node[wn](w3) at (-30:2){$w_3$}; \node[wn](w4) at (-90:2){$w_4$};
  \node[wn](w5) at (-150:2){$w_5$};\node[wn](w6) at (150:2){$w_6$};
  \node[wn](u) at (-0.55,0.62){$u$}; \node[wn](v) at (0.55,0.62){$v$};
  \draw[red!75,very thick] (0,0) circle (3.15);
  \node[red!75,font=\small] at (90:3.5){$G[R]=\beta_6$};
  \foreach \a in {90,30,-30,-90,-150,150} \draw[redge] (\a:2)--(\a:3.15);
  \draw[ed](w1)--(w2)--(w3)--(w4)--(w5)--(w6)--(w1);
  \draw[ed](u)--(w1); \draw[ed](u)--(w3); \draw[ed](u)--(w5); \draw[ed](u)--(v);
  \draw[ed](v)--(w2); \draw[ed](v)--(w4); \draw[ed](v)--(w6);
  \draw[gray,dotted,thick] (w3)--(w6)
       node[midway,fill=white,font=\small,inner sep=1pt]{no edge};
  \begin{pgfonlayer}{background}
    \node[bs,fit=(u)(w1),label={[red!75,font=\small]above left:$\beta_1$}]{};
    \node[bs,fit=(v)(w2),label={[red!75,font=\small]above right:$\beta_2$}]{};
    \node[bs,fit=(w3),label={[red!75,font=\small]right:$\beta_3$}]{};
    \node[bs,fit=(w4)(w5),label={[red!75,font=\small]below:$\beta_4$}]{};
    \node[bs,fit=(w6),label={[red!75,font=\small]left:$\beta_5$}]{};
  \end{pgfonlayer}
\end{tikzpicture}\endgroup}\par
            \textbf{(a)}
        \end{minipage}\hfill
        \begin{minipage}{0.48\textwidth}
            \centering
            \adjustbox{max width=0.8\linewidth,max height=4.6cm}{%
    \begingroup\setlength{\textwidth}{\linewidth}\begin{tikzpicture}[scale=1.15,
  wn/.style={circle,draw,fill=white,inner sep=1.6pt,font=\small},
  ed/.style={thick},
  redge/.style={gray!55},
  bs/.style={rounded corners,fill=red!10,draw=red!75,thick,inner sep=4pt}]
  \node[wn](w1) at (90:2){$w_1$};  \node[wn](w2) at (30:2){$w_2$};
  \node[wn](w3) at (-30:2){$w_3$}; \node[wn](w4) at (-90:2){$w_4$};
  \node[wn](w5) at (-150:2){$w_5$};\node[wn](w6) at (150:2){$w_6$};
  \node[wn](u) at (0,-0.85){$u$}; \node[wn](v) at (0.78,-0.30){$v$};
  \draw[red!75,very thick] (0,0) circle (3.15);
  \node[red!75,font=\small] at (90:3.5){$G[R]=\beta_3$};
  \foreach \a in {90,30,-30,-90,-150,150} \draw[redge] (\a:2)--(\a:3.15);
  \draw[ed](w1)--(w2)--(w3)--(w4)--(w5)--(w6)--(w1);
  \draw[ed](u)--(w1); \draw[ed](u)--(w2); \draw[ed](u)--(w4); \draw[ed](u)--(v);
  \draw[ed](v)--(w3); \draw[ed](v)--(w5); \draw[ed](v)--(w6);
  \draw[gray,dotted,thick] (w2)--(w5)
       node[pos=0.62,fill=white,font=\small,inner sep=1pt]{no edge};
  \begin{pgfonlayer}{background}
    \node[bs,fit=(u)(w4),label={[red!75,font=\small]left:$\beta_1$}]{};
    \node[bs,fit=(v)(w3),label={[red!75,font=\small]right:$\beta_2$}]{};
    \node[bs,fit=(w6)(w1),label={[red!75,font=\small]above:$\beta_4$}]{};
    \node[bs,fit=(w2),label={[red!75,font=\small]right:$\beta_5$}]{};
    \node[bs,fit=(w5),label={[red!75,font=\small]left:$\beta_6$}]{};
  \end{pgfonlayer}
\end{tikzpicture}\endgroup}\par
            \textbf{(b)}
        \end{minipage}
        \caption{The two $K_6^-$-models in Lemma~\ref{planaritypreserved}. In both panels,
        the red outer ring represents the connected branch set $G[R]$ adjacent to every
        vertex of $C_1$, and $\beta_1,\dots,\beta_6$ denote the branch sets in the order in which they are listed in the text. In (a), the owner-sequence is $(u,v,u,v,u,v)$ and
        $G[R]=\beta_6$; the only missing adjacency is between $\beta_3$ and $\beta_5$.
        In (b), the owner-sequence is $(u,u,v,u,v,v)$ and $G[R]=\beta_3$; the only
        missing adjacency is between $\beta_5$ and $\beta_6$.}
        \label{fig:planar-contraction-models}
    \end{figure}

\end{proof}

We will use the following contractible-edge theorem of Kawarabayashi to run the induction in the $K_4^-$-free case.

Kawarabayashi~\cite{Kawarabayashi2001} proved that, for every odd $k \ge 3$, every $k$-connected graph containing no $K_4^-$ has an edge whose contraction yields a $k$-connected graph. We use the case $k = 5$.

\begin{lemma}[Kawarabayashi~\cite{Kawarabayashi2001}]\label{kawarabayashi-contractible-edge}
    Let $G$ be a $5$-connected graph containing no $K_4^-$.
    Then there exists an edge $e \in E(G)$ such that $G/e$ is $5$-connected.
\end{lemma}

\begin{proof}[Proof of Theorem~\ref{thm:k6minus}]
    We argue by contradiction. Let $G$ be a counterexample with a minimum number of vertices; thus $G$ is $5$-connected and non-planar, and contains no $K_6^-$-minor.

    First suppose that $G$ contains a $K_4^-$, that is, four vertices $a_1,a_2,a_3,a_4$ spanning at least five edges, labeled so that $a_1a_2,a_1a_3,a_1a_4,a_2a_3,a_2a_4 \in E(G)$. By Lemma~\ref{K4minusgivesK6minus}, which applies to every such quadruple and labelling, $G - a_3$ has a planar embedding with $a_1,a_2,a_4$ on the outer face, and $G - a_4$ has one with $a_1,a_2,a_3$ on the outer face. If $a_3a_4 \notin E(G)$ for some such quadruple, then Lemma~\ref{lem:no-frame-k4minus} gives a $K_6^-$-minor. Otherwise, any four vertices of $G$ spanning at least five edges span a $K_4$. In particular, $\{a_1,a_2,a_3,a_4\}$ induces a $K_4$, each of its vertices can play the role of $a_3$, and so $G - a_i$ is planar for every $i$. By Lemma~\ref{lem:k4-pair}, after relabelling, $F^{a_4}_{a_3} = F^{a_3}_{a_4}$, and Lemma~\ref{lem:no-frame-k4minus} gives a $K_6^-$-minor. Hence $G$ is $K_4^-$-free.

    Since $G$ is $5$-connected and $K_4^-$-free, Lemma~\ref{kawarabayashi-contractible-edge} gives an edge $e$ such that $G' = G/e$ is $5$-connected. By Lemma~\ref{planaritypreserved}, $G'$ is non-planar. As $|V(G')| < |V(G)|$, the minimality of $G$ implies that $G'$ has a $K_6^-$-minor, which is also a minor of $G$, a contradiction.
\end{proof}

\section*{Concluding remarks}

The frame structure theorem shows that every frame rooted at four prescribed vertices of a $5$-connected graph forces one of five explicit configurations, independently of the edges among the roots and of any assumption on minors. When the roots span a $K_4^-$, each configuration yields a $K_6^-$-minor; with the planar argument for the case without a frame and the induction based on Kawarabayashi's contractible-edge theorem, this proves that every non-planar $5$-connected graph contains a $K_6^-$-minor.

We now explain how the frame structure theorem fits into an approach to J{\o}rgensen's conjecture. Let $G$ be a $6$-connected graph with no $K_6$-minor, and let $a_1,a_2,a_3,a_4$ span a $K_4$. By the rooted reformulation of the introduction, Conjecture~\ref{conj:jorgensen} for $G$ amounts to showing that $G-a_i$ has a planar embedding with the other three roots on a common face for some $i$. If this fails, it fails at $a_3$ or $a_4$, so $G$ contains a frame and hence one of the five configurations of Theorem~\ref{structure}; a configuration of type $A$, together with the $K_4$ on the roots and a path joining its two trees while avoiding the roots, yields a $K_6$-minor in $G$. One way forward is to show that, under $6$-connectivity and a $K_4$ on the roots, each configuration of types $B$--$E$ yields a $K_6$-minor or forces rooted planarity of some $G-a_i$. The types $A$--$E$, however, were chosen so that each yields a $K_6^-$-minor in the presence of a $K_4^-$; they are tailored to Theorem~\ref{thm:k6minus}, not to $K_6$. The proof of Theorem~\ref{structure} is a case analysis of the paths attached to a minimal frame, and the same analysis can be run toward a different list of terminal configurations, chosen so that each yields a $K_6$-minor in the presence of a $K_4$ or certifies rooted planarity. We believe that this is the natural route from frames to Conjecture~\ref{conj:jorgensen}. One caveat makes this question precise. If the roots span a $K_4$, then every frame together with the edge $a_3a_4$ already contains a subgraph of type $C$, obtained by taking $a'_2=a_2$, $u=a_3$, $v=a_4$ and $U$ the edge $a_3a_4$; so Theorem~\ref{structure} must be applied to $G-a_3a_4$, which is still $5$-connected and still contains the frame, and the configuration it produces must then be analyzed with the edge $a_3a_4$ restored. Graphs with no $K_4$ require a separate argument, as in the $K_4^-$-free case treated here.

\section*{Declaration of generative AI and AI-assisted technologies in the writing process}

During the preparation of this paper, the authors used large language model assistants (Anthropic's Claude and OpenAI's ChatGPT) to improve the wording and organization of the text, to edit the \LaTeX{} source and figures, and to look for gaps and inconsistencies in draft proofs. They were also used to write short scripts checking selected path configurations in Lemmas~\ref{S3toT}--\ref{S1pairs} and the branch-set models of Lemma~\ref{K4minusgivesK6minus} on small instances; these scripts served only as an additional sanity check, and no result of the paper relies on them. After using these tools, the authors reviewed and edited the content as needed. All mathematical results, proofs, and constructions were developed and verified by the authors, who take full responsibility for the content of the publication.


\begin{thebibliography}{99}

\bibitem{Diestel}
R. Diestel,
\emph{Graph Theory}, 5th ed., Graduate Texts in Mathematics 173, Springer, Berlin, 2017.

\bibitem{Hadwiger1943}
H. Hadwiger,
\emph{{\"U}ber eine Klassifikation der Streckenkomplexe},
Vierteljschr. Naturforsch. Ges. Z{\"u}rich \textbf{88} (1943), 133--142.

\bibitem{HKY2025}
K. Hayashi, K. Kawarabayashi, and Y. Yoo,
\emph{Chasing tripods to obtain a rooted subdivision},
SIAM J. Discrete Math. \textbf{39} (2025), 1683--1711.

\bibitem{HeWangYu2020}
D. He, Y. Wang, and X. Yu,
\emph{The Kelmans--Seymour conjecture IV: A proof},
J. Combin. Theory Ser. B \textbf{144} (2020), 309--358.

\bibitem{Jorgensen1994}
L. K. J{\o}rgensen,
\emph{Contractions to $K_8$},
J. Graph Theory \textbf{18} (1994), 431--448.

\bibitem{Kawarabayashi2001}
K. Kawarabayashi,
\emph{Note on $k$-contractible edges in $k$-connected graphs},
Australas. J. Combin. \textbf{24} (2001), 165--168.

\bibitem{KNTW2018}
K. Kawarabayashi, S. Norine, R. Thomas, and P. Wollan,
\emph{$K_6$ minors in large $6$-connected graphs},
J. Combin. Theory Ser. B \textbf{129} (2018), 158--203.

\bibitem{Lo2026}
O.-H.~S. Lo,
\emph{A characterization of graphs with no $K_{3,4}$ minor},
arXiv:2603.27973, 2026.

\bibitem{Perfect1968}
H. Perfect,
\emph{Applications of Menger's graph theorem},
J. Math. Anal. Appl. \textbf{22} (1968), 96--111.

\bibitem{RST1993}
N. Robertson, P. D. Seymour, and R. Thomas,
\emph{Hadwiger's conjecture for $K_6$-free graphs},
Combinatorica \textbf{13} (1993), 279--361.

\bibitem{Tutte1963}
W. T. Tutte,
\emph{How to draw a graph},
Proc. London Math. Soc. (3) \textbf{13} (1963), 743--767.

\bibitem{Whitney1933}
H. Whitney,
\emph{2-isomorphic graphs},
Amer. J. Math. \textbf{55} (1933), 245--254.

\end{thebibliography}
\end{document}